\documentclass[12pt]{amsart}

\usepackage[utf8]{inputenc}
\usepackage{graphicx}
\usepackage{graphicx,color}
\usepackage{latexsym}
\usepackage[all]{xy}
\usepackage{mathrsfs}
\usepackage{enumerate}
\usepackage[shortlabels]{enumitem}
\usepackage{multicol}
\usepackage{lipsum}
\usepackage{amssymb}
\usepackage{tikz}
\usetikzlibrary{cd}
\usetikzlibrary{patterns}
\definecolor{DimGray}{rgb}{0.41, 0.41, 0.41}
\usepackage{float}
\usepackage{bm}
\usepackage{mathptmx}
\usetikzlibrary{shapes.geometric,positioning}
\usepackage{dsfont}
\usepackage{cjhebrew}
\usepackage[cal=dutchcal]{mathalfa}
\usepackage[alpine]{ifsym}
\usepackage{xpicture}
\usepackage{calculator}
\usepackage{graphicx}
\usepackage{pgf,tikz,pgfplots}
\pgfplotsset{compat=1.15}
\usetikzlibrary{arrows}

\usepackage{amsmath,mathtools}
\usepackage{amsfonts}
\usepackage{amssymb,enumerate}
\usepackage{amsthm}
\usepackage{ mathrsfs }
\usepackage{fancyhdr}
\usepackage{hyperref}
\usepackage{cleveref}

\usepackage{caption}
\usepackage{subcaption}

\def\sideremark#1{\ifvmode\leavevmode\fi\vadjust{\vbox to0pt{\vss 
			\hbox to 0pt{\hskip\hsize\hskip1em           
				\vbox{\hsize2cm\tiny\raggedright\pretolerance10000
					\noindent #1\hfill}\hss}\vbox to8pt{\vfil}\vss}}} %
\theoremstyle{plain}
\newtheorem{theorem}{Theorem}[section]
\newtheorem{proposition}[theorem]{Proposition}

\newtheorem*{theorem*}{Theorem}
\newtheorem{lemma}[theorem]{Lemma}
\newtheorem{corollary}[theorem]{Corollary}

\theoremstyle{definition}
\newtheorem{defin}[theorem]{Definition}

\newtheorem{ex}[theorem]{Example}

\theoremstyle{remark}
\newtheorem{rem}[theorem]{Remark}
\newcommand{\ord}{\operatorname{ord}}

\def\lcm{\operatorname{lcm}}

\newcommand{\Pds}{\mathscr{P}}
\newcommand{\Lds}{\mathscr{L}}
\newcommand{\longitude}{\mathcal{p}}
\newcommand{\meridian}{\mathcal{m}}

\counterwithin{figure}{section}
\numberwithin{equation}{section}

\title[Complete intersection semigroups and end curves ]{Complete intersection semigroups and end curves on singularities of splice type}

\author{Patricio Almir\'on}
\author{Julio-Jos\'e Moyano-Fern\'andez}
\author{Adri\'an Olivares-Fern\'andez}



\thanks{The authors are partially funded by MCIN/AEI/10.13039/501100011033 and by ``ERDF -- A way of making Europe", grants PID2022-138906NB-C21 and PID2022-138906NB-C22.The first author also acknowledge financial support of the Department of Education of the Junta de Castilla y León and FEDER Funds (Reference: CLU-2025-1-02- IMUVA). The second author has been partially supported also by Universitat Jaume I, grants UJI-B2021-02, GACUJIMA/2023/06 and GACUJIMB/2023/03. The third author has been partially suported by Conselleria d’Educació, Cultura, Universitats i Ocupació of Generalitat Valenciana and European Social Fund Plus (ESF+), grant CIACIF/2024/504. }

\address{Departamento de \'Algebra, An\'alisis Matem\'atico, Geometr\'ia y Topolog\'ia; IMUVA (Instituto de Investigaci\'on en Matem\'aticas), Universidad de Valladolid\\ 
	Paseo de Bel\'en 7\\
	47011 Valladolid, Spain.}
\email{palmiron@uva.es}

\address{Universitat Jaume I, Campus de Riu Sec, Departamento de Matem\'aticas \& Institut Universitari de Matem\`atiques i Aplicacions de Castell\'o, 12071
	Caste\-ll\'on de la Plana, Spain}
\email{moyano@uji.es}
\email{aolivare@uji.es}

\begin{document}
	
	\begin{abstract}
		We prove that any complete intersection numerical semigroup can be realized as the value semigroup of an end curve in a splice quotient singularity. Therefore, providing a converse to the semigroup-theoretic statement of the End Curve Theorem of Neumann and Wahl for integral homology sphere links.  For a certain family of curves embedded on splice quotient singularities with $\mathbb{Z}HS^3$ link, we show the coincidence between the Alexander polynomial of the knot cutted out by the curve and its Poincaré series. Hence providing a generalization of the result due to Campillo, Delgado and Gusein-Zade for plane curve singularities. Moreover, for any complete-intersection numerical semigroup, we define a family of knots such that the Alexander polynomial of the knots coincides with the generating series of the semigroup.
	\end{abstract}
	
	\maketitle

	\tableofcontents
	\section{Introduction}

	Let $(X,0)$ be a normal complex surface singularity whose link $\Sigma$ is an integral homology sphere, and let $\Lambda$ denote its splice diagram. The End Curve Theorem of Neumann and Wahl \cite{neumann2005complex,neumann2010end} characterizes singularities of splice type in terms of the realization of the end knots of $\Lambda$. More precisely, $(X,0)$ is of splice type if and only if every end knot associated with a leaf of $\Lambda$ is realized as the link of a hypersurface section. An \emph{end-curve function} associated with a leaf $w_i$ of $\Lambda$ is a holomorphic germ $f_i\in\mathcal{O}_{X,0}$ such that \( C_i: f_i=0 \) is an irreducible curve germ whose strict transform under a good resolution of $(X,0)$ is smooth, meets transversely the exceptional component corresponding to $w_i$ at a single point, and is disjoint from all other exceptional components. Equivalently, $C_i\cap\Sigma$ is the end knot represented by $w_i$. 
	\medskip
	
	Associated with each end curve $C_i$ are two numerical semigroups: its analytic value semigroup and a topologically defined semigroup obtained from the linking numbers encoded by $\Lambda$. A fundamental ingredient in the proof of the End Curve Theorem is the identification of these two semigroups. This identification implies the semigroup condition for $\Lambda$ and leads to splice-diagram equations defining $(X,0)$ as a complete intersection of splice type. Moreover, the end curves arising in this setting are complete-intersection curves, and hence their value semigroups are complete-intersection numerical semigroups. Consequently, for a singularity of splice type with integral homology sphere link, the value semigroup of each end curve is determined entirely by the corresponding splice-diagram. This leads naturally to the converse problem: which complete-intersection numerical semigroups occur as value semigroups of end curves of such surface singularities? 
	\medskip
	
	Recall that a \emph{numerical semigroup} is an additive submonoid $S\subseteq\mathbb{N}$ with finite complement \(|\mathbb{N}\setminus S|<\infty\). If \(S=\langle a_0,\ldots,a_g\rangle, \) its semigroup algebra is the graded monomial algebra \[ \mathbb{C}[S] :=\mathbb{C}[t^{a_0},\ldots,t^{a_g}] =\bigoplus_{\alpha\in S}\mathbb{C}t^\alpha. \] The semigroup $S$ is called a \emph{complete-intersection numerical semigroup} if $\mathbb{C}[S]$ is a complete-intersection ring. Delorme \cite{Delormegluing} characterized such semigroups combinatorially in terms of a recursive gluing construction. We recall the form of this characterization in Section~\ref{section:three}; see also \cite{AlmironCI} and the references therein.
	\medskip
	
	Our first main result gives a converse to the semigroup-theoretic part of the End Curve Theorem for integral homology spheres. We prove that every complete-intersection numerical semigroup occurs as the value semigroup of an end curve on a normal surface singularity with integral homology sphere link. More precisely, given a complete-intersection numerical semigroup $S$, we construct a splice diagram with a distinguished root whose associated topological semigroup is $S$ and which satisfies the semigroup and edge determinant conditions. Consequently, this diagram is realized by a singularity of splice type in which the end curve associated with the root has value semigroup $S$ (Theorem~\ref{teo:realizationCI}).
	\medskip
	
	Our result further develops the construction of Neumann and Wahl \cite{neumann2005complex}, which closely parallels Delorme's recursive construction at the level of the semigroups associated with end knots. Delorme's recursive decomposition is used in \cite[Lemma~3.3(6)]{neumann2005complex} to show that the topological semigroup $S(\Lambda_w)$ associated with an end knot $w$ is a complete-intersection numerical semigroup. The same recursive structure then underlies the construction of explicit complete-intersection equations for a curve realizing the prescribed end knot. The proof of \cite[Lemma~3.3(6)]{neumann2005complex}, however, does not make the recursive decomposition of $S(\Lambda_w)$ fully explicit. In Theorem~\ref{teo:semigroup_rooted}, we give a direct and explicit proof that identifies the precise steps of this decomposition relevant to the converse realization problem. This provides the starting point for our construction and reveals the principal obstruction to a naive converse. 
	\medskip

	Indeed, as shown in Subsection~\ref{ssec:construction_knot}, reversing this construction naively produces a splice diagram whose weights need not satisfy the edge determinant condition, so that it is not the splice diagram of a normal complex surface singularity \cite[Section~24]{EN}. In Subsection~\ref{subsec:negativedefinite} we address this obstruction by modifying only the weights pointing towards the distinguished root. Since modifying these weights does not alter the linking numbers with the root, the prescribed semigroup is preserved, while all edge determinants become positive. This leads to Theorem~\ref{teo:realizationCI} and, in particular, to a topological characterization of complete-intersection numerical semigroups as the semigroups of distinguished roots of splice diagrams satisfying the semigroup and edge determinant conditions.
	\medskip
	
	While the End Curve Theorem extends to normal surface singularities with rational homology sphere links \cite{neumann2010end}, its formulation in that setting involves the action of the discriminant group on the universal abelian cover. The resulting end-curve therefore carries equivariant information, and it remains unclear in what form Delorme's recursive construction can be incorporated into this framework. The converse to \cite[Lemma~3.3(6)]{neumann2005complex} established in subsection \ref{subsec:negativedefinite} highly depends on the semigroup-theoretic part of the construction and may thus provide a natural starting point for considering analogous realization problems in this broader setting. However, the integral homology sphere case is naturally adapted to the approach developed in this paper. Indeed, if $\Sigma$ is an integral homology sphere, then $H_1(\Sigma;\mathbb{Z})=0$, so the discriminant group is trivial and the universal abelian cover coincides with the singularity itself. Consequently, the notions of splice quotient singularity and singularity of splice type coincide, and the end-curve realization problem can be formulated directly in terms of the numerical semigroups encoded by the splice diagram making Delorme's recursive construction directly applicable.
	\medskip
	
	Our topological characterization also allows us to revisit the classical relation among value semigroups, Poincaré series, and Alexander polynomials. More generally, the same combinatorial framework applies to reduced curve germs $(C,0)\subset (X,0)$ whose irreducible components are end curves. Our second main result establishes that the Poincaré series associated with the curve valuations agrees with the Alexander polynomial of the corresponding curve link.
	\medskip
	
	This identity originates in the work of Campillo, Delgado, and Gusein-Zade \cite{CDGduke,CDG99b,CDG99a,CDG00}, who established, in several forms, the equality between the Poincaré series of a plane curve singularity and the Alexander polynomial of its link. They later extended this correspondence to rational surface singularities \cite{CDG04} as well as curves embedded on them \cite{CDG05}. These results initiated a broader study of identities relating analytic invariants defined by valuation filtrations to topological invariants extracted from the resolution graph. 
	\medskip
	
	For normal surface singularities, the appropriate formulation replaces the Alexander polynomial by a topologically defined multivariable series $\zeta$ of A'Campo type. Némethi \cite{ZPNemethi,nemethi2012cohomology} placed these identities in a unified framework by comparing $\zeta$ with an analytic Poincaré series $P$. When the link is an integral homology sphere, the relevant relative specialization of $ \zeta$ agrees with the Alexander polynomial of the curve link. It is worth noting that this equality is not automatic: there exist curves embedded in splice-quotient singularities for which the Poincaré series does not agree with the Alexander polynomial; see, for instance, \cite[Section~7.1]{ZPNemethi}. As a consequence, identifying geometrically meaningful families of curves for which the two invariants agree is a nontrivial problem, as it depends on the analytic structure of \((X,0).\)
	\medskip
	
	Our approach builds on that developed by the first and third authors in \cite{TopologicalAlmironolivares}. As in that work, the guiding idea is to formulate the question as an inverse problem and to exploit the relation between splice formulas for the Alexander polynomial and the way Delorme's recursive structure is reflected in the Poincaré series of an irreducible curve germ.
	\medskip
	
	In Subsection~\ref{ssec:construction_knot}, we associate with every complete-intersection numerical semigroup a family of knots sharing the same Alexander polynomial. This polynomial agrees, up to the factor $1-t$, with the generating series of the semigroup, and as a consequence we extend the results of \cite{TopologicalAlmironolivares}. We then show that the modifications of the splice-diagram weights required to satisfy the edge determinant condition leave the Alexander polynomial of the distinguished end knot unchanged. Once the resulting splice diagram is realized by a normal surface singularity, the prescribed semigroup becomes the value semigroup of an end curve, and its generating series becomes the Poincaré series of the corresponding curve valuation. The resulting series therefore admits both a combinatorial description and an analytic interpretation.
	\medskip
	
	In Section~\ref{section:four}, we focus on the interplay between the
	combinatorics of the splice diagram and the analytic filtrations of
	$\mathcal{O}_{X,0}$, namely the divisorial filtration associated with an exceptional divisor and the filtration defined by an embedded curve, through their semigroups of values and their Poincaré series. First, we exploit the results from Section \ref{section:three} to determine the semigroup of integer multiplicities of functions along a divisor of the resolution (Theorem \ref{thm:semigroup_divisors}). Second, from the divisorial identities supplied by our splice-diagram calculations, we pass to the corresponding curve valuations through a sequence of divisorial approximations. Applying recent results of Némethi and Veys \cite{Nemethi-Veys_filtrations}, we will establish the relative identity for the resulting reduced curve germs obtained as limit of divisorial valuations. As a particular case, this recovers the result of Campillo, Delgado and Gusein-Zade for plane curve singularities \cite{CDGduke} (see Remark~\ref{rem:plane_curves}). 
	\medskip
	
	The general conclusion of this work is that Delorme's recursive structure admits simultaneous combinatorial, topological, and analytic realizations. They also reinforce the close relation between the gluing operation for numerical semigroups and the splice operation for rational homology sphere links, previously explored by the authors in \cite{AlmironCI,AMknots,TopologicalAlmironolivares}. Starting from a complete-intersection numerical semigroup, we construct a splice diagram satisfying the required topological conditions, realize the prescribed semigroup as the value semigroup of an end curve, and identify its generating series with both the Poincaré series of the corresponding curve valuation and the Alexander polynomial of the end knot. The same principle extends to reduced curve germs whose irreducible components are end curves, yielding new families for which the identity $ \zeta=P$ holds.
	
	
	\section{Preliminaries}\label{section:pre}
	Let \(L=S_1\cup \cdots \cup S_n\) be a link in a homology \(3\)--sphere \(\Sigma,\) i.e. \(n\) smooth disjoint oriented copies of simple closed curves $S_i$ inside a homology 3-sphere $\Sigma$; this is therefore defined by a couple denoted $(\Sigma,L).$  For two given links $(\Sigma,L)$ and $(\Sigma',L')$,  Eisenbud and Neumann introduced in \cite{EN} an operation called \emph{splicing}, which encompasses some classical constructions, as cabling and connected sum of links, as particular cases. We briefly recall its definition and refer to \cite{EN} for further properties.
	Let $S$ and $S$ be components of $L$ and $L'$, respectively. Let $N(S)$ and $N(S')$ denote tubular neighborhoods of $S$ and $S'$. The link exterior of $L$ will be referred as $\Sigma\setminus N(L)$. Note that the boundary $\partial N(S)$ of $N(S)$ (resp. $\partial N(S')$ of $N(S')$) is homeomorphic to a torus, therefore we can choose two distinguised simple closed curves $(\meridian,\longitude)$ (resp. $(\meridian',\longitude')$), determined up to isotopy by the following properties:
	\vspace{1mm}
	
	\begin{equation*}
		\begin{aligned}
			\meridian &\sim 0,\qquad
			\longitude \sim S
			&& \text{in } H_1(N(S));\\
			\ell(\longitude,S)&=0,\qquad
			\ell(\meridian,S)=1.
		\end{aligned}\\
	\end{equation*}
	\vspace{1mm}
	
	where $\ell(\cdot ,\cdot )$ denotes the linking number. We refer to the curves $\meridian$ (resp. $\longitude$) as topological standard meridian (resp. topological standard longitude). Eisenbud and Neumann construct a closed $3$-manifold $\Sigma''$ by gluing $\Sigma\setminus \operatorname{int}(N(S))$ and $\Sigma'\setminus \operatorname{int}(N(S'))$ along their boundary tori so that $\meridian$ is identified with $\longitude'$ and $\meridian'$ with $\longitude$.  The pair
	\[
	\left(\Sigma'',(L\setminus S)\cup(L'\setminus S')\right)
	\]
	is called the \emph{splice} of $(\Sigma,L)$ and $(\Sigma',L')$ along the components $S$ and $S'$.
	\medskip
	
	Although the construction apply to arbitrary links, our main interest lies in the class of \emph{graph links}, namely those links whose exterior has a JSJ decomposition consisting entirely of Seifert fibered pieces. In the particular case where the JSJ decomposition consists of a single Seifert fibered piece, the pair $(\Sigma,K)$ is called a \emph{Seifert link}; equivalently, $\Sigma$ admits an effective circle action and $K$ is a union of fibers of the corresponding Seifert fibration. In this situation, we simply refer to $\Sigma$ as a Seifert fibered space.
	\medskip
	
	The classification of Seifert fibered spaces that are $\mathbb{Z}$-homology $3$-spheres was obtained by Seifert in \cite[Section~12]{seifert1933topologie} (see also \cite[Section 7]{EN}). They are precisely the manifolds with (unnormalized) Seifert invariants \((0;(d_1,\Bar{d}_1),\dots,(d_n,\Bar{d}_n)),
	\) where $d_1,\dots,d_n$ are pairwise coprime integers satisfying $d_i\ge2$, $n\ge3$, and
	\begin{equation*}
		\sum_{i=1}^{n}
		\Bar{d}_i\,
		d_1\cdots\widehat{d_i}\cdots d_n
		=
		1.
	\end{equation*}
	Such a manifold will be denoted by \(\Sigma(d_1,\dots,d_n).\) Accordingly, the notation
	\[
	(\Sigma(d_1,\dots,d_n),\,S_1\cup\cdots\cup S_n)
	\]
	will always denote the link formed by the singular fibers of the corresponding Seifert fibration, where $S_i$ is the singular fiber of multiplicity $d_i$. 
	
	\subsection{Splice diagrams}\label{ssec:Splice_diagrams}
	Given the structure of the manifolds described above, Eisenbud and Neumann introduced in \cite[Section~8]{EN} a graphical tool, called a \emph{splice diagram}, to describe and characterize the topology of graph links in integer homology spheres. A splice diagram consists of the following elements.
	
	\begin{itemize}
		
		\item \textbf{Nodes.} These are vertices of valence at least $3$, decorated with a sign $\varepsilon=\pm1$. Each node represents a Seifert fibered piece of the JSJ decomposition of the ambient $\mathbb{Z}$-homology $3$-sphere, while the sign records its orientation.
		
		\item \textbf{Boundary vertices.} These are vertices of the form
		[
		\begin{tikzpicture}[baseline=-0.5ex]
			\draw (0,0) -- (1.5,0)
			node[pos=1,fill=black,circle,minimum size=5pt,inner sep=0pt]{};
		\end{tikzpicture}
		].
		They represent the tubular neighborhood of a fiber of the Seifert fibered piece corresponding to the node to which they are attached. These vertices have valency one and will be called \textit{``leaves''.}
		
		\item \textbf{Arrowhead vertices.} These are vertices of the form
		[
		\begin{tikzpicture}[baseline=-0.5ex]
			\draw[-{Classical TikZ Rightarrow[length=1.5mm]}] (0,0) -- (1.5,0);
		\end{tikzpicture}
		].
		They correspond to the components of the link.
		
	\end{itemize}
	
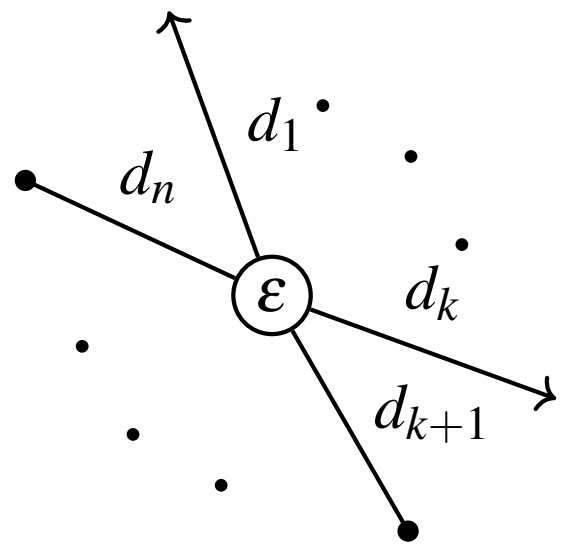
\begin{figure}[H]
		\centering
		\begin{tikzpicture}[scale=1, thick]
			\node[draw, circle, inner sep=2pt] (eps) at (0,0) {$\varepsilon$};
			
			\draw[->] (eps) -- (110:2) node[midway, right=2pt, yshift=2pt] {$d_1$};
			
			\draw (eps) -- (155:1.8) coordinate (dn);
			\node[above=2pt] at (155:0.9) {$d_n$};
			\fill (dn) circle (2pt);
			
			\draw[->] (eps) -- (340:2) node[midway, above=2pt] {$d_k$};
			
			\draw (eps) -- (300:1.8) coordinate (dk1);
			\node[right=2pt] at (300:0.9) {$d_{k+1}$};
			\fill (dk1) circle (2pt);
			
			\fill (75:1.3) circle (1.2pt);
			\fill (45:1.3) circle (1.2pt);
			\fill (15:1.3) circle (1.2pt);
			
			\fill (255:1.3) circle (1.2pt);
			\fill (225:1.3) circle (1.2pt);
			\fill (195:1.3) circle (1.2pt);
		\end{tikzpicture}
		\caption{Splice diagram of a Seifert link.}
		\label{fig:spliceseifertlink}
	\end{figure}

	Every edge either joins two nodes, in which case it represents a separating torus in the JSJ decomposition, or joins a node to a boundary vertex or an arrowhead vertex. For a vertex $v$, we denote its valence by $\delta_v$. Since we restrict our attention to $\mathbb{Z}$-homology $3$-spheres whose JSJ decomposition consists entirely of Seifert fibered pieces, the fundamental building blocks of splice diagrams are those consisting of a single node, namely Seifert fibered $\mathbb{Z}$-homology $3$-spheres. Accordingly, the Seifert link
	\[
	(\Sigma(d_1,\dots,d_n),,S_1\cup\cdots\cup S_k)
	\]
	is represented by a splice diagram with a single node and $n$ adjacent edges, the $i$-th edge being decorated with the weight $d_i$ (see Figure~\ref{fig:spliceseifertlink}). We extend the definition to admit $d_i=1,$ in which case the component of the link is represented by a non-singular fiber. 
	
	\begin{rem}
		There is a unique orientation for which the above manifold is the link of a complex surface singularity. Throughout this paper, we assume that this is the positive orientation (that is, $\varepsilon=+$). Consequently, unless explicitly stated otherwise, every node will be assumed to have positive sign, and the sign will be omitted from the splice diagrams.
	\end{rem}
	
\begin{figure}[H]
		\centering
		\begin{tikzpicture}[
			large circle/.style={
				circle,
				draw,
				minimum size=1.3cm, 
				inner sep=0pt,
				text width=1.2cm,   
				text centered,
				font=\fontsize{16pt}{18pt}\selectfont
			},
			thick
			]
			
			\node[large circle] (L1) {$\Lambda$};
			
			\node[large circle, right=4cm of L1] (L1prime) {$\Lambda'$};
			
			\draw[->] (L1.east) -- ++(1.5cm, 0); 
			
			\draw[->] (L1prime.west) -- ++(-1.5cm, 0); 
			
			\node[large circle, right=7cm of L1] (L2) {$\Lambda$};
			\node[large circle, right=4cm of L2] (L2prime) {$\Lambda'$};
			
			\draw (L2) -- (L2prime);
			
		\end{tikzpicture}
		\caption{The splicing operation on splice diagrams.}
		\label{fig:splicediagsplicing}
	\end{figure}
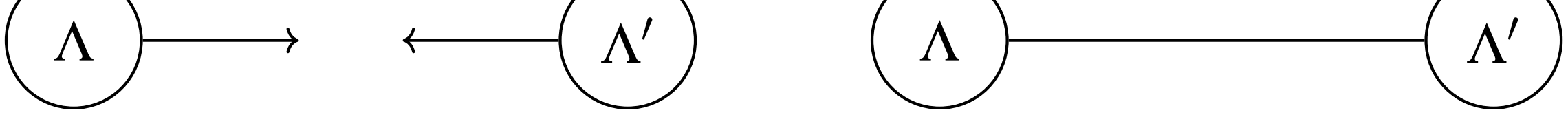

	Having described the splice diagram associated with a Seifert piece, we now explain how splice diagrams behave under the splicing operation. Let $(\Sigma,K)$ and $(\Sigma',K')$ be two links with splice diagrams $\Lambda$ and $\Lambda'$, respectively, and let $S\subset K$ and $S'\subset K'$ be distinguished components. The splice of $(\Sigma,K)$ and $(\Sigma',K')$ along $S$ and $S'$ is represented by the operation shown in Figure~\ref{fig:splicediagsplicing}, where the arrowheads indicate the corresponding link components.
	
	In general, a given graph link may admit several different splice diagrams. Nevertheless, any two such diagrams are related by a finite sequence of local moves. In the following theorem we describe the move that will be used throughout this paper; for the complete list we refer the reader to \cite[Theorem~8.1]{EN}.

	\begin{theorem}\label{teo:splice_movements}
		The following local move preserve the represented link. Assume that
		\[
		d_0d_0'
		=
		\varepsilon\varepsilon'
		d_1\cdots d_r
		d_1'\cdots d_s'.
		\]
		Then,
		
		
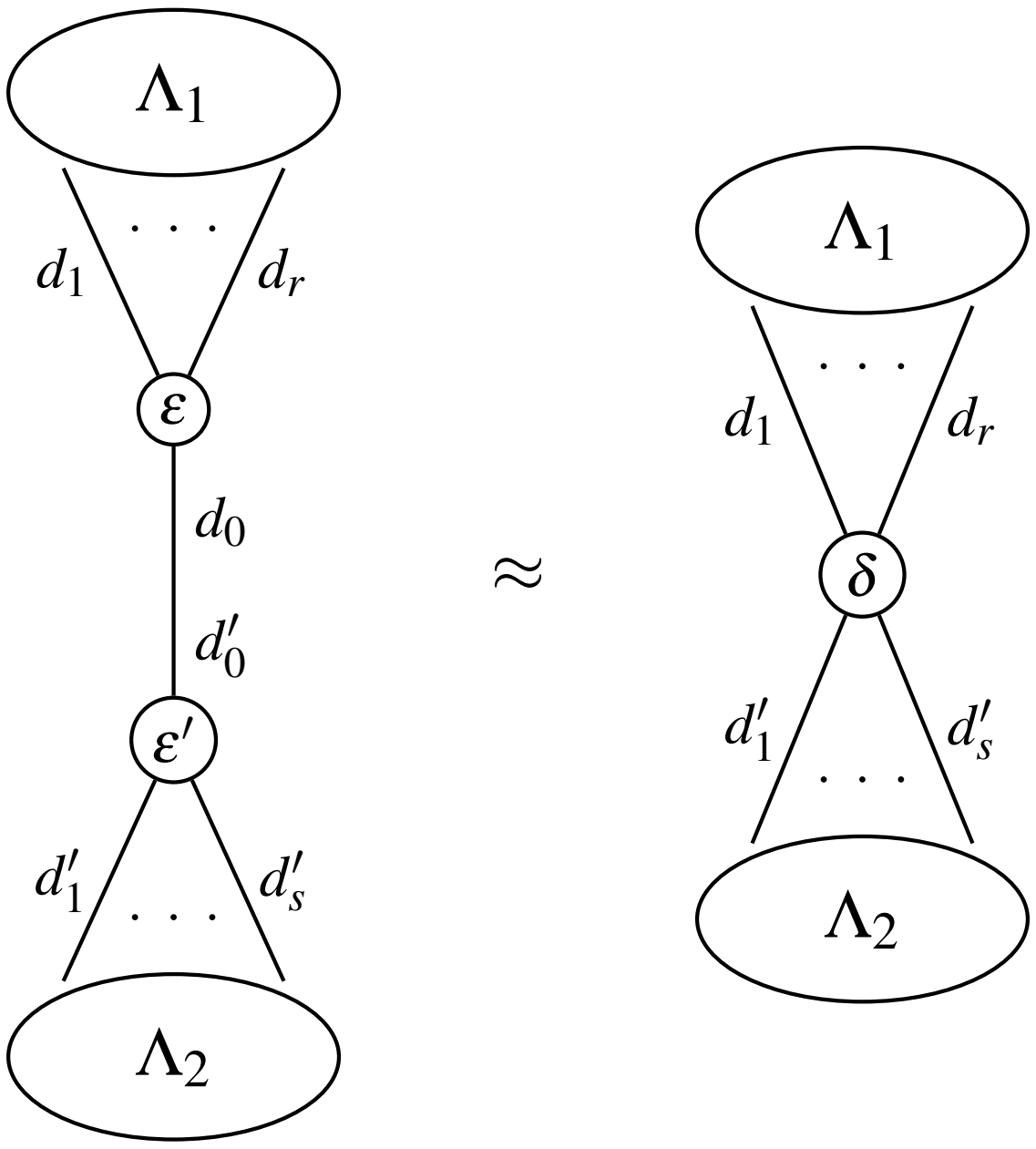
\begin{figure}[H]
			\centering
			\begin{tikzpicture}[scale=1, thick]
				
				
				\node[circle, draw, inner sep=2pt] (e) at (0, 1.2) {$\varepsilon$};
				\node[circle, draw, inner sep=1.5pt] (ep) at (0, -1.2) {$\varepsilon'$};
				
				\draw (0, 3.5) ellipse (1.2 and 0.6);
				\node at (0, 3.5) {\large $\Lambda_1$};
				
				\draw (e) -- (-0.8, 2.95) node[midway, left] {$d_1$};
				\draw (e) -- (0.8, 2.95) node[midway, right] {$d_r$};
				\node at (0, 2.5) {$\cdot \ \cdot \ \cdot$};
				
				\draw (e) -- (ep) node[pos=0.3, right] {$d_0$} node[pos=0.8, right] {$d'_0$};
				
				\draw (0, -3.5) ellipse (1.2 and 0.6);
				\node at (0, -3.5) {\large $\Lambda_2$};
				
				\draw (ep) -- (-0.8, -2.95) node[midway, left] {$d'_1$};
				\draw (ep) -- (0.8, -2.95) node[midway, right] {$d'_s$};
				\node at (0, -2.5) {$\cdot \ \cdot \ \cdot$};
				
				\node at (2.5, 0) {\large $\approx$};
				
				
				\node[circle, draw, inner sep=2pt] (delta) at (5, 0) {$\delta$};
				
				\draw (5, 2.5) ellipse (1.2 and 0.6);
				\node at (5, 2.5) {\large $\Lambda_1$};
				
				\draw (delta) -- (4.2, 1.95) node[midway, left] {$d_1$};
				\draw (delta) -- (5.8, 1.95) node[midway, right] {$d_r$};
				\node at (5, 1.5) {$\cdot \ \cdot \ \cdot$};
				
				\draw (5, -2.5) ellipse (1.2 and 0.6);
				\node at (5, -2.5) {\large $\Lambda_2$};
				
				\draw (delta) -- (4.2, -1.95) node[midway, left] {$d'_1$};
				\draw (delta) -- (5.8, -1.95) node[midway, right] {$d'_s$};
				\node at (5, -1.5) {$\cdot \ \cdot \ \cdot$};
				
			\end{tikzpicture}
		\end{figure}
		
		where $\delta=\pm1$ is chosen so that
		\[
		d_0=\delta\varepsilon\,d_1'\cdots d_s',
		\qquad
		d_0'
		=
		\delta\varepsilon'\,d_1\cdots d_r.
		\]
		
	\end{theorem}

	An alternative description of graph manifolds is provided by plumbing graphs. Indeed, every graph manifold which is a $\mathbb{Z}$-homology $3$-sphere is a plumbed manifold and can be obtained by plumbing a collection of $2$-spheres according to a weighted tree (see \cite{hirzebruch1971differentiable,milnor1958differentiable} for the definition of plumbing). Eisenbud and Neumann showed in \cite[Chapter~V]{EN} that plumbing graphs and splice diagrams are equivalent descriptions in the case of an integral homology sphere: they gave an explicit algorithm to construct the splice diagram associated with a plumbing graph, as well as its inverse (see \cite[Section 9]{neumann2005complex}). Consequently, each of these objects determines the other. A crucial feature of this correspondence is the unimodularity of the intersection matrix associated with the plumbing graph.
	\medskip
	
	Relaxing the unimodularity of the intersection matrix naturally leads to the broader class of \emph{rational homology $3$-spheres}, which we shall abbreviate by $\mathbb{Q}HS^3$. Recall that a closed $3$-manifold $\Sigma$ is a $\mathbb{Q}HS^3$ if $H_1(\Sigma,\mathbb{Z})$ is finite. As in the integral case, every graph manifold $\mathbb{Q}HS^3$ is a plumbed manifold whose plumbing graph is a weighted tree of rational curves. Moreover, if $A$ denotes the intersection matrix of the plumbing graph, then
	\[
	|H_1(\Sigma,\mathbb{Z})|=|\det(A)|.
	\]
	
	The close analogy between the plumbing descriptions of graph manifold $\mathbb{Z}$- and $\mathbb{Q}$-homology $3$-spheres led Neumann and Wahl to extend the notion of splice diagram to the rational setting; see \cite{neumann2002universal} and \cite[Section~12]{neumann2005complete}. In this more general setting, the edge weights adjacent to a node need not be pairwise coprime, and a splice diagram no longer determines the underlying plumbing graph uniquely. For instance, there are non-homeomorphic $\mathbb{Q}$-homology $3$-spheres with the same splice diagram but with different plumbing graphs. In any case, plumbing graphs provide a complete description of graph manifold $\mathbb{Q}$-homology $3$-spheres from a purely topological perspective.
	\medskip
	
	Although the topological classification of these manifolds is of independent interest, their importance is further enhanced by the fact that they arise naturally as links of normal surface singularities. A fundamental question is therefore to determine which graph manifold $\mathbb{Q}$-homology $3$-spheres can occur in this way. As proved in \cite{neumann1977brieskorn}, every Seifert manifold occurring as a node of a splice diagram is the link of a normal surface singularity, and the defining equations of these singularities are given explicitly. For general graph manifolds, the corresponding characterization can be formulated either in terms of plumbing graphs or splice diagrams. From the plumbing point of view, a graph manifold $\mathbb{Q}$-homology $3$-sphere is the link of a normal surface singularity if and only if the associated intersection matrix is negative definite (\cite{grauert1962modifikationen,mumford1961topology}). Within the framework of splice diagrams, however, the analogous criterion requires the introduction of several additional notions, which we now recall.
	\medskip
	\begin{defin}
		For any subdiagram of a splice diagram $\Lambda$ of the form shown in Figure~\ref{fig:edge determinant}, the \emph{edge determinant} is defined to be
		\[
		d_0d_0'-\varepsilon\varepsilon' d_1\cdots d_k d_1'\cdots d_r'.
		\]
		
		\begin{figure}[H]
			\centering
			\begin{tikzpicture}[scale=1, thick]
				
				\node[circle, draw, minimum size=0.7cm, inner sep=0pt] (e) at (0,0) {$\varepsilon$};
				\node[circle, draw, minimum size=0.7cm, inner sep=0pt] (ep) at (3,0) {$\varepsilon'$};
				
				\draw (e) -- (ep) node[pos=0.2, above] {$d_0$} node[pos=0.8, above] {$d'_0$};
				
				\draw (e) -- (-1.5, 1) node[pos=0.5, above=1pt] {$d_1$};
				\draw (e) -- (-1.5, -1) node[pos=0.5, below=1pt] {$d_k$};
				
				\fill (-1.1,  0.35) circle (1.2pt);
				\fill (-1.2,  0.00) circle (1.2pt);
				\fill (-1.1, -0.35) circle (1.2pt);
				
				\draw (ep) -- (4.5, 1) node[pos=0.5, above=1pt] {$d'_1$};
				\draw (ep) -- (4.5, -1) node[pos=0.5, below=1pt] {$d'_r$};
				
				\fill (4.1,  0.35) circle (1.2pt);
				\fill (4.2,  0.00) circle (1.2pt);
				\fill (4.1, -0.35) circle (1.2pt);
				
			\end{tikzpicture}
			\caption{}
			\label{fig:edge determinant}
		\end{figure}
	\end{defin}
	
	\begin{theorem}[{\cite[Theorem~4.1]{pedersen2011splice}}]\label{teo:singularity_cond}
		Let $\Sigma$ be a graph manifold $\mathbb{Q}HS^3$ and let $\Lambda$ be the splice diagram defined with the plumbing graph of $\Sigma$. Then, $\Sigma$ is the link of a complex surface singularity if and only if the following conditions hold:
		
		\begin{enumerate}
			\item Every edge determinant in $\Lambda$ is positive.
			
			\item Every weight in $\Lambda$ decorating a node or an edge is positive.
		\end{enumerate}
	\end{theorem}

	As a consequence, there is an obvious close connection between the edge determinant conditions and the negative definiteness of the intersection matrix of the corresponding plumbing graph. However, to recover the intersection matrix of the plumbing graph may be not immediate from the splice diagram. To solve this problem, a particular type of splice diagram, namely the \emph{maximal splice diagram},  was introduced in \cite[section 12]{neumann2005complete}. We refer to \cite[section 12]{neumann2005complete} for the details about its construction.
	\medskip
	
	The maximal splice diagram has the same underlying graph as the plumbing graph —that is, it has the same vertices and edges— but each edge is decorated with weights at both of its ends. This modification makes it possible to define linking numbers exactly as for ordinary splice diagrams. Recall that
	
	\begin{defin}\label{def:linking_number}
		Given two vertices $v_i$ and $v_j$ of $\Lambda$, let $\sigma_{ij}$ denote the geodesic (unique simple path) joining (and including) them in $\Lambda$. The \emph{linking number} of $v_i$ and $v_j$ is defined by
		
		$$\begin{array}{rl}
			\ell_{ij}= & \left(\text{product of all signs of nodes on }\sigma_{ij}\right)\cdot(\text{product of all edge } \\
			& \text{ weights adjacent to these nodes but not on }\sigma_{ij}) \\
			=&\displaystyle\left(\prod_{\footnotesize\text{nodes on }\sigma_{ij}}\varepsilon\right)
			\left(\prod_{\footnotesize\substack{\text{edge weights adjacent}\\
					\text{to those nodes but}\\
					\text{not lying on }\sigma_{ij}}}d_w\right).
		\end{array}$$

		Furthermore, we denote by $\delta_j$ the valence of the vertex $v_j$. Thus, $\delta_j=1$ for arrowhead and boundary vertices, whereas $\delta_j\ge3$ for nodes.The linking numbers between any two vertices of valence other than two of a maximal splice diagram coincides with the corresponding linking number in the associated splice diagram.
	\end{defin}
	\begin{theorem}[{\cite[Theorem 12.2]{neumann2005complete}}]\label{teo:adjoint_matrix}
		Let $\Gamma$ be the plumbing graph for a $\mathbb{Q}HS^3$ and let $A(\Gamma)$ denote its intersection matrix. Consider $\Lambda$ the maximal splice diagram obtained from $\Gamma$ and let $L:=(\ell_{vw})$ be the matrix whose entries are the linking numbers from $\Lambda$. Then, 
		$$A(\Gamma)L=-\det(-A(\Gamma))I$$
		where $I$ denotes the identity matrix.
	\end{theorem}
	
	\subsection*{Alexander polynomial} Splice diagrams behave particularly well with respect to classical invariants associated to links. One of them is the Alexander polynomial, originally introduced by J.~W. Alexander in \cite{alexander1928topological} as an invariant of classical links. For classical links, the Alexander polynomial can be computed by several different methods, here we adopt a definition that extends naturally to links embedded in arbitrary $\mathbb{Z}$-homology $3$-spheres.
	\medskip
	
	Let $(\Sigma,K=S_1\cup\cdots\cup S_n)$ be a link, and let $\Sigma_0$ denote its exterior. Since $\Sigma$ is a $\mathbb{Z}$-homology $3$-sphere,
	\[
	H_1(\Sigma_0)\cong\mathbb{Z}^n.
	\]
	Let \(\pi:\widetilde{\Sigma}_0\longrightarrow\Sigma_0\) be the universal abelian covering of $\Sigma_0$, let $b\in\Sigma_0$ be a base point, and write $\widetilde{b}=\pi^{-1}(b)$ for its fiber. The group of covering transformations is naturally isomorphic to $H_1(\Sigma_0)$. Choosing $t_i$ to be the generator corresponding to an oriented meridian of the component $S_i$ endows \(H_1(\widetilde{\Sigma}_0,\widetilde{b})\) with the structure of a $\mathbb{Z}[t_1^{\pm1},\ldots,t_n^{\pm1}]$-module. The first Fitting ideal of $H_1(\widetilde{\Sigma}_0,\widetilde{b})$
	is generated by the $(g-1)\times(g-1)$ minors of a presentation matrix of the module, where $g$ denotes the number of generators. The Alexander polynomial
	\[
	\Delta_K(t_1,\ldots,t_n)
	\]
	is then defined as the greatest common divisor of the generators of this ideal.
	
	\begin{rem}\label{rem:Alexnormalized}
		The Alexander polynomial is determined only up to multiplication by a unit of
		\[
		\mathbb{Z}[t_1^{\pm1},\ldots,t_n^{\pm1}].
		\]
	\end{rem}
	
	The Alexander polynomial can be obtained directly from the combinatorics of the splice diagram and the linking numbers defined in Definition \ref{def:linking_number}. Let
	
	\[ v_1,\ldots,v_n,v_{n+1},\ldots,v_k \]
	be the vertices of a splice diagram $\Lambda$, where $v_1,\ldots,v_n$ are the arrowhead vertices and $v_{n+1},\ldots,v_k$ are the remaining vertices. Then,
	
	\begin{theorem}\cite[Theorem 12.1]{EN}\label{teo:alexpol_splice}
		Assume $\Lambda$ is connected. Let $L$ be the link represented by $\Lambda.$ Then,
		\[ \Delta_L(t_1,\dots, t_n)=\left\{\begin{array}{lr}
			\displaystyle\prod_{j=n+1}^k \left(t_1^{\ell_{1j}}t_2^{\ell_{2j}}\cdots t_n^{\ell_{nj}}-1\right)^{\delta_j-2},&  n\geq2\\
			& \\
			(t_1-1)\displaystyle\prod_{j=2}^k \left(t_1^{\ell_{1j}}-1\right)^{\delta_j-2},& n=1
		\end{array}\right.\]
		where $\Delta_L(\cdot)$ denotes the Alexander polynomial for L.
	\end{theorem}
	
	\subsection{Splice-quotient singularities}\label{ssec:splice-quotient}
	
	Splice-quotient singularities (see Definition~\ref{def:splice-quotient}) were introduced by Neumann and Wahl in \cite{neumann2005complete}. They are a family of complex surface singularities whose splice diagram satisfy some special combinatorial conditions, the \emph{Semigroup Condition} and the \emph{Congruence Condition}. Let us recall their definition.
	\medskip
	
	Let $(X,0)$ be the germ of a complex analytic surface with an isolated singularity whose link $\Sigma$ is a $\mathbb{Q}$-homology $3$-sphere (i.e. a $\mathbb{Q}HS^3$), and let $\Lambda$ be its associated splice diagram. Let $v$ be one of its nodes. For each edge $e$ adjacent to $v$, denote by $d_{ve}$ the weight decorating $v$ along the edge $e$, and let
	\[
	d_v=\prod_{e} d_{ve}
	\]
	be the product of all weights adjacent to $v$. For an edge $e$ adjacent to $v$, let $\Lambda_{ve}$ denote the subdiagram of $\Lambda$ lying beyond $e$, that is, in the direction determined by $e$ away from $v$ (see Figure~\ref{fig:splice_edge}). Finally, for every leaf $w$ of $\Lambda$, let $\ell_{vw}$ be the linking number defined in Definition~\ref{def:linking_number}, and let $\ell'_{vw}$ be defined analogously, except that the weights adjacent to $v$ and pointing away from $w$ are omitted from the product.
	
	\begin{defin}\label{def:semigroup condition}
		
		We say that an adjacent edge $e$ of a node $v$ of $\Lambda$ satisfies the semigroup condition if $d_{ve}$ is in the following semigroup:
		\begin{equation*}
			d_{ve}\in \mathbb{N}\langle \ell_{vw}' : w\text{ a leaf in } \Lambda_{ve}\rangle
		\end{equation*}
		Moreover, $\Lambda$ satisfies the semigroup condition if for every node, all the edges satisfy the semigroup condition.
	\end{defin}
	\begin{rem}
		Observe that the semigroup condition can be defined for any splice diagram, not necessarily one associated to a surface singularity. \end{rem}
	\begin{figure}[H]
		\centering
		\begin{tikzpicture}[scale=1, thick]
			
			\tikzset{
				dot/.style={fill=black, circle, inner sep=0pt, minimum size=3pt},
				node dot/.style={fill=black, circle, inner sep=0pt, minimum size=4.5pt},
				squigarrow/.style={->, decorate, decoration={snake, amplitude=.5mm, segment length=3mm, post length=1.5mm}}
			}
			
			\node[node dot] (node) at (0,0) {};
			
			\node (v) at (0.3, 1.3) {$v$};
			\draw[squigarrow] (v) -- (0, 0.15); 
			
			\draw (node) -- (-1.8, 1.2);
			\draw (node) -- (-1.8, -1.2);

			\draw[dotted] (-1.8, 1.2) -- (-2.3, 1.5);
			\draw[dotted] (-1.8, -1.2) -- (-2.3, -1.5);
			
			\node[dot] at (-1.3, 0.4) {};
			\node[dot] at (-1.4, 0) {};
			\node[dot] at (-1.3, -0.4) {};
			
			\node[below] at (0.4, -0.05) {$e$};
			
			\draw (node) -- (1.5, 0);
			\draw[dotted] (1.5, 0) -- (3, 0);
			
			\draw[decorate, decoration={brace, amplitude=10pt, mirror}] 
			(0.2, -0.7) -- (2.8, -0.7)
			node[midway, below=15pt] {$\Lambda_{ve}$};
			
		\end{tikzpicture}
		\caption{}
		\label{fig:splice_edge}
	\end{figure}
	In other words, the semigroup condition states that for each edge $e$ at a node $v$ we can write 
	\begin{equation}\label{eq:semigroup_cond_1}
		d_{ve}=\sum \alpha_{vw}\ell_{vw}' \hspace{1cm} \alpha_{vw}\in \mathbb{Z}_{\geq0},
	\end{equation}
	where the sum goes over all the leaves $w$ in $\Lambda_{ve}$. As a consequence,
	\begin{equation*}
		d_{v}=\sum \alpha_{vw}\ell_{vw}.\\
	\end{equation*}
	As \(\Lambda\) is the splice diagram associated to a surface singularity then by Theorem \ref{teo:singularity_cond} it satisfies the positive edge determinant conditions. If we further assume that it satisfies the semigroup condition, then we may associate to each leaf $w$ in $\Lambda$ a variable $z_w$. At a node $v$ we define the filtration (called $v-$filtration) in the polynomial ring  that assigns to the variable $z_w$ the $v-$weight $\ell_{vw}$. Then, we say that a monomial $M_{ve}=\prod_{w\in\Lambda_{ve}}z_w^{\alpha_{vw}}$ is an admissible monomial (for $e$ at $v$) if the exponents satisfy Equation \ref{eq:semigroup_cond_1}, such a monomial is weighted homogeneous (with respect to the $v-$filtration) of degree $d_v.$  
	
	\begin{defin}\label{def:splice-quotient equations}
		Let $\Lambda$ be an splice diagram as above. At each node $v$ and adjacent edge $e$ choose an admissible monomial $M_{ve}$. Splice diagram equations for $\Lambda$ consist of a collection of equations of the form 
		\begin{equation*}
			\sum_{e}a_{vie}M_{ve}+H_{vi},\hspace{1cm}\text{ $v$ a node, \hspace{1.5mm}$i=1,\dots,\delta_v-2$}
		\end{equation*}
		
		where
		\begin{itemize}
			\item For every $v$, all maximal minors of the $(\delta_v-2)\times \delta_v$-matrix $\left(a_{vie}\right)$ have maximum rank.\\
			
			\item $H_{vi}$ is a convergent power series in $z_w$ all whose monomials have v-weight $>d_v$
		\end{itemize}
	\end{defin}
	
	Definition \ref{def:splice-quotient equations} provide candidate equations for the singularity, however not any equation as above can be used; specific admissible monomials, determined by the congruence condition, are required. To do so, we need first to recall a few facts about the resolution of \((X,0).\) 
	
	Let $\pi:(\Tilde{X},0)\rightarrow (X,0)$ be a resolution. Associated with the resolution is a decorated graph, called the \emph{dual graph} (or \emph{resolution graph}) \(\Gamma\), which encodes the geometry of the exceptional curve. Its vertices correspond to the exceptional divisors, and two vertices are joined by an edge whenever the corresponding divisors intersect. Each vertex is further decorated with the self-intersection number of the corresponding exceptional divisor. Besides describing the geometry of the resolution, the dual graph also provides a plumbing graph for the link, whose intersection matrix is negative definite. In particular, the link is a $\mathbb{Q}HS^3$ if and only if the plumbing graph is a tree and every irreducible component has genus zero. Throughout the remainder of this paper, all resolution graphs will be assumed to satisfy these conditions. 
	\medskip

	Denote by $\mathcal{V}$ the set of vertices of \(\Gamma\) and by $\mathcal{E}$ the set of leaves (vertices with valence one). The lattice $L$ generated over $\mathbb{Z}$ by the exceptional divisors has an associated intersection form $A(\Gamma)$ and, since $\Tilde{X}$ has the homotopy type of $E$, $L$ may be identified with $H_2(\Tilde{X},\mathbb{Z})$. If one considers $L'$ its dual, $A(\Gamma)$ provides an inclusion of abelian groups $L\subset L'$ whose quotient is the discriminant group $D(\Gamma):= L'/L$. It has order $\det(-A(\Gamma))$ and can be identified with $H_1(\Sigma,\mathbb{Z})$. In particular, if $\{\bar{E}_v\}\subset L'$ is the dual basis to $\{E_v\}$ then, the image of the $\bar{E}_w$ such that $w\in \mathcal{E}$ generate $D(\Gamma).$\\
	\medskip
	
	Let $n=|\mathcal{E}|$ and denote $w_1,\dots,w_n$ the $n$ leaves of the diagram, an element $\tau\in D(\Gamma)$ acts on $\mathbb{C}^{n}$ by the diagonal action $[\exp(2\pi i(\tau,\bar{E}_{w_1})),\dots,\exp(2\pi i(\tau,\bar{E}_{w_n}))]$, where $(,)$ is the intersection form of the lattice. Then, if one considers the $z_{w_i}$ as the projection to coordinate $i$ on $\mathbb{C}^{n}$, the discriminant group acts on the polynomial ring $\mathbb{C}[z_{w_1},\dots,z_{w_n}]$ multiplying each monomial $\prod z_{w_i}^{\alpha_{w_i}} $ by $\exp(-2\pi i\sum (\tau,\bar{E}_{w_i})\alpha_{w_i}).$ 
	
	\begin{defin}
		A splice diagram satisfies the \emph{congruence condition} if for each node $v$ of its splice diagram we can choose for each edge $e$ an admissible monomial $M_{ve}$ such that the discriminant group of the resolution graph acts on each $M_{ve}$ multiplying by the same number.
	\end{defin}
	
	\begin{theorem}\cite[Theorem 7.2]{neumann2005complete}\label{teo:splice-quotient}
		Assume $\Gamma$ is a resolution graph satisfying the semigroup and congruence conditions. Fix splice diagram equations for $\Gamma$ with the monomials from the congruence condition. Then, these equations define an isolated complete intersection singularity $(X,0)\subset(\mathbb{C}^{|\mathcal{E}|},0)$, the discriminant group acts freely off the origin and $(X,0)$ is the universal abelian cover of the quotient by this action. 
	\end{theorem}
	The previous theorem motivated the definition of splice-quotient singularities:
	\begin{defin}\label{def:splice-quotient}
		A singularity $(Y,0)$ is a splice-quotient singularity if its universal abelian cover is as in Theorem \ref{teo:splice-quotient}, i.e. if it satisfies both the semigroup and congruence condition (respect to $(Y,0)$).
	\end{defin}
	
	Although this is the original definition of splice-quotient singularities, two alternative characterizations are now available. The first is based on Okuma's \emph{Monomial Condition} (see \cite{okuma2006universal,okuma2008geometric} and \cite[Section~13]{neumann2005complete}), while the second is given in terms of the existence of \emph{end curve functions} \cite{neumann2010end}.
	\medskip
	
	For our purposes, the relevant ingredients are the Semigroup Condition and the theory of end curve functions. Accordingly, we shall adopt Neumann and Wahl's original point of view together with the End Curve Theorem. Fix a good resolution \(\pi:(\widetilde{X},0)\longrightarrow (X,0)\) of the normal surface singularity $(X,0)$ and let $\Gamma$ be its dual graph. Let $C_i$ be a smooth irreducible curve on $\widetilde{X}$ that intersects the exceptional curve transversely at a single point of the exceptional divisor $E_i$ and is disjoint from all the remaining exceptional divisors. The curve $C_i$ induces a knot $K_i$ in the link $\Sigma$, corresponding to the divisor $E_i$. If $E_i$ is a leaf of $\Gamma$, then $C_i$ is called an \emph{end curve}, and $K_i$ the corresponding \emph{end knot}. An \emph{end curve function} for $C_i$ is a function
	\[
	f\in\mathcal{O}_{X,0}
	\]
	whose zero locus is precisely the image of $C_i$ in $X$.
	
	\begin{defin}\label{def:end-curve condition}
		A resolution graph $\Gamma$ is said to satisfy the \emph{End Curve Condition} if every leaf of $\Gamma$ admits an end curve function.
	\end{defin}
	
	When a splice diagram satisfies the End Curve Condition, its leaves no longer represent merely topological objects (namely, embedded knots in the link), but also admit an analytic realization as the zero loci of algebraic curves on the surface. According to the previous discussion, it is easy to see that the coordinate functions (or some power of them) play the role of the end curve functions for the respective leaf of the splice diagram. Additionally, if the link of the singularity is a $\mathbb{Q}HS^3$, the End Curve Theorem \cite{neumann2010end} from Neumann and Wahl states that satisfying the end curve condition is equivalent to be a splice-type singularity.
	
	\begin{theorem}[End Cuve Theorem]
		Let $(Y,0)$ be a normal surface singularity with $\mathbb{Q}HS$ link $\Sigma$. Suppose that for each leaf of the resolution diagram $\Gamma$ there exists a corresponding end curve function $z_i:(X,0)\rightarrow (\mathbb{C},0)$ which cuts an end knot $K_i\subset \Sigma$ (or end curve) for that leaf. Then, $(Y,0)$ is a splice quotient singularity and a choice of a suitable root $x_i$ of $z_i$ for each $i$ gives coordinates for the splice quotient description.
	\end{theorem}
	
	Hence, end curves and end curve functions are essential for splice-quotient singularities. In regard to that relevance, it is natural to consider each leaf as a distinguished element in the splice diagram.

	\begin{defin}
		Let $\Lambda$ be a splice diagram and $w$ a leaf in it. We define the rooted splice diagram $(\Lambda,w)$ as the splice diagram obtained from $\Lambda$ by deleting all the weights that point towards $w.$ 
	\end{defin}
	
	When considering a rooted splice diagram $(\Lambda,w)$ we lose some information, it still codifies the relevant data for the leaf $w$ nevertheless. As can be seen from Definition \ref{def:linking_number} the weights pointing in the direction of $w$ do not affect the linking numbers with $w$. Hence, the linking numbers $\ell_{wv}$ of the rooted splice diagram $(\Lambda,w)$ and the general splice diagram $\Lambda$ coincide, therefore the object below is well defined.\\
	
	\begin{defin}
		We define the semigroup of a rooted splice diagram $(\Lambda,w)$ as the set:
		\begin{equation*}
			S(\Lambda_w)=\mathbb{N}\langle \ell_{ww'}\hspace{1mm}|\hspace{1mm} w' \text{ a leaf in }  (\Lambda,w)\rangle
		\end{equation*}
		
		In general, for a leaf $w$ in a splice diagram $\Lambda$ we define its semigroup as the semigroup of $(\Lambda,w)$. Furthermore, we say that $(\Lambda,w)$ satisfies the semigroup condition if the  edge weights pointing away from $w$ satisfy the condition.
	\end{defin}
	\begin{rem}\label{rem:numerical}
		Observe that if $\Lambda$ has exactly two vertices of valency one, namely
		$w$ and a leaf $x$, then $\Lambda$ has no nodes, since a tree with a node has
		at least three vertices of valency one. Thus $\Lambda$ consists of $w$ joined
		to $x$ by an edge, $\ell_{wx}=1$  and
		$S(\Lambda_w)=\mathbb{N}$. 
		
		More generally, assume that the weights of $\Lambda$ are positive and pairwise coprime at
		each node. If the
		vertex adjacent to $w$ is a leaf $x$, then $\ell_{wx}=1$. Otherwise, let $v$
		be the node adjacent to $w$ and $e_1,\dots,e_m$ the edges at $v$ pointing
		away from $w$, with weights $d_j=d_{ve_{j}}$. Then
		$\ell_{wx}=\left(\prod_{i\neq j}d_i\right)\ell'_{vx}$ for every leaf $x$ of
		$\Lambda_{ve_j}$, and the products $\prod_{i\neq j}d_i$, $j=1,\dots,m$, are
		coprime. Iterating this argument along each $\Lambda_{ve_{j}}$ shows that
		$\gcd S(\Lambda_w)=1$. As a consequence, the assumption that the weights of $\Lambda$ are positive and pairwise coprime at each node is sufficient to show that $S(\Lambda_w)$ and  $S(\Lambda_{ve})$ are numerical semigroups.
	\end{rem}
	
	To conclude the background on splice quotient singularities, we examine how the blow up of a smooth point in a good resolution modifies the maximal splice diagram. Let $(X,0)$ be a splice quotient singularity with link $\Sigma$ a $\mathbb{Z}HS^3$ and denote $\Lambda$ the splice diagram of $\Sigma$. Then, in \cite{neumann2005complete} is presented an algorithm to build the maximal splice diagram $\Lambda^{\text{max}}_\Gamma$ from a dual graph $\Gamma$ (corresponding to a good resolution). Conversely, Neumann and Wahl in \cite{neumann2005complex} provide an algorithm to recover from $\Lambda$ the maximal splice diagram $\Lambda^{max}_{\Gamma_\text{min}}$ associated to the minimal good resolution. Recall that any good resolution $\Gamma$ can be obtained by a series of blow ups on the minimal good resolution and that blowing up a smooth point does not change the topology of the link (i.e. the splice diagram $\Lambda$ and the link $\Sigma$ remain fixed). However, since the maximal splice diagram mirrors the structure of the resolution graph, $\Lambda^{\text{max}}_\Gamma$ and $\Lambda^{\text{max}}_{\Gamma_{\text{min}}}$ generally do not coincide. Thus, we analyze how a blow up (at a smooth point) affects the maximal splice diagram. Furthermore, we proof that if $\Tilde{\Gamma}$ is obtained from $\Gamma$ by blowing up smooth points (not belonging to the intersection of exceptional divisors), then, $\Lambda^{\text{max}}_\Gamma$ arises as a subdiagram of $\Lambda^{\text{max}}_{\Tilde{\Gamma}}.$\\
	
	The argument is supported on the following property:
	\begin{theorem}\cite[Theorem 12.2]{neumann2005complete}\label{teo:maxsplice_edge_determinant}
		Let $(X,0)$ be a normal surface singularity with $\mathbb{Q}HS^3-$link $\Sigma$. Let $\pi:\Tilde{X}\rightarrow X$ be a good resolution and let $\Gamma$ be its dual graph. Then, every edge determinant of the maximal splice diagram $\Lambda^{\text{max}}_\Gamma$ is $det(\Gamma)=\lvert H_1\left(\Sigma\right)\lvert.$
	\end{theorem}
	
\begin{figure}[H]
		\centering
		\begin{tikzpicture}[scale=1, thick]
			
			\node[circle, fill, inner sep=1.5pt, label=above:{\small $e_w - 1$}] (v1) at (0.3,0) {};
			\node[circle, fill, inner sep=1.5pt, label=above:{\small $-1$}] (v2) at (1.8,0) {};
			
			\draw (v1) -- (v2);
			
			\draw (v1) -- (-1.2, 0.7);
			\draw[dotted] (-1.2, 0.7) -- (-1.7, 0.95);
			
			\draw (v1) -- (-1.2, -0.7);
			\draw[dotted] (-1.2, -0.7) -- (-1.7, -0.95);
			
			\fill (-1.0,  0.3) circle (1.2pt);
			\fill (-1.1,  0.0) circle (1.2pt);
			\fill (-1.0, -0.3) circle (1.2pt);
			
			\draw[decorate, decoration={brace, amplitude=8pt}] 
			(-2.0, -1.1) -- (-2.0, 1.1) 
			node[midway, left=10pt] {\small $\delta$ edges};
			
		\end{tikzpicture}
		\caption{$\Tilde{\Gamma}$ near $w.$}
		\label{fig:blow_up_graph}
	\end{figure}
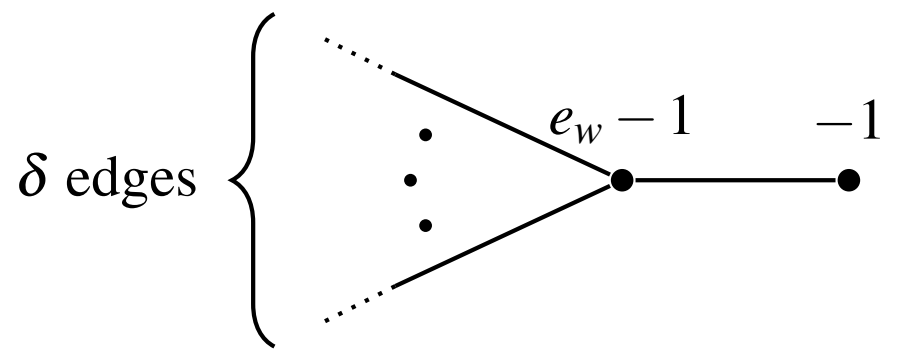

	Fix $(X,0)$ a splice quotient singularity with $\mathbb{Z}HS^3-$link and let $\Gamma$ be the dual graph from a good resolution. Choose $w$ a vertex in $\Gamma$ (simultaneously, an exceptional divisor $E_w$) and denote $\delta$ the valence of $w$. Then, blow up a smooth point of $E_w$ (we additionally suppose that such point does not belong to any other exceptional divisor). This operation provides a new good resolution for the surface singularity whose dual graph will be denoted $\Tilde{\Gamma}$. Since we have blown up a smooth point $\Tilde{\Gamma}$ has the shape of Figure \ref{fig:blow_up_graph} near $w$ ($e_w$ denotes the self-intersection number of $E_w$ in $\Gamma$) whereas the rest of $\Tilde{\Gamma}$ is identical to $\Gamma$.
	
	Compute $\Lambda^{\text{max}}_{\Gamma}$ and $\Lambda^{\text{max}}_{\Tilde{\Gamma}}$ applying the procedure from \cite[Appendix 1]{neumann2005complete}. Since blowing up a smooth point does not change the topology, the determinant of the subdiagrams do not change. Then, it is straightforward that $\Lambda^{\text{max}}_{\Gamma}$ appears as a subdiagram of $\Lambda^{\text{max}}_{\Tilde{\Gamma}}$ and  $\Lambda^{\text{max}}_{\Tilde{\Gamma}}$ has one more edge and vertex than $\Lambda^{\text{max}}_{\Gamma}$, see Figure \ref{fig:blow_up_maximal}. Moreover, the weight of the new edge at $w$ is $1$, the other weight will be called $\ell_{\Tilde{w}\Tilde{w}}$ for now.
	
	\begin{figure}[H]
		\centering
		\begin{tikzpicture}[scale=1.4, thick]
			
			\fill (0,0) circle (1.5pt);
			\fill (1.5,0) circle (1.5pt);
			
			\draw (0,0) -- (1.5,0) 
			node[pos=0.2, above] {$1$} 
			node[pos=0.8, above] { $\ell_{\widetilde{w}\widetilde{w}}$};
			
			\draw (0,0) -- (-1, 0.4) node[midway, above right=-4pt] {$d_1$};
			\draw[dotted] (-1, 0.4) -- (-1.5, 0.6);
			
			\draw (0,0) -- (-1, -0.4) node[midway, below right=-4pt] {$d_\delta$};
			\draw[dotted] (-1, -0.4) -- (-1.5, -0.6);
			
			\fill (-1.1,  0.15) circle (0.8pt);
			\fill (-1.2,  0.00) circle (0.8pt);
			\fill (-1.1, -0.15) circle (0.8pt);
			
			\draw[decorate, decoration={brace, amplitude=8pt, mirror}] 
			(-1.6, -0.8) -- (0, -0.8) 
			node[midway, below=10pt] {$\Lambda_{\Gamma}^{\max}$};
			
		\end{tikzpicture}
		\caption{$\Lambda^{\text{max}}_{\Tilde{\Gamma}}$ }
		\label{fig:blow_up_maximal}
	\end{figure}
	
	Therefore, in order to complete $\Lambda^{\text{max}}_{\Tilde{\Gamma}}$ remains to compute $\ell_{\Tilde{w}\Tilde{w}}$. Recall that we have fixed $(X,0)$ to be a splice quotient singularity whose link $\Sigma$ is a $\mathbb{Z}HS^3$. Hence, Theorem \ref{teo:maxsplice_edge_determinant} implies that all edge determinants from $\Lambda^{\text{max}}_{\Gamma}$ and $\Lambda^{\text{max}}_{\Tilde{\Gamma}}$ are $1.$ Consequently, if $d_1,\dots,d_\delta$ are the adjacent weights of $w$ in $\Lambda^{\text{max}}_{\Gamma}$ we obtain the equality
	\begin{equation*}
		\ell_{\Tilde{w}\Tilde{w}}-d_1\cdots d_\delta=1,
	\end{equation*}
	equivalently,
	\begin{equation*}
		\ell_{\Tilde{w}\Tilde{w}}=d_1\cdots d_\delta+1.
	\end{equation*}
	
	Iterating the previous procedure it is possible to compute $\Lambda^{\text{max}}_{\Gamma}$ for non-minimal (good) resolution graphs obtained by blowing up smooth points not belonging to the intersection of exceptional divisors. Therefore, \cite[Theorem 9.1]{neumann2005complex} and the discussion above provide the following
	
	\begin{lemma}\label{lemma:linking_ww}
		Under the aforementioned conditions. If $\Tilde{w}$ is a leaf in $\Lambda^{\text{max}}_{\Gamma}$, then
		\begin{equation*}
			\ell_{\Tilde{w}\Tilde{w}}=\lceil d_v/d^2 \rceil +s,
		\end{equation*}
		
		where
		\begin{itemize}
			\item $d_v$ is the product of all the edge weights adjacent to $v$ (the nearest node to $\Tilde{w}$ in $\Lambda^{\text{max}}_{\Gamma}$).\\
			
			\item $d$ is the edge weight pointing towards $\Tilde{w}$ at $v$.\\
			
			\item $s$ is: 
			\begin{enumerate}
				\item If $d=1$: the number of blow ups in the resolution between $v$ and $\Tilde{w}$; $\Tilde{w}$ included\\
				
				\item If $d>1$: the number of blow ups, starting from the minimal good resolution, performed in the direction of $\Tilde{w}$ to obtain $\Tilde{w}$.  \\
			\end{enumerate} 
		\end{itemize}
		In particular, $\lceil d_v/d^2 \rceil$ can be computed through the splice diagram $\Lambda.$
	\end{lemma} 
	
	\section{Complete intersection semigroups and End-Curves}\label{section:three}
	
	Let $\pi: (\Tilde{X},0)\rightarrow (X,0)$ be a good resolution with dual graph $\Gamma$ and $\mathbb{Z}HS^3$ link $\Sigma$. Consider $\Lambda$ the splice diagram defined by $\Gamma$. For a leaf $w \text{ in } \Lambda$, recall from Section \ref{ssec:splice-quotient} that an end curve for $w$ is a smooth irreducible curve $C_w$ transversal to the exceptional curve intersecting only the component $E_w$ in just one point. Moreover, $K_w:=C_w\cap \Sigma$ defines an embedded knot in $\Sigma$ which is represented in $\Lambda$ replacing the vertex $w$ with an arrow. This interaction between $\Lambda$ and $\Sigma$ allows to associate a numerical semigroup $S(\Lambda_w)$ to the knot $K_w$ , i.e, the semigroup associated to the rooted splice diagram $(\Lambda,w).$ Generally, identifying the class of $S(\Lambda_w)$ constitutes a non-trivial task. However, rooted splice diagrams satisfying the semigroup condition represent the same type.
	\medskip
	
	Neumann and Wahl stated in \cite[Lemma 3.3 (6)]{neumann2005complex} that $S(\Lambda_w)$ is a complete intersection semigroup. We give here a self-contained proof of a combinatorial generalization of this result, which makes explicit its underlying combinatorial structure. In addition, we will show that any complete intersection numerical semigroup can be realized as the semigroup of rooted splice diagram which satisfies the semigroup condition. However, these splice diagram are not in general associated to the resolution graph of any singularity. Nevertheless, we will also construct a modification of such splice diagram to obtain links of normal surface singularities with prescribed complete intersection semigroups as end curve semigroups.
	\medskip
	
	\subsection*{Delorme's characterization of complete intersection numerical semigroups}\label{ssec:complete intersections}
	Our study of end curves will heavily rely on the combinatorial properties of complete intersection numerical semigroups provided by Delorme in \cite{Delormegluing}. We will briefly recall some of those properties, for further details we refer to \cite{AlmironCI,Delormegluing,herzog1970generators} and the references therein.
	\medskip
	
	Let $S$ be a finitely generated subsemigroup of $\mathbb{N}$ such that $\gcd(S)=1$, i.e, $S$ is a numerical semigroup. Let us denote by 
	\[S=\langle a_0,\dots,a_g\rangle=\{n\in\mathbb{N}\;\mid\; n=\epsilon_0a_0+\cdots+\epsilon_{g}a_g\;\;\text{with}\;\;\epsilon_0,\dots,\epsilon_g\in\mathbb{N}\}.\]
	As $\gcd(S)=1$, then there exists a positive integer $c(S)$ such that $c(S)+\mathbb{N}\subset S$ and $c(S)-1\notin S$, we will refer to $c(S)$ as the conductor of $S.$
	
	Recall that a numerical semigroup is said to be a \emph{complete intersection numerical semigroup} if its associated graded algebra \(\mathbb{C}[S]\) is a complete intersection. The finite generation of \(S\) also provides a presentation of \(\mathbb{C}[S]\) as a quotient of the polynomial ring, wich is also useful. Consider, $\mathbb{C}[x_0,\dots,x_g]$ the polynomial ring in $g+1$ variables and the graduation  $\deg(x_i)=a_i,$ then we have a \(0\)--graded $\mathbb{C}-$algebra epimorphism defined by 
	\begin{equation*}
		\begin{aligned}
			\varphi: \mathbb{C}[x_0,\dots,x_g] &\longrightarrow \mathbb{C}[S] := \bigoplus_{\alpha\in S} \mathbb{C}t^{\alpha} \\
			x_i &\mapsto t^{a_i}.
		\end{aligned}
	\end{equation*}
	
	As a consequence, \(\mathbb{C}[S]\cong \mathbb{C}[x_0,\dots,x_g] /\operatorname{ker}(\varphi)\) and \(S\) is a complete intersection numerical semigroup if \(\operatorname{ker}(\varphi)\) is generated by a regular sequence. In \cite{Delormegluing}, Delorme provided a combinatorial characterization of complete intersection numerical semigroups by characterizing the possible degrees of the generators of \(\operatorname{ker}(\varphi)\). This is done as follows: let $G=\{a_0,\dots,a_g\}$ be a finite non-empty family of positive integers. For $L\subseteq G$ we denote by $S_L$ the submonoid of $\mathbb{N}$ generated by $L$ and $M(S_L)$ the affine graded $\mathbb{C}$-algebra with gradings defined by $L.$
	\medskip
	
	A distinguished sequence over $G$ is a pair $\left(\Pds,\Lds\right)$ such that
	\begin{itemize}
		\item $\Pds=\left(P_0,\dots,P_g\right)$ is a sequence of partitions of $G$ into non-empty sets defined recursively as follows:
		\begin{itemize}
			\item $P_g=\{\{a_i\} \hspace{1mm} | \hspace{1mm} 0\leq i\leq g\}.$
			\item For each $1\leq i\leq g$, there are $L_i,L_i'\in P_i$ such that $L_i\neq L_i'$ and:
			$$P_{i-1}=\left(P_i\setminus \{L_i,L_i'\}\right)\cup \left(\{a_i\in L_i\cup L_i' \}\right).$$
		\end{itemize}
		\item $\Lds=\left(Z_1,\dots, Z_g\right)\subset M(S)$ is a sequence of binomials associated to $\Pds$ such that if $L_i\text{ and } L_i'$ are the sets used to obtain $P_{i-1} \text{ from } P_i$ then
		\begin{center}
			$Z_i=\mathbf{x}-\mathbf{y}$\hspace{2mm} with\hspace{2mm} $\mathbf{x}\in M(S_{L_i}),\hspace{1mm} \mathbf{y}\in M(S_{L_i'}) $\hspace{2mm} and\hspace{2mm} $\deg(\mathbf{x})=\deg(\mathbf{y}).$
		\end{center}
	\end{itemize}
	\begin{rem}\label{rem:families}
		The subsets of $G$ considered above (the elements of the partitions
		$P_i$, the sets $L_i,L'_i$, etc.) are understood as
		subsets of the index set $\{0,\dots,g\}$, and cardinalities are counted with
		multiplicity; accordingly, $M(S_L)$ is the polynomial ring with one variable $x_i$ for each
		index $i\in L$. In fact, repetitions may occur, as we will see in the case of the families of linking numbers attached to rooted splice diagrams (see
		Corollary~\ref{cor:teosemigroup_rooted}).
	\end{rem}
	
	It is natural to relate distinguished sequences with complete intersection semigroups. Indeed, \cite[Lemme 6]{Delormegluing} states that $S_G$ is a complete intersection semigroup if and only if there exists a distinguished sequence $\left(\Pds,\Lds\right)$ for $G$ such that $\Lds$ generates $I(G)$. Furthermore, the following results show that the above property only depends on the arithmetics of the set of generators $G.$
	
	\begin{proposition}\cite[Lemme 7]{Delormegluing}\label{prop:delorme_ds_1}
		Let \((\Pds,\Lds)\) be a distinguished sequence over \(G\) with \(g\geq 1.\) Let us denote \(e:=\gcd(L_1),\) \(e':=\gcd(L'_1),\) \(m:=\lcm(e,e')\) and by \(r=\deg(Z_1).\) Then, we have the following properties:\\
		
		\begin{enumerate}
			\item Let \(F, F'\) be minimal generators of \(I(L_1)\) and \(I(L'_1)\) respectively. Then, \(F,F',Z_1\) is a minimal system of generators of \((I(L_1),I(L'_1),Z_1).\)\\
			
			\item \((I(L_1),I(L'_1),Z_1)\) generates \(I(G)\) if and only if \(m=r.\)
		\end{enumerate}
	\end{proposition}
	
	\begin{corollary}\cite[Lemme 8]{Delormegluing}\label{coro:delorme_ds_2}
		Let \(\left(\Pds,\Lds\right)\) be a distinguished sequence over \(G\). \(\Lds\) generates \(I(G)\) if and only if the following condition is satisfied
		
		\begin{equation*}
			\text{for all}\quad 1\leq i\leq g,\quad \deg(Z_i)=\lcm(\gcd(L_i),\gcd(L'_i)).
		\end{equation*}
		
	\end{corollary}
	
	As a consequence, one can obtain a splitting theorem for complete intersection numerical semigroups. This is usually interpreted as an operation called gluing. Consider $S_1,S_2$ two numerical semigroups, let $\alpha,\beta$ be two strictly positive coprime integers such that $\alpha\in S_2\text{ and }\beta\in S_1$, the gluing of $S_1 \text{ and } S_2$ is the numerical semigroup $S:=\alpha S_1+\beta S_2.$ Moreover, if $G_1 (\text{resp. } G_2)$ is a set of generators for $S_1 (\text{resp. }S_2)$, is immediate that $\alpha G_1\cup \beta G_2$ is a set of generators for $S.$ Thus, it is straightforward to connect the operation of gluing with distinguished sequences. Indeed, let $\left(\Pds,\Lds\right)$ be a distinguished sequence for $S$ and $r=\gcd(L_1), r'=\gcd(L_1')$, it follows that $S=rS_{L_1/r}+r'S_{L_1'/r'}$. For complete intersections, gluing not only arises naturally but also characterize them completely as can be seen from Delorme's results:
	
	\begin{theorem}\cite[Proposition 9]{Delormegluing}
		Let $G$ be the minimal generator from the numerical semigroup $S_G$. Then, $S_G$ is a complete intersection if and only if one of the following conditions is satisfied:
		\begin{enumerate}
			\item $G={1}$ (hence $S_G=\mathbb{N})$.\\
			
			\item $G$ is the disjoint union of two non-empty partitions of the form $a_1G_1\text{ and }a_2G_2$, where $a_1,a_2$ are coprime, $G_1,G_2$ are the minimal generators of two complete intersection numerical semigroups satisfying the conditions $a_1\in S_{G_2}\text{ and } a_1\notin G_2,\hspace{1mm} a_2\in S_{G_1}\text{ and } a_2\notin G_1$
		\end{enumerate}
	\end{theorem}
	\begin{theorem}\cite[Proposition 10]{Delormegluing}\label{teo:gluing_CI} 
		Let $S:=\alpha S_1+\beta S_2$ be the gluing of the semigroups $S_1,S_2$. Then,
		
		\begin{enumerate}
			\item  $c(S)=\alpha c(S_1)+\beta c(S_2) +(\alpha-1)(\beta-1).$\\
			
			\item $S$ is a complete intersection if and only if $S_1\text{ and } S_2$ are complete intersection.
		\end{enumerate}
		
	\end{theorem}
	
	\subsection{Semigroup of an End-Curve}\label{subsection:31}
	Once Delorme's results have been recalled, we give a self-contained proof of a combinatorial generalization of \cite[Lemma 3.3]{neumann2005complex}, which makes explicit the gluing structure used in the sequel. This proof will also show how Delorme's partition can be translated into a particular shape of the splice diagram associated to a rooted splice diagram $(\Lambda,w)$.

	\begin{theorem}\label{teo:semigroup_rooted}
		Let $(\Lambda,w)$ be a rooted splice diagram with positive weights, pairwise coprime at each node, satisfying the semigroup condition. Then $S(\Lambda_w)$ is a complete intersection numerical semigroup.
	\end{theorem}
	\begin{proof}
		We first observe that if \(\Lambda\) has only two vertices of valency one, then both are leaves and it is trivial to see that \(S(\Lambda_w)=\mathbb{N}.\) Thus, let us first assume that $(\Lambda, w)$ with one node and $n+1\geq 3$ leaves as in Figure \ref{fig:base_case_CI}. 
		\begin{figure}[H]
			\centering
			\begin{tikzpicture}[scale=1.1, thick]
				
				
				\fill (0,0) circle (2pt); 
				
				\fill (0, 1) circle (2pt) node[above=2pt] {$w$}; 
				
				\fill (-1, -0.8) circle (2pt); 
				
				\fill (1, -0.8) circle (2pt); 
				
				
				\draw (0,0) -- (0,1);
				
				\draw (0,0) -- (-1, -0.8) node[midway, above left=2pt] {$d_1$};
				
				\draw (0,0) -- (1, -0.8) node[midway, above right=2pt] {$d_n$};
				
				\fill (-0.4, -0.85) circle (1pt);
				\fill (0,    -0.95) circle (1pt);
				\fill (0.4,  -0.85) circle (1pt);
				
			\end{tikzpicture}
			\caption{One node and $n+1$ edges}
			\label{fig:base_case_CI}
		\end{figure}
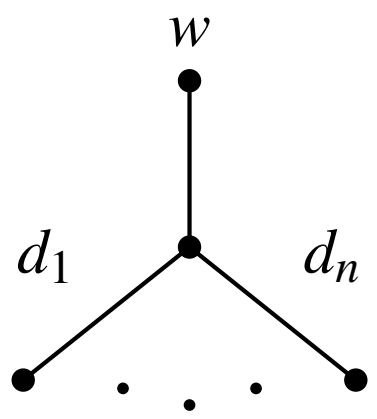
		Choose any edge of $(\Lambda, w)$ pointing away from $w,$ for example the edge of weight \(d_1\) in Figure \ref{fig:base_case_CI} and apply Theorem \ref{teo:splice_movements} to $\Lambda$ so that we obtain a new splice diagram with one extra node as in Figure \ref{fig:base_splice_proof_CI_1}. It is straightforward to see that this extra node in the new diagram provides the following decomposition at the semigroup level
		\begin{equation}\label{eq:gluing_CI_1}
			S(\Lambda_w)=d_2\cdots d_n\mathbb{N}+d_1 S(\Lambda_{ve}), 
		\end{equation} 
		where $v$ is the node adjacent to $w$ and $e$ is the new edge. As \((\Lambda,w)\) satisfies the semigroup condition then \(d_2\cdots d_n\in S(\Lambda_{ve}) \) and we have \(\gcd(d_1,d_2\cdots d_n)=1.\) As a consequence, Equation \eqref{eq:gluing_CI_1} provides a decomposition of \(S(\Lambda_w)\) as gluing of \(\mathbb{N}\) and \(S(\Lambda_{ve}).\) Thus, in virtue of Theorem \ref{teo:gluing_CI},  \(S(\Lambda_w)\) is a complete intersection numerical semigroup if and only if \(S(\Lambda_{ve})\) is.
		
		Furthermore, $\Lambda_{ve}$ is a splice diagram with only one node but one leaf less than $\Lambda$. Therefore, we may repeat the same procedure as with $\Lambda$ to obtain the equality
		\begin{equation}\label{eq:gluing_CI_2}
			S(\Lambda_{ve})=d_3\cdots d_n\mathbb{N}+d_2 S(\Lambda_{ve}'),
		\end{equation}
		where $\Lambda_{ve}'$ is constructed from $\Lambda_{ve}$ by the same procedure done in the previous case with $\Lambda.$ By the same reasons, Eq.~
		\ref{eq:gluing_CI_2} yields a gluing of semigroups. Iterating this procedure we arrive to the splice diagram in Figure \ref{fig:base_splice_proof_CI_2} and, consequently, $S(\Lambda_w)$ is a complete intersection if and only if $d_{n-1}\mathbb{N}+d_n\mathbb{N}$ is a complete itersection. As a consequence, \(S(\Lambda_w)\) is a complete intersection numerical semigroup.

		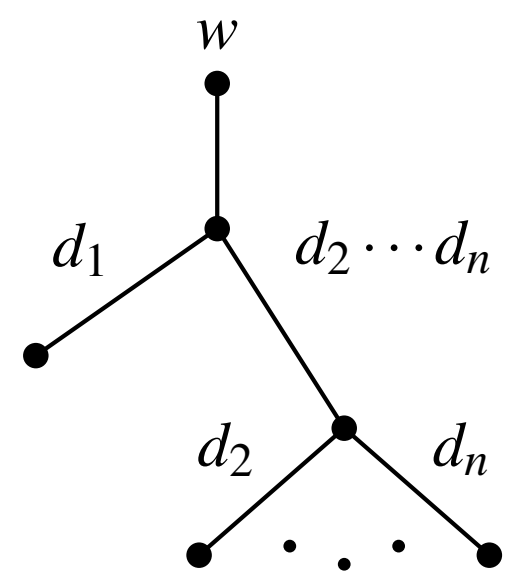
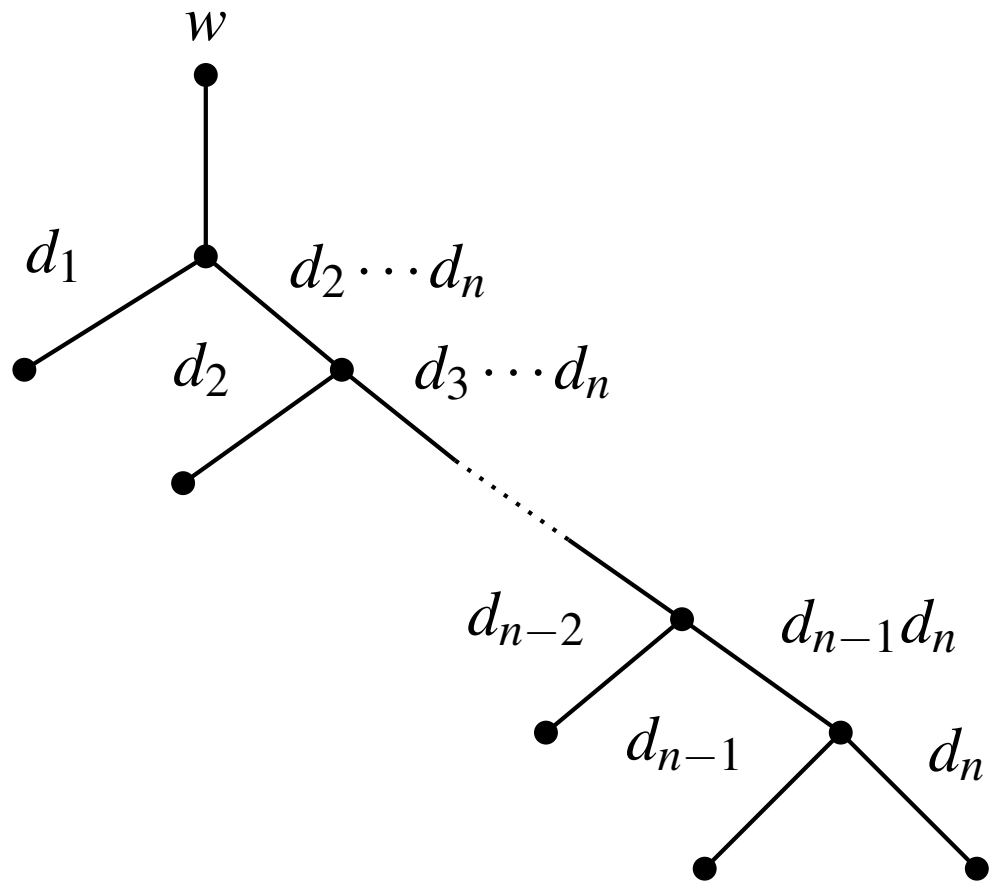
\begin{figure}[H]
			\centering
			\begin{subfigure}{0.4\linewidth}
				\centering
				\begin{tikzpicture}[scale=1.2, thick]
					
					\fill (0, 1.2) circle (2pt) node[above=2pt] {$w$}; 
					
					\fill (0, 0.4) circle (2pt); 
					
					\fill (-1, -0.3) circle (2pt); 
					
					\fill (0.7, -0.7) circle (2pt);
					
					\fill (-0.1, -1.4) circle (2pt);
					
					\fill (1.5, -1.4) circle (2pt);
					
					\draw (0, 0.4) -- (0, 1.2);
					
					\draw (0, 0.4) -- (-1, -0.3) node[midway, above left=-2pt] {$d_1$};
					
					\draw (0, 0.4) -- (0.7, -0.7);
					
					\node[right] at (0.3, 0.3) {$d_2 \cdots d_n$};
					
					\draw (0.7, -0.7) -- (-0.1, -1.4) node[midway, above left=-2pt] {$d_2$};
					
					\draw (0.7, -0.7) -- (1.5, -1.4) node[midway, above right=-2pt] {$d_n$};
					
					\fill (0.4, -1.35) circle (1pt);
					\fill (0.7, -1.45) circle (1pt);
					\fill (1.0, -1.35) circle (1pt);
					
				\end{tikzpicture}
				\caption{First step of the decomposition of \(S(\Lambda_w)\) with only one node}
				\label{fig:base_splice_proof_CI_1}
			\end{subfigure}
			\begin{subfigure}{0.4\linewidth}
				\centering
				\begin{tikzpicture}[scale=1.5, thick]
					
					
					\fill (0, 1.2) circle (1.5pt) node[above=2pt] {$w$}; 
					\fill (0, 0.4) circle (1.5pt);                       
					\fill (0.6, -0.1) circle (1.5pt);                    
					\fill (2.1, -1.2) circle (1.5pt);                    
					\fill (2.8, -1.7) circle (1.5pt);                    
					
					\fill (-0.8, -0.1) circle (1.5pt);                   
					\fill (-0.1, -0.6) circle (1.5pt);                   
					\fill (1.5, -1.7) circle (1.5pt);                    
					\fill (2.2, -2.3) circle (1.5pt);                    
					\fill (3.4, -2.3) circle (1.5pt);                    
					
					\draw (0, 0.4) -- (0, 1.2);
					
					\draw (0, 0.4) -- (-0.8, -0.1) node[midway, above left=2pt] {$d_1$};
					\draw (0, 0.4) -- (0.6, -0.1) node[midway, above right=-2pt] {$d_2 \cdots d_n$};
					
					\draw (0.6, -0.1) -- (-0.1, -0.6) node[midway, above left=2pt] {$d_2$};
					\draw (0.6, -0.1) -- (1.1, -0.5) node[midway, above right=-2pt] {$d_3 \cdots d_n$};
					
					\draw[dotted] (1.1, -0.5) -- (1.6, -0.85);
					\draw (1.6, -0.85) -- (2.1, -1.2);
					
					\draw (2.1, -1.2) -- (1.5, -1.7) node[midway, above left=1pt] {$d_{n-2}$};
					\draw (2.1, -1.2) -- (2.8, -1.7) node[midway, above right=-1pt] {$d_{n-1} d_n$};
					
					\draw (2.8, -1.7) -- (2.2, -2.3) node[midway, above left=1pt] {$d_{n-1}$};
					\draw (2.8, -1.7) -- (3.4, -2.3) node[midway, above right=-1pt] {$d_n$};
					
				\end{tikzpicture}
				\caption{Last step of the decomposition of \(S(\Lambda_w)\) with only one node}
				\label{fig:base_splice_proof_CI_2}
			\end{subfigure}
			
			\caption{Splice diagrams of the decompositions of \(S(\Lambda_w)\) with only one node} \label{figure32}
			
		\end{figure}
		
		Let us now move to the general case. We consider a splice diagram with a distinguished leaf $(\Lambda, w)$ with $k$ nodes (see Figure \ref{fig:general_splice_proof_CI_1}). Let \(v_1\) be the node adjacent to \(w,\) then we are in similar situation as before, just changing leaves by splice diagrams \(\Lambda_1,\dots, \Lambda_n.\) As a consequence, we can proceed as before, applying Theorem \ref{teo:splice_movements} to the node \(v\) to obtain a new splice diagram with one extra node which induces the semigroup decomposition
		\begin{equation}\label{eq:gluing_CI_3}
			S(\Lambda_{w})=d_2\cdots d_n S(\Lambda_{1})+d_1 S(\Lambda_{v_1e_1}),
		\end{equation}
		where $e_1$ the edge with weight $d_2\cdots d_n$ and $S(\Lambda_{1})$ is the semigroup computed from Figure \ref{fig:general_splice_proof_CI_2} for the case $i=1,$ and setting $w'$ as the root vertex. Recall that weights pointing towards $w$ do not affect in the computation of the linking numbers with $w$ and as a consequence the rooted splice $(\Lambda_{1}, w')$ satisfies the same properties as $(\Lambda, w),$ i.e it satisfies the hypothesis of the statement. As before, the Equation \eqref{eq:gluing_CI_3} provides a gluing decomposition of the semigroup \(S(\Lambda_{w}).\) Thus, Theorem \ref{teo:gluing_CI} implies that \(S(\Lambda_{w})\) is a complete intersection if and only if \(S(\Lambda_{ve})\) and \(S(\Lambda_{1})\) are complete intersections. As in the case of only one node, a recursive application of Theorem \ref{teo:splice_movements} (see Figure \ref{fig:general_splice_proof_CI_3}) allows us to provide the semigroup decompositions of the form
		\begin{equation*}
			S(\Lambda_{v_{i}e_{i}})=d_{i+2}\cdots d_n S(\Lambda_{i+1})+d_{i+1} S(\Lambda_{v_{i+1}e_{i+1}}).
		\end{equation*}
		As $(\Lambda, w)$ satisfies the semigroup condition then we have \(d_i\in S(\Lambda_{i}) \) and in this way \(d_{i+2}\cdots d_n\in S(\Lambda_{v_{i+1}e_{i+1}}),\) since \(d_{i+3}\cdots d_n S(\Lambda_{i+2})\subset S(\Lambda_{v_{i+1}e_{i+1}}).\) Therefore, it is a gluing decomposition. As a consequence, Theorem \ref{teo:gluing_CI} implies that \(S(\Lambda_{w})\) is a complete intersection if and only if \(S(\Lambda_{1}),\dots, S(\Lambda_{n})\) are complete intersection numerical semigroups. Finally, as \(\Lambda_{1},\dots, \Lambda_{n}\) are splice diagrams with at least one node less than \(\Lambda,\) we can iterate the process until we arrive to the previously considered situation of only one node, from which we deduce the required result.
		
		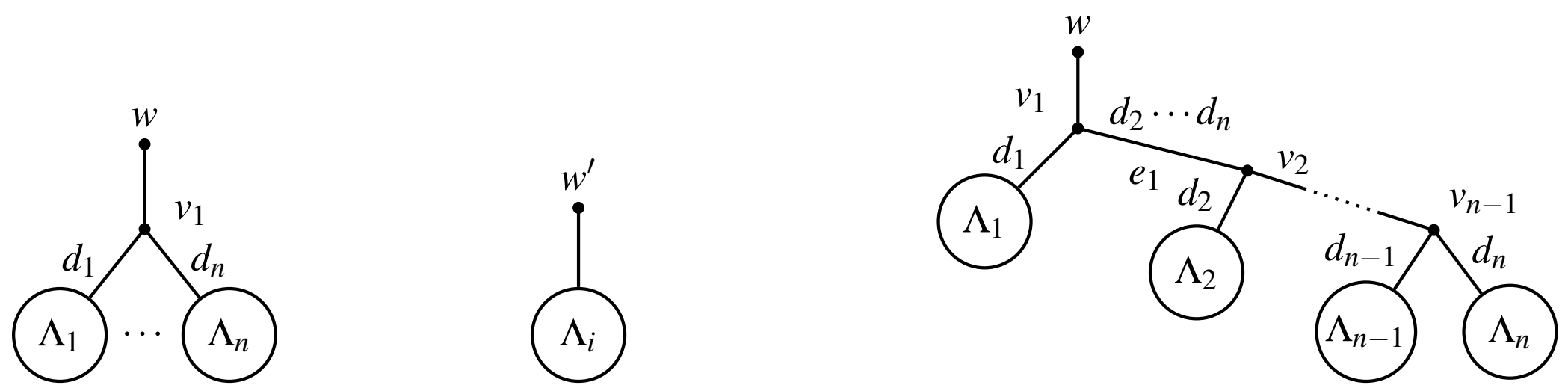
\begin{figure}[H]
			\centering
			\begin{subfigure}[b]{0.26\linewidth}
				\centering
				\begin{tikzpicture}[scale=0.75, thick, every node/.style={font=\footnotesize}]
					\coordinate (w) at (0,1);
					\coordinate (v1) at (0,0);
					\node[circle,draw,minimum size=0.82cm,inner sep=1pt] (L1) at (-1,-1.25) {$\Lambda_1$};
					\node[circle,draw,minimum size=0.82cm,inner sep=1pt] (Ln) at (1,-1.25) {$\Lambda_n$};
					\draw (w) -- (v1);
					\draw (v1) -- (L1) node[pos=0.45,left=1pt] {$d_1$};
					\draw (v1) -- (Ln) node[pos=0.45,right=1pt] {$d_n$};
					\node at (0,-1.25) {$\cdots$};
					\fill (w) circle (2pt) node[above] {$w$};
					\fill (v1) circle (2pt) node[right=3pt,yshift=4pt] {$v_1$};
				\end{tikzpicture}
				\caption{$(\Lambda,w)$}
				\label{fig:general_splice_proof_CI_1}
			\end{subfigure}
			\hfill
			\begin{subfigure}[b]{0.14\linewidth}
				\centering
				\begin{tikzpicture}[scale=0.75, thick, every node/.style={font=\footnotesize}]
					\coordinate (w) at (0,1);
					\node[circle,draw,minimum size=0.82cm,inner sep=1pt] (Li) at (0,-0.5) {$\Lambda_i$};
					\draw (w) -- (Li);
					\fill (w) circle (2pt) node[above] {$w'$};
				\end{tikzpicture}
				\caption{$(\Lambda_i,w')$}
				\label{fig:general_splice_proof_CI_2}
			\end{subfigure}
			\hfill
			\begin{subfigure}[b]{0.52\linewidth}
				\centering
				\begin{tikzpicture}[scale=0.75, thick, every node/.style={font=\footnotesize}]
					\coordinate (w)  at (0,0.9);
					\coordinate (v1) at (0,0);
					\coordinate (v2) at (2.0,-0.5);
					\coordinate (vk) at (4.2,-1.2);
					\node[circle,draw,minimum size=0.82cm,inner sep=1pt] (L1)  at (-1.1,-1.1) {$\Lambda_1$};
					\node[circle,draw,minimum size=0.82cm,inner sep=1pt] (L2)  at (1.4,-1.7)  {$\Lambda_2$};
					\node[circle,draw,minimum size=0.82cm,inner sep=1pt] (Lk1) at (3.4,-2.4)  {$\Lambda_{n-1}$};
					\node[circle,draw,minimum size=0.82cm,inner sep=1pt] (Ln)  at (5.1,-2.4)  {$\Lambda_n$};
					\draw (w) -- (v1);
					\draw (v1) -- (L1) node[pos=0.4,left=2pt] {$d_1$};
					\draw (v1) -- (v2) node[pos=0.55,above=1pt] {$d_2\cdots d_n$}
					node[pos=0.4,below=1pt] {$e_1$};
					\draw (v2) -- (L2) node[pos=0.4,left=1pt] {$d_2$};
					\draw (v2) -- (2.66,-0.71);
					\draw[dotted] (2.66,-0.71) -- (3.54,-0.99);
					\draw (3.54,-0.99) -- (vk);
					\draw (vk) -- (Lk1) node[pos=0.35,left=1pt] {$d_{n-1}$};
					\draw (vk) -- (Ln)  node[pos=0.35,right=1pt] {$d_n$};
					\fill (w)  circle (2pt) node[above] {$w$};
					\fill (v1) circle (2pt) node[left=3pt,yshift=7pt] {$v_1$};
					\fill (v2) circle (2pt) node[right=3pt,yshift=2pt] {$v_2$};
					\fill (vk) circle (2pt) node[above right=-1pt] {$v_{n-1}$};
				\end{tikzpicture}
				\caption{Splitting of the node $v_1$}
				\label{fig:general_splice_proof_CI_3}
			\end{subfigure}
			\caption{Splice diagrams used in the general case.}
			\label{fig:CI_proof}
		\end{figure}
	\end{proof}
	
	Obviously, we have the following corollary already stablished in Neumann and Wahl paper \cite{neumann2005complex}.
	\begin{corollary}
		Let  \((X,0)\) be a splice quotient singularity with \(\mathbb{Z}\)HS link and \(\Lambda\) be its splice diagram. Let \(w\) be a leaf of \(\Lambda.\) Then, \(S(\Lambda_w)\) is a complete intersection numerical semigroup.
	\end{corollary}
	
	The proof of Theorem~\ref{teo:semigroup_rooted} not only establishes that $S(\Lambda_w)$ is a complete intersection, but also allows us to track Delorme's partitions in terms of the movements of the splice diagram. Indeed, applying recursively at every node of valency bigger or equal than four of $(\Lambda,w)$  the splitting of Figure~\ref{fig:general_splice_proof_CI_3}, we obtain a splice diagram equivalent to $(\Lambda,w)$ all of whose nodes have valency three. We call it a \emph{Delorme form} of $(\Lambda,w)$. It is not unique: at each node of
	valency at least four one may choose which edge is split off first, and different choices yield different distinguished sequences for $S(\Lambda_w)$ (see Figure \ref{fig:example_Delorme form}).

	\begin{corollary}\label{cor:teosemigroup_rooted}
		Let $(\Lambda,w)$ be as in Theorem~\ref{teo:semigroup_rooted}, let
		$G$ be the family of linking numbers \(\ell_{wx}\) of $w$ with the other
		leaves $x$ of $\Lambda$, and fix a Delorme form of $(\Lambda,w)$. Then there
		is a distinguished sequence $(\Pds,\Lds)$ for $G$ satisfying the condition
		of Corollary~\ref{coro:delorme_ds_2} such that:
		\begin{enumerate}
			\item at each node $v$ of the Delorme form, with edges $u,u'$ pointing
			away from $w$, the two sets merged are the linking numbers of $w$ with
			the leaves beyond $u$ and beyond $u'$, respectively;
			\item the corresponding syzygy has degree $\ell_{wv}$.
		\end{enumerate}
		In particular, the syzygy degrees of $(\Pds,\Lds)$ are the linking numbers
		of $w$ with the nodes of the Delorme form.
	\end{corollary}
	\begin{proof}
		
		Before proceeding with the proof we present some notation and we describe a partition of the generating set of $S(\Lambda_w)$. First, number the leaves $w_0,\dots,w_g$ and denote $\ell_{i}=\ell_{ww_i},\text{ for } i=0,\dots,g.$ Then, $S(\Lambda_w)$ is generated by $\ell_{0},\dots,\ell_{g}$. Set $P_g:=\{\{\ell_i\}|0\leq i \leq g\}\text{ and }P_0:=\{\ell_0,\dots,\ell_g\}.$ Consider the node $v$ adjacent to $w$; since $v$ has valence three, there are exactly two edges $e_1,e_2$ pointing away from $w$. Set also $L_1:=\{\ell_i | w_i \text{ a leaf in } \Lambda_{ve_1}\}$, $L_2:=\{\ell_i | w_i \text{ a leaf in } \Lambda_{ve_2}\}$ and a partition $P_1:=\{L_1,L_2\}$. Repeating this process with $\Lambda_{ve_1} \text{ and } \Lambda_{ve_2}$ provides the rest of the partitions and the set $\Pds$.

		Once the partition is done, we construct $\Lds$ as follows. To each leaf other than $w$, we associate a variable $z_i$ with weight \(\ell_i\), as in Section~\ref{ssec:splice-quotient}. Let $v$, $e_1$, $e_2$, $L_1$, and $L_2$ be as above. Recall that, for every node $v$ with edges $e_1,e_2$ pointing away from $w$ and every leaf $x$ of $\Lambda_{ve_i}$, a straightforward computation from
		Definition~\ref{def:linking_number} shows that the linking numbers generating $S(\Lambda_w)$ and $S(\Lambda_{ve_i})$ satisfy the relation $\ell_{wx}=b_i\,\ell'_{vx}$, where $b_i:=\ell_{wv}/d_{ve_i}$. As a
		consequence, $S_{L_i}=b_i\,S(\Lambda_{ve_i})$ and, since $S(\Lambda_{ve_i})$
		is a numerical semigroup (see remark \ref{rem:numerical}), we have
		$b_i=\gcd(L_i)$ and $S(\Lambda_{ve_i})$ is generated by $L_i/b_i$.

		Consider an arbitrary node $v$. Let $e_1$ and $e_2$ be the two edges incident to $v$ that point away from $w$, and let $L_1$ and $L_2$ be defined as above. Set as before
		$ 
		b_i=\gcd(L_i)=\frac{\ell_{wv}}{d_{ve_i}}, \qquad i=1,2.
		$
		
		Then $S(\Lambda_{ve_i})$ is generated by the elements of $L_i/b_i$. Since $(\Lambda,w)$ satisfies the semigroup condition, we have $d_{ve_i}\in S(\Lambda_{ve_i})$, and therefore we may choose an admissible monomial $M_{v,i}$ for each $i$. As the variable $z_i$ has weight $\ell_i$, it follows that
		$$
		\deg(M_{v,i})=d_{ve_i}b_i.
		$$
		Using the identity $b_i=\ell_{wv}/d_{ve_i}$, we obtain
		\(
		\deg(M_{v,1})=\deg(M_{v,2})=\ell_{wv}=\operatorname{lcm}\bigl(\gcd(L_1),\gcd(L_2)\bigr).
		\)
		Consequently, the binomial
		$$
		Z_v:=M_{v,1}-M_{v,2}
		$$
		satisfies the hypotheses of Corollary~\ref{coro:delorme_ds_2}. Finally, setting
		$$
		\Lds:=\{Z_v \mid v \text{ is a node of }(\Lambda,w)\},
		$$
		we obtain a distinguished sequence $(\Pds,\Lds)$ satisfying the assumptions of Corollary \ref{coro:delorme_ds_2}.
	\end{proof}
	
	To illustrate the role of the Delorme form in Theorem~\ref{teo:semigroup_rooted}
	and Corollary~\ref{cor:teosemigroup_rooted}, as well as the possible occurrence of repetitions in the family of generators \(G\), we finish this subsection with the following two examples.
	
	\begin{ex}\label{example:several delorme}
		Consider the splice diagram
		$\Lambda$ of Figure~\ref{fig:example_Delorme form}, rooted at the leaf $w$, and
		denote by $w_0,w_1,w_2$ the leaves with weights $3,5,7$ at the node $v$. The
		weight $1$ points towards $w$, so it plays no role, and
		\[
		\ell_{ww_0}=5\cdot 7=35,\quad \ell_{ww_1}=3\cdot 7=21,\quad
		\ell_{ww_2}=3\cdot 5=15,\quad\text{so}\quad S(\Lambda_w)=\langle 15,21,35\rangle .
		\]
		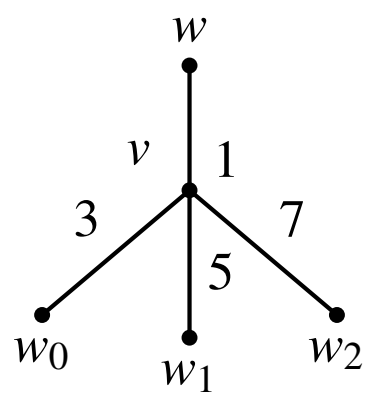
\begin{figure}[H]
			\centering
			\begin{tikzpicture}[scale=0.75, thick, every node/.style={font=\footnotesize}]
				\coordinate (w) at (0,1.1);    \coordinate (v) at (0,0);
				\coordinate (a) at (-1.3,-1.1); \coordinate (b) at (0,-1.3);
				\coordinate (c) at (1.3,-1.1);
				\draw (w) -- (v) node[pos=0.75,right] {$1$};
				\draw (v) -- (a) node[pos=0.5,above left=-2pt] {$3$};
				\draw (v) -- (b) node[pos=0.55,right=-1pt] {$5$};
				\draw (v) -- (c) node[pos=0.5,above right=-2pt] {$7$};
				\foreach \P in {w,v,a,b,c} \fill (\P) circle (2pt);
				\node[above] at (w) {$w$};  \node[left=3pt,yshift=7pt] at (v) {$v$};
				\node[below] at (a) {$w_0$}; \node[below] at (b) {$w_1$};
				\node[below] at (c) {$w_2$};
			\end{tikzpicture}
			\caption{Splice diagram with non unique Delorme form.}
			\label{fig:example_Delorme form}
		\end{figure}
		Since $v$ has three edges pointing away from $w$, Theorem~\ref{teo:splice_movements} produces the three Delorme forms of
		Figure~\ref{fig:three_Delorme_forms}: in each of them the remaining two leaves
		are attached to a new node $v'$, joined to $v$ by an edge whose weight at $v$
		is the product of their weights.
		
		\begin{figure}[H]
			\centering
			\begin{subfigure}[b]{0.3\linewidth}
				\centering
				\begin{tikzpicture}[scale=0.75, thick, every node/.style={font=\footnotesize}]
					\coordinate (w) at (0,1.1); \coordinate (v) at (0,0); \coordinate (vp) at (1,-1);
					\coordinate (a) at (-1,-1); \coordinate (b) at (0.3,-2.1); \coordinate (c) at (1.7,-2.1);
					\draw (w)--(v) node[pos=0.75,right] {$1$};
					\draw (v)--(a) node[pos=0.6,above left=-2pt] {$3$};
					\draw (v)--(vp) node[pos=0.3,above right=-2pt] {$35$} node[pos=0.8,above right=-2pt] {$3$};
					\draw (vp)--(b) node[pos=0.4,left] {$5$};
					\draw (vp)--(c) node[pos=0.4,right] {$7$};
					\foreach \P in {w,v,vp,a,b,c} \fill (\P) circle (2pt);
					\node[above] at (w) {$w$}; \node[left=3pt,yshift=7pt] at (v) {$v$}; \node[right=2pt] at (vp) {$v'$};
					\node[below] at (a) {$w_0$}; \node[below] at (b) {$w_1$}; \node[below] at (c) {$w_2$};
				\end{tikzpicture}
				\caption{}
				\label{fig:Delorme_3}
			\end{subfigure}
			\hfill
			\begin{subfigure}[b]{0.3\linewidth}
				\centering
				\begin{tikzpicture}[scale=0.75, thick, every node/.style={font=\footnotesize}]
					\coordinate (w) at (0,1.1); \coordinate (v) at (0,0); \coordinate (vp) at (1,-1);
					\coordinate (a) at (-1,-1); \coordinate (b) at (0.3,-2.1); \coordinate (c) at (1.7,-2.1);
					\draw (w)--(v) node[pos=0.75,right] {$1$};
					\draw (v)--(a) node[pos=0.6,above left=-2pt] {$5$};
					\draw (v)--(vp) node[pos=0.3,above right=-2pt] {$21$} node[pos=0.8,above right=-2pt] {$5$};
					\draw (vp)--(b) node[pos=0.4,left] {$3$};
					\draw (vp)--(c) node[pos=0.4,right] {$7$};
					\foreach \P in {w,v,vp,a,b,c} \fill (\P) circle (2pt);
					\node[above] at (w) {$w$}; \node[left=3pt,yshift=7pt] at (v) {$v$}; \node[right=2pt] at (vp) {$v'$};
					\node[below] at (a) {$w_1$}; \node[below] at (b) {$w_0$}; \node[below] at (c) {$w_2$};
				\end{tikzpicture}
				\caption{}
				\label{fig:Delorme_5}
			\end{subfigure}
			\hfill
			\begin{subfigure}[b]{0.3\linewidth}
				\centering
				\begin{tikzpicture}[scale=0.75, thick, every node/.style={font=\footnotesize}]
					\coordinate (w) at (0,1.1); \coordinate (v) at (0,0); \coordinate (vp) at (1,-1);
					\coordinate (a) at (-1,-1); \coordinate (b) at (0.3,-2.1); \coordinate (c) at (1.7,-2.1);
					\draw (w)--(v) node[pos=0.75,right] {$1$};
					\draw (v)--(a) node[pos=0.6,above left=-2pt] {$7$};
					\draw (v)--(vp) node[pos=0.3,above right=-2pt] {$15$} node[pos=0.8,above right=-2pt] {$7$};
					\draw (vp)--(b) node[pos=0.4,left] {$3$};
					\draw (vp)--(c) node[pos=0.4,right] {$5$};
					\foreach \P in {w,v,vp,a,b,c} \fill (\P) circle (2pt);
					\node[above] at (w) {$w$}; \node[left=3pt,yshift=7pt] at (v) {$v$}; \node[right=2pt] at (vp) {$v'$};
					\node[below] at (a) {$w_2$}; \node[below] at (b) {$w_0$}; \node[below] at (c) {$w_1$};
				\end{tikzpicture}
				\caption{}
				\label{fig:Delorme_7}
			\end{subfigure}
			\caption{The three Delorme forms of the diagram of
				Figure~\ref{fig:example_Delorme form}.}
			\label{fig:three_Delorme_forms}
		\end{figure}
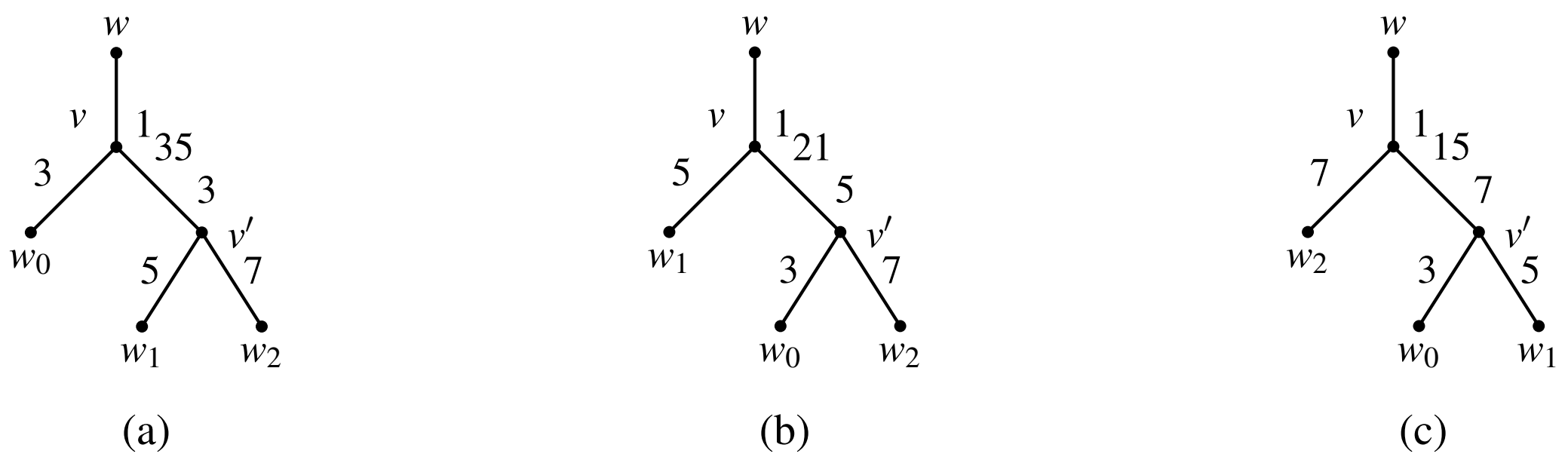
		By Corollary~\ref{cor:teosemigroup_rooted}, each Delorme form yields a
		distinguished sequence for $G=\{15,21,35\}$. Let $z_0,z_1,z_2$ be the variables
		of $w_0,w_1,w_2$, of weights $35,21,15$. In all three cases
		$\ell_{wv}=\ell_{wv'}=105$, so both syzygies have degree $105$, but the
		partitions are different:
		\begin{center}
			\begin{table}[h!]
				\centering
				\renewcommand{\arraystretch}{1.5}
				\begin{tabular}{||c c c c||}
					\hline
					Delorme form & $P_1$ & $Z_{v'}$ & $Z_v$ \\ [0.5ex]
					\hline\hline
					\ref{fig:Delorme_3} & $\{\{35\},\{21,15\}\}$ & $z_1^5-z_2^7$ & $z_0^3-z_1^5$\\
					\hline
					\ref{fig:Delorme_5} & $\{\{21\},\{35,15\}\}$ & $z_0^3-z_2^7$ & $z_1^5-z_0^3$\\
					\hline
					\ref{fig:Delorme_7} & $\{\{15\},\{35,21\}\}$ & $z_0^3-z_1^5$ & $z_2^7-z_0^3$\\
					\hline
				\end{tabular}
				\caption{Distinguished sequences attached to the Delorme forms of
					Figure~\ref{fig:three_Delorme_forms}.}
				\label{table:Delorme_forms}
			\end{table}
		\end{center}
		In this way, we see that the
		three Delorme forms correspond to the three decompositions of $S(\Lambda_w)$ as
		a gluing,
		\[
		S(\Lambda_w)=35\,\mathbb{N}+3\,\langle 5,7\rangle
		=21\,\mathbb{N}+5\,\langle 3,7\rangle
		=15\,\mathbb{N}+7\,\langle 3,5\rangle .
		\]
	\end{ex}
	\begin{ex}\label{ex:twosplicerepeated}
		Consider the non equivalent splice diagrams $\Lambda_A$ and $\Lambda_B$ of
		Figure~\ref{fig:examples}, rooted at the leaf $w$. 
		
		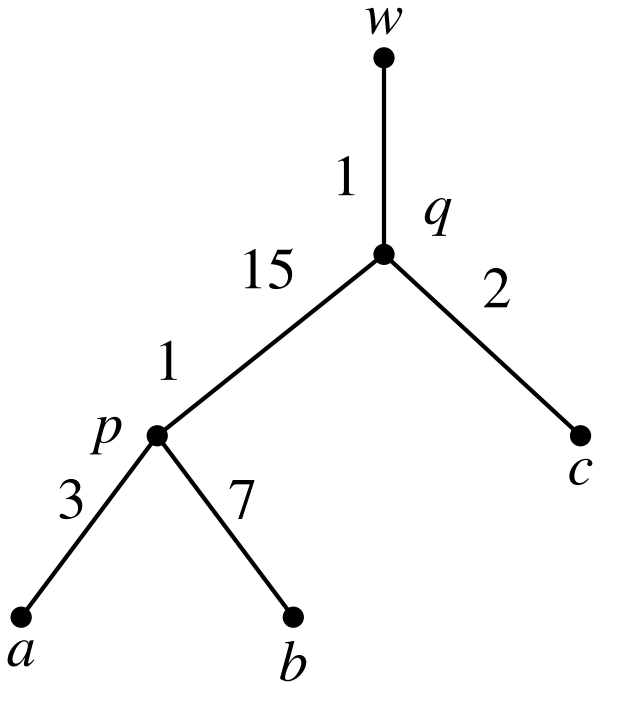
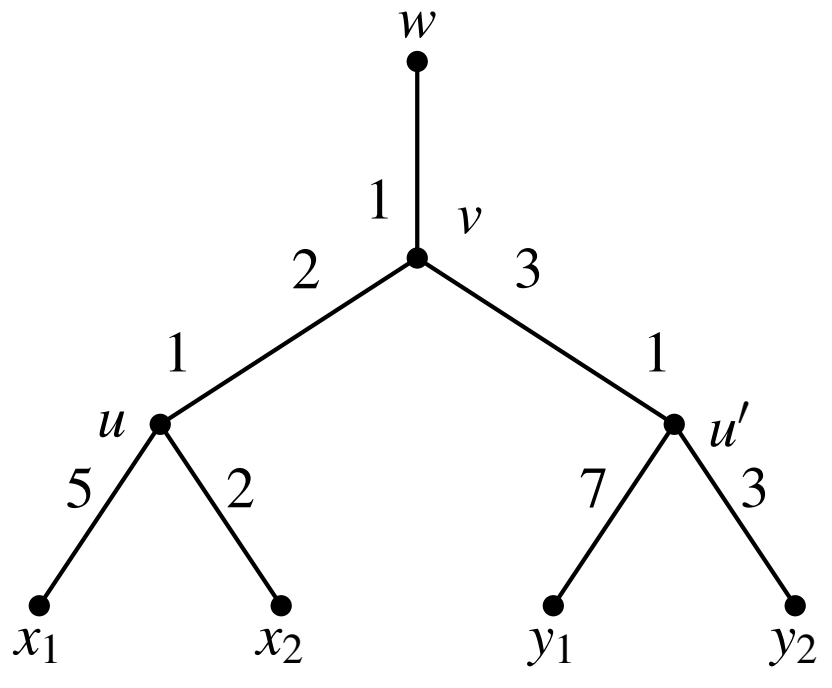
\begin{figure}
			\centering
			\begin{subfigure}[b]{0.42\linewidth}
				\centering
				\begin{tikzpicture}[thick, every node/.style={font=\small}]
					\coordinate (w) at (0,1.3);   \coordinate (q) at (0,0);
					\coordinate (p) at (-1.5,-1.2); \coordinate (c) at (1.3,-1.2);
					\coordinate (a) at (-2.4,-2.4); \coordinate (b) at (-0.6,-2.4);
					\draw (w) -- (q) node[pos=0.6,left] {$1$};
					\draw (q) -- (p) node[pos=0.3,above left=-1pt] {$15$}
					node[pos=0.8,above left=-1pt] {$1$};
					\draw (q) -- (c) node[pos=0.4,above right=-1pt] {$2$};
					\draw (p) -- (a) node[pos=0.35,left] {$3$};
					\draw (p) -- (b) node[pos=0.35,right] {$7$};
					\foreach \P in {w,q,p,c,a,b} \fill (\P) circle (2pt);
					\node[above] at (w) {$w$};    \node[right=3pt,yshift=7pt] at (q) {$q$};
					\node[left=2pt] at (p) {$p$};
					\node[below] at (a) {$a$};    \node[below] at (b) {$b$};
					\node[below] at (c) {$c$};
				\end{tikzpicture}
				\caption{$\Lambda_A$}
				\label{fig:LambdaA}
			\end{subfigure}
			\hfill
			\begin{subfigure}[b]{0.5\linewidth}
				\centering
				\begin{tikzpicture}[thick, every node/.style={font=\small}]
					\coordinate (w) at (0,1.3);    \coordinate (v) at (0,0);
					\coordinate (u) at (-1.7,-1.1); \coordinate (up) at (1.7,-1.1);
					\coordinate (x1) at (-2.5,-2.3); \coordinate (x2) at (-0.9,-2.3);
					\coordinate (y1) at (0.9,-2.3);  \coordinate (y2) at (2.5,-2.3);
					\draw (w) -- (v) node[pos=0.7,left] {$1$};
					\draw (v) -- (u)  node[pos=0.3,above left=-1pt] {$2$}
					node[pos=0.8,above left=-1pt] {$1$};
					\draw (v) -- (up) node[pos=0.3,above right=-1pt] {$3$}
					node[pos=0.8,above right=-1pt] {$1$};
					\draw (u) -- (x1)  node[pos=0.35,left]  {$5$};
					\draw (u) -- (x2)  node[pos=0.35,right] {$2$};
					\draw (up) -- (y1) node[pos=0.35,left]  {$7$};
					\draw (up) -- (y2) node[pos=0.35,right] {$3$};
					\foreach \P in {w,v,u,up,x1,x2,y1,y2} \fill (\P) circle (2pt);
					\node[above] at (w) {$w$};     \node[right=3pt,yshift=7pt] at (v) {$v$};
					\node[left=2pt] at (u) {$u$};  \node[right=2pt] at (up) {$u'$};
					\node[below] at (x1) {$x_1$};  \node[below] at (x2) {$x_2$};
					\node[below] at (y1) {$y_1$};  \node[below] at (y2) {$y_2$};
				\end{tikzpicture}
				\caption{$\Lambda_B$}
				\label{fig:LambdaB}
			\end{subfigure}
			\caption{Two splice diagrams rooted at $w$.}
			\label{fig:examples}
		\end{figure}
		
		Following Figure \ref{fig:examples}, in $\Lambda_A$, the linking numbers of $w$ with the leaves and with the nodes
		are
		\[
		\ell_{wa}=2\cdot 7=14,\quad \ell_{wb}=2\cdot 3=6,\quad \ell_{wc}=15,
		\qquad
		\ell_{wp}=2\cdot 21=42,\quad \ell_{wq}=2\cdot 15=30,
		\]
		whereas in $\Lambda_B$ they are
		\[
		\ell_{wx_1}=3\cdot 2=6,\quad \ell_{wx_2}=3\cdot 5=15,\quad
		\ell_{wy_1}=2\cdot 3=6,\quad \ell_{wy_2}=2\cdot 7=14,
		\]
		\[
		\ell_{wu}=3\cdot 10=30,\quad \ell_{wu'}=2\cdot 21=42,\quad \ell_{wv}=2\cdot 3=6.
		\]
		Hence $S((\Lambda_A)_w)=S((\Lambda_B)_w)=\langle 6,14,15\rangle$, which is a complete intersection numerical semigroup. Obviously, both rooted diagrams
		satisfy the semigroup condition. Observe that in $\Lambda_B$ the leaves $x_1$ and $y_1$ have the same linking number with $w$, namely $6=3\cdot 2=2\cdot 3$, obtained from two different paths. Thus the family of generators attached to $\Lambda_B$ is $G_B=\{6,15,6,14\}$, with a repeated element, whereas $\Lambda_A$ gives the
		minimal system $G_A=\{6,14,15\}$. By Corollary~\ref{cor:teosemigroup_rooted},
		the syzygy degrees of the corresponding distinguished sequences are the linking numbers of $w$ with the nodes, namely $30,42$ for $\Lambda_A$ and $30,42,6$ for $\Lambda_B$; the extra syzygy of degree $6$ corresponds to the binomial $z_{x_1}-z_{y_1}$ relating the two repeated generators. This is why distinguished sequences must be allowed to have repetitions.
	\end{ex}
	\subsection{Knots associated to complete intersection semigroups }\label{ssec:construction_knot}
	
	In \cite{TopologicalAlmironolivares}, the first and third named authors provided a family of knots associated to a particular class of complete intersection numerical semigroups, called free numerical semigroups. Moreover, we showed that the Alexander polynomial of each knot of the family coincides with the generating series of the semigroup. Here, we will generalize such results to the entire class of complete intersection semigroups. To do so, we will use the distinguished sequence of the complete intersection semigroup to provide a splice diagram in Delorme's form such that there is a root of the diagram with the prescribed semigroup.
	\medskip
	
	The following proposition provides a converse to Corollary~\ref{cor:teosemigroup_rooted}.
	
	\begin{theorem}\label{tma:converse_rooted}
		Let $S$ be a complete intersection numerical semigroup, let
		$G=\{a_0,\dots,a_g\}$ be a family of generators of $S$, and let
		$(\Pds,\Lds)$ be a distinguished sequence for $G$ satisfying the condition of
		Corollary~\ref{coro:delorme_ds_2}. Then there exists a rooted splice diagram
		$(\Lambda,w)$ in Delorme's form, with positive weights, pairwise coprime at
		each node and equal to $1$ on every edge pointing towards $w$, satisfying the
		semigroup condition and such that:
		\begin{enumerate}
			\item its leaves $w_0,\dots,w_g$ satisfy $\ell_{ww_j}=a_j$; in particular
			$S(\Lambda_w)=S$;
			\item its nodes $v_1,\dots,v_g$ satisfy $\ell_{wv_i}=\deg(Z_i)$;
			\item the partitions induced by $(\Lambda,w)$ as in
			Corollary~\ref{cor:teosemigroup_rooted} are those of $\Pds$.
		\end{enumerate}
	\end{theorem}
	\begin{proof}
		
		Let $S$ be a complete intersection semigroup and $G=\{a_0,\dots,a_g\}$ be a
		finite non-empty family of positive integers generating $S$ and fix a
		distinguished sequence $(\Pds,\Lds)$ for $G$ satisfying the condition of
		Corollary~\ref{coro:delorme_ds_2}. In order to lighten the notation, we
		denote $e_i:=\gcd(L_i)$, $e_i':=\gcd(L_i')$ and $r_i:=\gcd(e_i,e_i')$ for
		$L_i,L_i'\in P_i,$ and $B_i:=L_i\cup L'_i\in P_{i-1}$, so
		that $r_i=\gcd(B_i)$.
		\medskip
		
		We follow the partition sequence $P_g,P_{g-1},\dots,P_0$.
		For every $0\leq i\leq g$ and for every $L\in P_i$ we construct a rooted splice diagram $(\Lambda_L,w_L)$ whose root $w_L$ is an arrowhead with weight $1$ at its node, whose leaves $w_j$ are in bijection with the indices $j$ of $L$, and such that
		\[
		\ell_{w_Lw_j}=\frac{a_j}{\gcd(L)}\quad\text{for every } j\in L,
		\qquad\text{so that}\qquad S\bigl((\Lambda_L)_{w_L}\bigr)=S_{L/\gcd(L)} .
		\]
		For $L=\{a_j\}\in P_g$, $(\Lambda_L,w_L)$ is the diagram without nodes formed by the arrowhead $w_L$ joined to the leaf $w_j$, so that $\ell_{w_Lw_j}=1$. Following the notation of Section~\ref{ssec:complete intersections}, let
		$L_g,L_g'\in P_g$ be the sets defining $P_{g-1}$. Reordering if necessary, we may suppose without loss of generality that $L_g=\{a_0\}$ and $L_g'=\{a_1\}$. Then, to the semigroup $S_{g}:=\langle a_0/r_g,a_1/r_g\rangle$ we associate the rooted splice diagram $(\Lambda_{B_g},w_{B_g})$ of Figure~\ref{fig:Pgknot}. It has only one node $v_g$, two leaves $w_0,w_1$ and the root $w_{B_g}$, an arrowhead of weight $1$. The leaf $w_0$ carries the weight $a_1/r_g$ and the leaf
		$w_1$ the weight $a_0/r_g$, so that $\ell_{w_{B_g}w_0}=a_0/r_g$ and
		$\ell_{w_{B_g}w_1}=a_1/r_g$.
		
		\begin{figure}[H]
			\centering
			\begin{tikzpicture}[scale=0.75, thick, every node/.style={font=\footnotesize}]
				\coordinate (v) at (0,0); \coordinate (x0) at (-1.6,0); \coordinate (x1) at (0,-1.5);
				\draw (v) -- (x0) node[pos=0.35,above] {$a_1/r_g$};
				\draw (v) -- (x1) node[pos=0.4,right] {$a_0/r_g$};
				\draw[->] (v) -- (1.8,0) node[pos=0.25,above] {$1$};
				\foreach \P in {v,x0,x1} \fill (\P) circle (2pt);
				\node[below left=1pt] at (v) {$v_g$};
				\node[left] at (x0) {$w_0$}; \node[below] at (x1) {$w_1$};
				\node[right] at (1.8,0) {$w_{B_g}$};
			\end{tikzpicture}
			\caption{Rooted splice diagram $(\Lambda_{B_g},w_{B_g})$.}
			\label{fig:Pgknot}
		\end{figure}
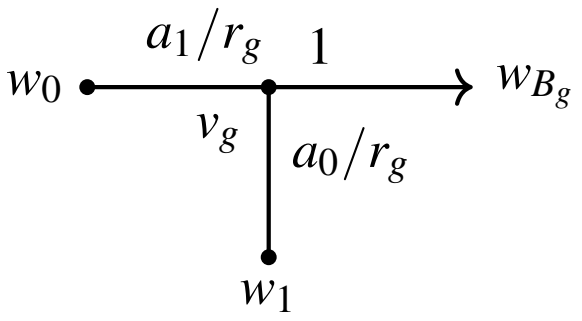

		At the step from $P_i$ to $P_{i-1}$ only the blocks $L_i,L'_i$ change: we construct $(\Lambda_{B_i},w_{B_i})$ from
		$(\Lambda_{L_i},w_{L_i})$ and $(\Lambda_{L'_i},w_{L'_i})$, and the diagrams of the remaining blocks of $P_i$ are kept unchanged in $P_{i-1}$. Assume that the rooted splice diagrams of all the blocks of
		$P_i$ have been constructed. Then, to construct $(\Lambda_{B_i},w_{B_i})$ we
		proceed according to the cardinality of the sets $L_i,L_i'\in P_i$. There are
		three different cases to be considered:
		
		\begin{itemize}
			\item $|L_i|=1=|L_i'|$:
			
			Suppose $L_i=\{a_l\}$ and $L_i'=\{a_s\}$ for some $l\neq s$. This case is
			analogous to the one for $P_g$: to the semigroup
			$S_{B_i}:=\langle a_l/r_i,a_s/r_i\rangle$ we associate the rooted splice
			diagram $(\Lambda_{B_i},w_{B_i})$ of Figure~\ref{fig:case1constr}, with one
			node $v_i$, leaves $w_l,w_s$ of weights $a_s/r_i$ and $a_l/r_i$, and the
			root $w_{B_i}$ of weight $1$. Here, it is important to notice that no splicing with
			previously constructed diagrams takes place.
			
			\begin{figure}[H]
				\centering
				\begin{subfigure}[b]{0.3\linewidth}
					\centering
					\begin{tikzpicture}[scale=0.75, thick, every node/.style={font=\footnotesize}]
						\coordinate (v) at (0,0); \coordinate (xl) at (-1.6,0); \coordinate (xs) at (0,-1.5);
						\draw (v) -- (xl) node[pos=0.35,above] {$a_s/r_i$};
						\draw (v) -- (xs) node[pos=0.4,right] {$a_l/r_i$};
						\draw[->] (v) -- (1.6,0) node[pos=0.3,above] {$1$};
						\foreach \P in {v,xl,xs} \fill (\P) circle (2pt);
						\node[below left=1pt] at (v) {$v_i$};
						\node[left] at (xl) {$w_l$}; \node[below] at (xs) {$w_s$};
						\node[right] at (1.6,0) {$w_{B_i}$};
					\end{tikzpicture}
					\caption{Case $|L_i|=|L'_i|=1$}
					\label{fig:case1constr}
				\end{subfigure}
				\hfill
				\begin{subfigure}[b]{0.3\linewidth}
					\centering
					\begin{tikzpicture}[scale=0.75, thick, every node/.style={font=\footnotesize}]
						\coordinate (v) at (0,0); \coordinate (xs) at (0,-1.5);
						\draw[->] (v) -- (-1.6,0) node[pos=0.4,above] {$a_s/r_i$};
						\draw (v) -- (xs) node[pos=0.4,right] {$e_i/r_i$};
						\draw[->] (v) -- (1.6,0) node[pos=0.3,above] {$1$};
						\foreach \P in {v,xs} \fill (\P) circle (2pt);
						\node[below left=1pt] at (v) {$v_i$};
						\node[below] at (xs) {$w_s$};
						\node[right] at (1.6,0) {$w_{B_i}$};
					\end{tikzpicture}
					\caption{Seifert link}
					\label{fig:case2constr1}
				\end{subfigure}
				\hfill
				\begin{subfigure}[b]{0.34\linewidth}
					\centering
					\begin{tikzpicture}[scale=0.75, thick, every node/.style={font=\footnotesize}]
						\coordinate (v) at (0,0); \coordinate (xs) at (0,-1.5);
						\node[circle,draw,minimum size=0.9cm,inner sep=1pt] (L) at (-2.1,0) {$\Lambda_{L_i}$};
						\draw (L) -- (v) node[pos=0.7,above] {$a_s/r_i$}
						node[pos=0.12,above] {$1$};
						\draw (v) -- (xs) node[pos=0.4,right] {$e_i/r_i$};
						\draw[->] (v) -- (1.6,0) node[pos=0.3,above] {$1$};
						\foreach \P in {v,xs} \fill (\P) circle (2pt);
						\node[below left=1pt] at (v) {$v_i$};
						\node[below] at (xs) {$w_s$};
						\node[right] at (1.6,0) {$w_{B_i}$};
					\end{tikzpicture}
					\caption{$(\Lambda_{B_i},w_{B_i})$}
					\label{fig:case2constr2}
				\end{subfigure}
				\caption{Construction for $|L'_i|=1$.}
				\label{fig:constr12}
			\end{figure}
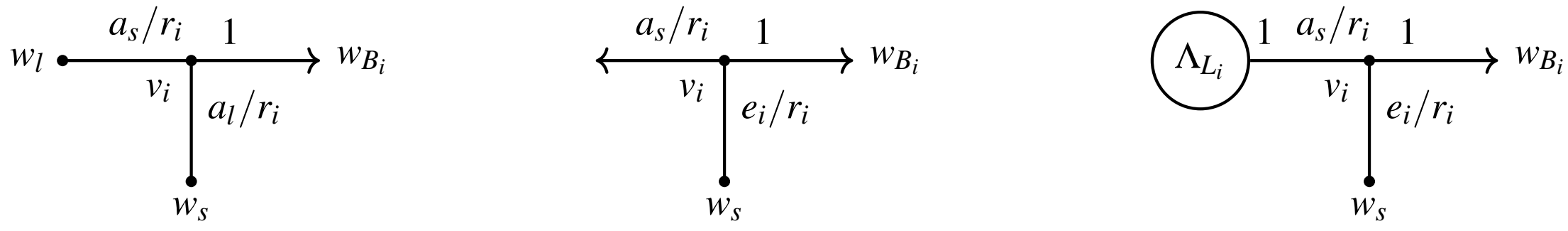
			
			\item $|L_i|>1$ and $|L_i'|=1$:
			
			Suppose $L_i'=\{a_s\}$, so that $e'_i=a_s$. Consider the Seifert link of
			Figure~\ref{fig:case2constr1}, with one node $v_i$, a leaf $w_s$ of weight
			$e_i/r_i$, an arrowhead of weight $a_s/r_i$ and an arrowhead $w_{B_i}$ of
			weight $1$. Splicing $(\Lambda_{L_i},w_{L_i})$ and this link along the root
			$w_{L_i}$ and the arrowhead of weight $a_s/r_i$, we obtain the rooted
			splice diagram $(\Lambda_{B_i},w_{B_i})$ of Figure~\ref{fig:case2constr2}.
			For every $j\in L_i$ the path from $w_{B_i}$ to $w_j$
			passes through $v_i$, where the only weight not on this path is $e_i/r_i$, and then
			follows the path from $w_{L_i}$ to $w_j$; hence
			\[
			\ell_{w_{B_i}w_j}=\frac{e_i}{r_i}\cdot\frac{a_j}{e_i}=\frac{a_j}{r_i},
			\qquad \ell_{w_{B_i}w_s}=\frac{a_s}{r_i}.
			\]
			As a consequence, we have 
			\[S\left((\Lambda_{B_i})_{w_{B_i}}\right)=\frac{e_i}{r_i}S\left((\Lambda_{L_i})_{w_{L_i}}\right)+\frac{a_s}{r_i}\mathbb{N}\]
			\item $|L_i|>1$ and $|L_i'|>1$:
			
			Consider the Seifert link of Figure~\ref{fig:case3constr1}, with one node
			$v_i$ and three arrowheads of weights $e'_i/r_i$, $e_i/r_i$ and $1$.
			Splicing $(\Lambda_{L_i},w_{L_i})$ along $w_{L_i}$ and the arrowhead of
			weight $e'_i/r_i$, and $(\Lambda_{L'_i},w_{L'_i})$ along $w_{L'_i}$ and the
			arrowhead of weight $e_i/r_i$, we obtain the rooted splice diagram
			$(\Lambda_{B_i},w_{B_i})$ of Figure~\ref{fig:case3constr2}, whose root
			$w_{B_i}$ is the remaining arrowhead. As in the previous
			case,
			\[
			\ell_{w_{B_i}w_j}=\frac{e_i}{r_i}\cdot\frac{a_j}{e_i}=\frac{a_j}{r_i}
			\ \ (j\in L_i),\qquad
			\ell_{w_{B_i}w_j}=\frac{e'_i}{r_i}\cdot\frac{a_j}{e'_i}=\frac{a_j}{r_i}
			\ \ (j\in L'_i)
			\]
			and again we obtain the gluing operation
			\[S\left((\Lambda_{B_i})_{w_{B_i}}\right)=\frac{e_i}{r_i}S\left((\Lambda_{L_i})_{w_{L_i}}\right)+\frac{e'_i}{r_i}S\left((\Lambda_{L'_i})_{w_{L'_i}}\right)\]
		\end{itemize}
		
		\begin{figure}[H]
			\centering
			\begin{subfigure}[b]{0.4\linewidth}
				\centering
				\begin{tikzpicture}[scale=0.75, thick, every node/.style={font=\footnotesize}]
					\coordinate (v) at (0,0);
					\draw[->] (v) -- (-1.4,-1.3) node[pos=0.45,above left=-2pt] {$e'_i/r_i$};
					\draw[->] (v) -- (1.4,-1.3)  node[pos=0.45,above right=-2pt] {$e_i/r_i$};
					\draw[->] (v) -- (0,1.3)     node[pos=0.35,right] {$1$};
					\fill (v) circle (2pt);
					\node[below=5pt] at (v) {$v_i$};
					\node[above] at (0,1.3) {$w_{M_i}$};
				\end{tikzpicture}
				\caption{Seifert link}
				\label{fig:case3constr1}
			\end{subfigure}
			\hfill
			\begin{subfigure}[b]{0.45\linewidth}
				\centering
				\begin{tikzpicture}[scale=0.75, thick, every node/.style={font=\footnotesize}]
					\coordinate (v) at (0,0);
					\node[circle,draw,minimum size=0.9cm,inner sep=1pt] (L)  at (-1.7,-1.6) {$\Lambda_{L_i}$};
					\node[circle,draw,minimum size=0.9cm,inner sep=1pt] (Lp) at (1.7,-1.6)  {$\Lambda_{L'_i}$};
					\draw (v) -- (L)  node[pos=0.3,above left=-2pt] {$e'_i/r_i$}
					node[pos=0.8,above left=-2pt] {$1$};
					\draw (v) -- (Lp) node[pos=0.3,above right=-2pt] {$e_i/r_i$}
					node[pos=0.8,above right=-2pt] {$1$};
					\draw[->] (v) -- (0,1.3) node[pos=0.35,right] {$1$};
					\fill (v) circle (2pt);
					\node[below=5pt] at (v) {$v_i$};
					\node[above] at (0,1.3) {$w_{M_i}$};
				\end{tikzpicture}
				\caption{$(\Lambda_{M_i},w_{M_i})$}
				\label{fig:case3constr2}
			\end{subfigure}
			\caption{Construction for $|L_i|,|L'_i|>1$.}
			\label{fig:constr3}
		\end{figure}
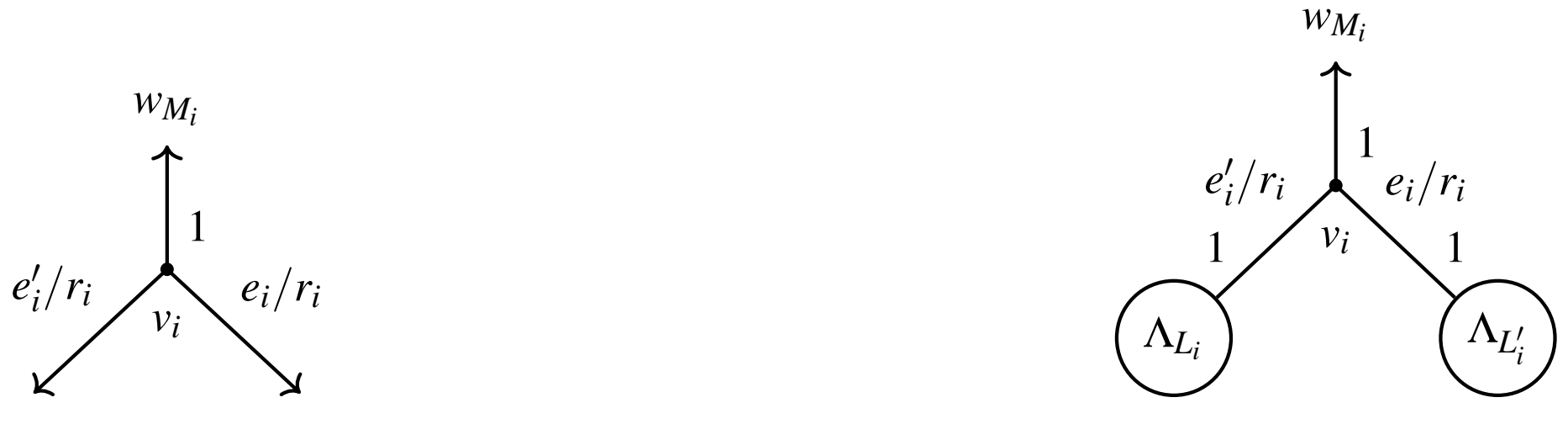
		
		In the three cases, the new node $v_i$ has weight
		$e'_i/r_i$ on the edge towards $L_i$, weight $e_i/r_i$ on the edge towards
		$L'_i$, and weight $1$ towards the new root; when $L_i=\{a_l\}$ (resp.
		$L'_i=\{a_s\}$) the corresponding edge ends at the leaf $w_l$ (resp. $w_s$).
		These weights are positive and pairwise coprime, and the linking number of
		the new root with the new node is
		$\ell_{w_{B_i}v_i}=e_ie'_i/r_i^2=\deg(Z_i)/r_i$. Moreover, write
		$Z_i=\mathbf{x}-\mathbf{y}$ with $\mathbf{x}\in M(S_{L_i})$; since
		$\deg(\mathbf{x})=\operatorname{lcm}(e_i,e'_i)=e_ie'_i/r_i$, then $e'_i/r_i\in S_{L_i/e_i}=S\bigl((\Lambda_{L_i})_{w_{L_i}}\bigr)$,
		and similarly $e_i/r_i\in S_{L'_i/e'_i}$. Hence, every
		$(\Lambda_{B_i},w_{B_i})$ satisfies the semigroup condition. Finally, since
		$P_0=\{G\}$ and $\gcd(G)=1$, the rooted splice diagram
		$(\Lambda,w):=(\Lambda_{B_1},w_{B_1})$ is in Delorme's form, satisfies
		$S(\Lambda_w)=S$, and the partitions it induces are those of $\Pds$. 
		
	\end{proof}
	Moreover, observe that the splice sequence obtained from the distinguished sequence \((\Pds,\Lds)\) provides a complete description of the gluing decomposition of the complete intersection numerical semigroup \(S\) and viceversa. Therefore, the previous construction provides a topological realization of the combinatorial properties of a complete intersection numerical semigroup. 
	
	By construction, all the rooted splice diagrams constructed by means of Theorem \ref{tma:converse_rooted} give rise to embedded knots in an integral homology sphere. Moreover, by \cite[Theorem 11.2]{EN} these embedded knots are fibered. Thus, for every complete intersection numerical semigroup we can provide the following topological families associted to it.
	
	\begin{defin}
		Let $S$ be a complete intersection semigroup. The family of knots associated with $S$ is defined as
		\begin{equation*}
			\mathcal{K}_S:=\left\{ K \hspace{1mm}\mid \hspace{1mm} \begin{array}{c}
				K \text{ is the embedded knot associated to  a distinguished sequence}  \left( \Pds, \Lds \right )\\
				\text{ for }S  \text{  that fulfills the condition from Corollary \ref{coro:delorme_ds_2}. }   
			\end{array}  \right\}
		\end{equation*}
		
		The family of splice diagrams associated with \(S\) is defined as 
		\begin{equation*}
			\mathcal{T}_S:=\left\{(\Lambda,w) \hspace{1mm}\mid \hspace{1mm} \begin{array}{c}
				(\Lambda,w) \text{ is the rooted splice associated to  a distinguished sequence}  \left( \Pds, \Lds \right )\\
				\text{ provided by Theorem \ref{tma:converse_rooted} }
			\end{array}  \right\}
		\end{equation*}
	\end{defin}
	These families generalize the families of a free numerical semigroup provided in \cite{TopologicalAlmironolivares}. However, the construction given here is more general, even for free semigroups: the construction of \cite{TopologicalAlmironolivares} only uses the free condition, whereas here we use the full Delorme decomposition and allow repetitions in the family of
	generators. As a consequence, $\mathcal{K}_S$ may contain knots which are not iterated torus knots. For instance, the semigroup of Example~\ref{ex:twosplicerepeated} is free, and the knot associated with $\Lambda_B$ is not an iterated torus knot, since its splice diagram is
	minimal and is not a chain.

	As mentioned at the beginning of this subsection, the main goal of this subsection is to associate a family of fibered knots to a given complete intersection semigroup such that the Alexander polynomial of each knot of the family coincides with the generating series of the semigroup modulo a \((1-t)\) factor; thus generalizing the results of \cite{TopologicalAlmironolivares}. Let us now move to show that any knot in the family has the same Alexander polynomial and it coincides with the generating series of the semigroup \(S.\)
	
	\begin{theorem}\label{teo:alexpoly_generating}
		Let $S$ be a complete intersection semigroup. Fix $G:=\{a_0,\dots, a_g\}$ a family of generators and consider $\left(\mathscr{P}\hspace{1mm}, \hspace{1mm} \mathscr{L}\right)$ a distinguished sequence for $G$ satisfying the condition from Corollary \ref{coro:delorme_ds_2}. Denote $K$ the knot associated to $\left(\mathscr{P}\hspace{1mm}, \hspace{1mm} \mathscr{L}\right)$. Then, the Alexander polynomial for $K$ is 
		
		\begin{equation}\label{eq:alexpoly_construction}
			\Delta_K\left(t\right)=\frac{\left(t-1\right)\prod_{i=1}^{g}\left(t^{m_i}-1\right)}{\prod_{i=0}^{g}\left(t^{a_i}-1\right)}=(1-t)\left(\sum_{n\in S}t^n\right),
		\end{equation}
		
		where the $m_i$ denote the degrees of the syzygies of the distinguished sequence. 
	\end{theorem}
	
	\begin{proof}
		Let $G=\{a_0,\dots,a_g\}$ be a set of generators for $S$ and $m_1,\dots,m_{g}$ the syzygies from a distinguished sequence. Then, it is a straightforward consequence of \cite[Theorem 3.1]{St78hilbert} (see also \cite[Section 4.1]{AlmironCI}) that
		\begin{center}
			$\sum_{n\in S}t^n=\frac{\prod_{i=1}^{g}(1-t^{m_i})}{\prod_{i=0}^{g}(1-t^{a_i})}.$
		\end{center}
		On the other hand, according to Theorem \ref{teo:alexpol_splice}, we only need to compute the linking numbers of the arrow with the vertices of the splice diagram. By Theorem \ref{tma:converse_rooted}, we have that the linking numbers with the nodes agree with the syzygy degrees and the linking numbers with the leaves agree with the generators of the semigroup. In addition, the splice diagram we are considering is in Delorme's form (each node has valence exactly $3$). Hence, Equation \ref{eq:alexpoly_construction} translates as:
		\[\Delta_K\left(t\right)=\frac{\left(t-1\right)\prod_{i=1}^{g}\left(t^{m_i}-1\right)}{\prod_{i=0}^{g}\left(t^{a_i}-1\right)}\]
		from where the result follows.
	\end{proof}
	
	To finish, let us show with an example that, generally, the splice diagrams associated to the complete intersection numerical semigroup are not associated to any surface singularity. 
	
	\begin{ex}\label{ejemplo:negativedeterminant}
		Take $S$ the semigroup generated by the set $G:=\{10,14,15,21\}$. We assert that $S$ is a complete intersection. In particular, consider the partitions $\mathcal{P}:=\{P_0,P_1,P_2,P_3\}$ with
		\begin{center}
			$
			\begin{array}{lll}
				P_3:=\{\{10\},\{14\},\{15\},\{21\}\}, & L_3:=\{10\},& \hspace{2.5mm} L_3':=\{15\}  \\
				
				P_2:=\{\{10,15\}, \{14\},\{21\}\}, & L_2:=\{14\},& \hspace{2.5mm} L_2':=\{21\}\\
				
				P_1:=\{\{10,15\},\{14,21\}\}, & L_1:=\{10,15\},& \hspace{2.5mm} L_1':=\{14,21\}\\
				
				P_0:=\{10,14,15,21\},& &\\
			\end{array}
			$
		\end{center}
		
		and the binomials $\mathcal{Z}:=\{Z_1,Z_2,Z_3\}$ with
		
		\begin{center}
			$
			\begin{array}{lll}
				Z_3:= x_0^3-x_2^2,  &  Z_2:=x_1^3-x_3^2 , &  Z_1:= x_0^2x_2-x_1x_3.\\

			\end{array}
			$
		\end{center}
		Clearly the distinguished sequence $\left( \Pds, \Lds \right )$ satisfies Corollary \ref{coro:delorme_ds_2}. Thus, we can construct its associated knot.\\
		
		Note that the greatest common divisors are $r_3=e_1=5,\hspace{1mm} r_2=e_1'=7\text{ and } r_1=1.$ Hence, the preeceding construction yields the following elements: for $\{10,15\}\text{ and }\{14,21\}$ the $\left(2,3\right)-$ torus knot, see Figure \ref{fig:exone}. Then, since $P_0$ has $|L_1|=|L_1'|=2$, we splice the diagrams for $L_1
		\text{ and }L_1'$ with the link represented in Figure \ref{fig:extwo} along the arrows. Consequently, the splice diagram for the knot associated to $S$ via this distinguished sequence is the one ilustrated in Figure \ref{fig:exthree}
		
		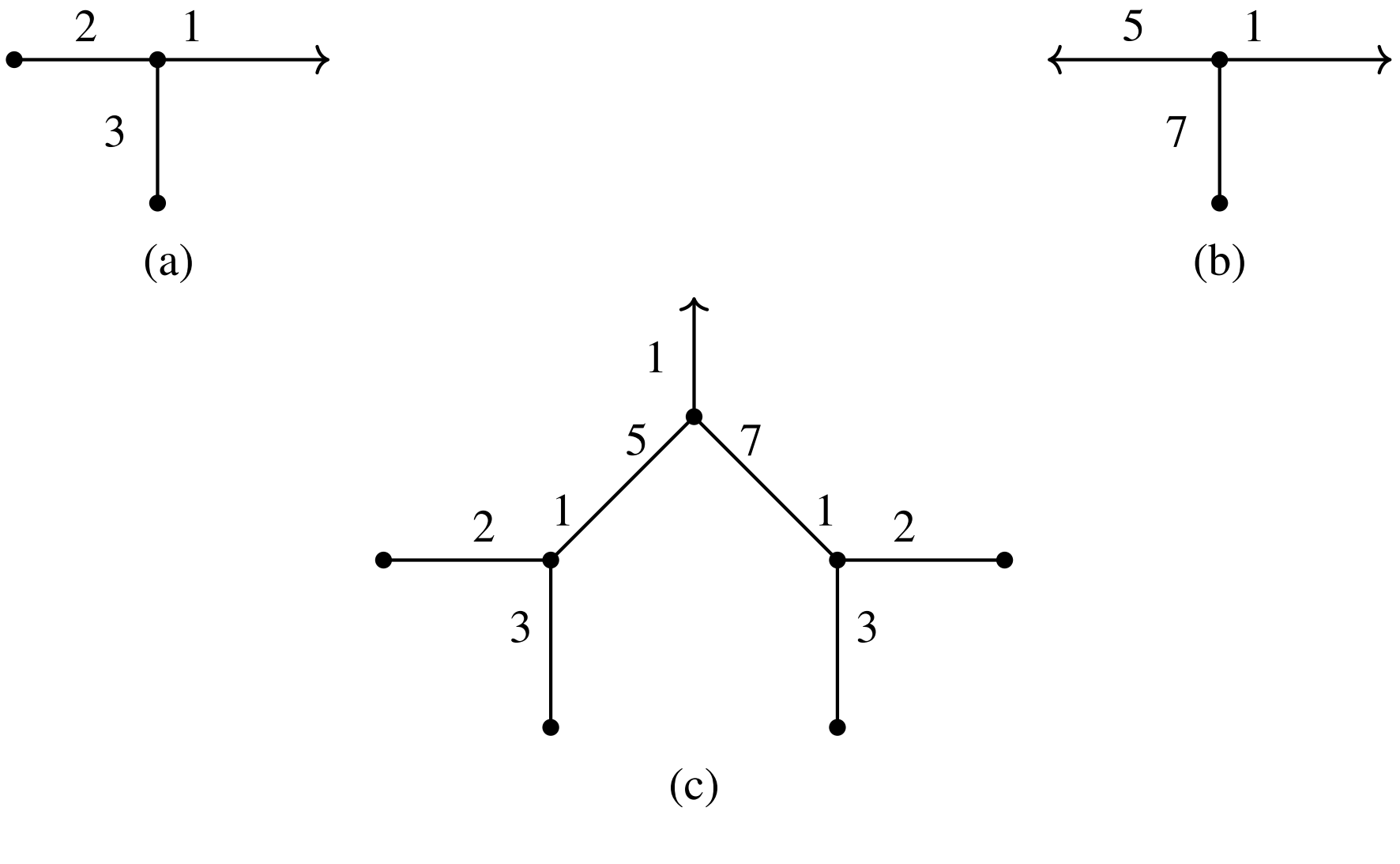
\begin{figure}[H]
			\centering
			\begin{subfigure}{0.45\textwidth}
				\centering
				\begin{tikzpicture}[scale=0.8, thick]
					
					\fill (-1.5, 0) circle (2.5pt); 
					\fill (0, 0) circle (2.5pt);    
					\fill (0, -1.5) circle (2.5pt); 
					
					\draw (-1.5, 0) -- (0, 0) node[midway, above=0pt] {\small$2$};
					
					\draw (0, 0) -- (0, -1.5) node[midway, left=3pt] {\small$3$};
					
					\draw[->] (0, 0) -- (1.8, 0) node[pos=0.2, above=0pt] {\small$1$};
					
				\end{tikzpicture}
				\subcaption{}
				\label{fig:exone}
			\end{subfigure}
			\hfill
			\begin{subfigure}{0.45\textwidth}
				\centering
				\begin{tikzpicture}[scale=0.8, thick]
					
					\fill (0, 0) circle (2.5pt);    
					\fill (0, -1.5) circle (2.5pt); 
					
					\draw[<-] (-1.8, 0) -- (0, 0) node[midway, above=0pt] {\small$5$};
					
					\draw (0, 0) -- (0, -1.5) node[midway, left=3pt] {\small$7$};
					
					\draw[->] (0, 0) -- (1.8, 0) node[pos=0.2, above=0pt] {\small$1$};
					
				\end{tikzpicture}
				\subcaption{}
				\label{fig:extwo}
			\end{subfigure}
			\begin{subfigure}{0.4\textwidth}
				\centering
				\begin{tikzpicture}[scale=1, thick]
					
					\coordinate (Top) at (0, 1.2);
					
					\coordinate (L1) at (-1.2, 0);
					\coordinate (R1) at (1.2, 0);
					
					\coordinate (L1_L) at (-2.6, 0);
					\coordinate (L1_B) at (-1.2, -1.4);
					
					\coordinate (R1_R) at (2.6, 0);
					\coordinate (R1_B) at (1.2, -1.4);
					
					\draw[->] (Top) -- (0, 2.2) node[midway, left=2pt] {\small$1$};
					
					\draw (Top) -- (L1) node[pos=0.4, above =0pt] {\small$5$};
					\draw (Top) -- (R1) node[pos=0.4, above =0pt] {\small$7$};
					
					\draw (L1) -- (L1_L) node[pos=0.4, above =0pt] {\small$2$};
					\draw (L1) -- (L1_B) node[pos=0.4, left =0pt] {\small$3$};
					\node[above =1pt] at (-1.1, 0.1) {\small$1$}; 
					
					\draw (R1) -- (R1_R) node[pos=0.4, above =0pt] {\small$2$};
					\draw (R1) -- (R1_B) node[pos=0.4, right =0pt] {\small$3$};
					\node[above=1pt] at (1.1, 0.1) {\small$1$}; 
					
					\fill (Top)  circle (2pt);
					\fill (L1)   circle (2pt);
					\fill (R1)   circle (2pt);
					\fill (L1_L) circle (2pt);
					\fill (L1_B) circle (2pt);
					\fill (R1_R) circle (2pt);
					\fill (R1_B) circle (2pt);
					
				\end{tikzpicture}
				\subcaption{}
				\label{fig:exthree}
			\end{subfigure}
			\caption{Splice diagrams associated to the partition $\mathcal{P}$.}
			\label{fig:placeholder}
		\end{figure}
		
		Check that $M_2=\{14,21\}\text{ and } M_3=\{10,15\}$. Thus, the semigroup generated by $G_2 \text{ and }G_3$ is $S_2=S_3:=\langle 2,3\rangle$. The induced distinguished sequence for each $M_i$, $i=2,3$, is:\\
		\begin{itemize}
			\item $\Pds_i:=\left( P_1:=\{\{2\},\{3\}\}, P_0:=\{2,3\}\right)$\\
			
			\item $\Lds_i:=\left(x_0^3-x_1^2\right)$\\
		\end{itemize}
		
		and these distinguished sequences associates the knot in Figure \ref{fig:exone} to $S_2$ and $S_3$ . Finally, it easy to see that the edge determinant of the two central edges 
		\[1\cdot 5-2\cdot3\cdot7\cdot 1=-37 \quad \text{and}\quad 1\cdot 7-2\cdot3\cdot5\cdot 1=-23\]
		are not positive. Thus, by Theorem \ref{teo:singularity_cond} the splice diagram (C) is not associated to a link of a complex surface singularity. 
	\end{ex}

	\subsection{End-curves with prescribed semigroup}\label{subsec:negativedefinite}
	
	At the moment, we have proven on one side that the semigroup of a rooted splice diagram is a complete intersection numerical semigroup. In particular, as splice quotient singularities admit end-curves, the semigroup of an end-curve in a splice quotient singularity is always a complete intersection numerical semigroup. On the other side, we have also proven that any complete intersection numerical semigroup can be realized as the semigroup of a rooted splice diagram.  However, as Example \ref{ejemplo:negativedeterminant} shows, the splice diagram constructed in subsection \ref{ssec:construction_knot} with the prescribed semigroup for a certain root is not, in general, the splice diagram associated with the link of a complex surface singularity. Consequently, it need not represent the semigroup of an end-curve function. However, we can slightly modify some of the edge weights in such a way the main properties of the splice diagram are the same but now the semigroup being realized as the value semigroup of an end curve. First, we prove the following crucial theorem.

	\begin{theorem}\label{thm:alteration}
		Let $S$ be a complete intersection numerical semigroup, let $G$ be a family of
		generators of $S$ and let
		$(\Lambda,w)$ be the rooted splice diagram in Delorme's form associated to it by
		Theorem~\ref{tma:converse_rooted}. Then there exists a rooted splice diagram
		$(\widetilde{\Lambda},w)$ with the same underlying tree and the same weights as
		$(\Lambda,w)$, except for the weights pointing towards $w$ at the nodes not
		adjacent to $w$, such that:
		\begin{enumerate}
			\item its weights are positive and pairwise coprime at each node,
			\item all its edge determinants are positive,
			\item it satisfies the semigroup condition at every edge,
			\item the linking numbers of $w$ with all the vertices of
			$\widetilde{\Lambda}$ coincide with those in $\Lambda$; in particular
			$S(\widetilde{\Lambda}_w)=S$.
		\end{enumerate}
	\end{theorem}
	\begin{proof}
		Assume \(G=\{a_0,\dots,a_g\}\) is a family of generators of \(S\), let \((\Pds,\Lds)\) be a distinguished sequence for \(G\) satisfying the condition of Corollary~\ref{coro:delorme_ds_2}, and let \((\Lambda,w)\) be the rooted splice diagram in Delorme's form associated to \((\Pds,\Lds)\) by Theorem~\ref{tma:converse_rooted}. Thus \((\Lambda,w)\) has \(g\) nodes \(v_1,\dots,v_g\). At each of the nodes \(v_2,\dots,v_g\) we will substitute the weight \(1\) pointing towards \(w\) by a new weight, obtained recursively
		from a sequence of inequalities starting at the last step of the construction in the proof of Theorem~\ref{tma:converse_rooted}. In this way, we will make extensive use of the notation introduced in that proof. 
		
		First, we consider the node \(v_1\) created when passing from \(P_1\) to \(P_0.\) Observe that if \(|L_1|=|L'_1|=1\) then there is nothing to proof as there is no other edges with weight \(1\) and the semigroup have only two generators. Then let us assume \(|L'_1|>1.\) By construction \(v_1\), it has three edges: the edge
		towards the arrow $w$, whose weight $1$ is not modified, and two edges going down to $\Lambda_{L_1}$ and $\Lambda_{L'_1}$, with weights $d_1=e'_1/r_1$ and $d_2=e_1/r_1$, respectively. Recall that \(e_1=\gcd(L_1),\) \(e'_1=\gcd(L'_1)\) and \(r_1=\gcd(e_1,e'_1).\) There two cases to be considered: 
		\begin{itemize}
			\item Assume \(|L_1|=1\), say \(L_1=\{a_s\}\). Then \(e_1=a_s\), the edge of weight \(d_1=e'_1/r_1\) ends at the leaf \(w_s\), and the edge of weight \(d_2=a_s/r_1\) goes to the node \(v_2\) of \(\Lambda_{L'_1}\) adjacent to its root. The node
			\(v_2\) has weight \(\alpha_2=1\) on the edge pointing towards \(w\), which we are going to re-scale, and weights \(d'_1,d'_2\) on the edges going down (see Figure \ref{fig:rescale_base1}). The edge determinant of \(v_1,v_2\) is
			\[
			\alpha_2\,\frac{a_s}{r_1}-\frac{e'_1}{r_1}\,d'_1d'_2 .
			\]
			At this step, we choose \(\widetilde{\alpha}_2\equiv 1 \pmod{d'_1d'_2}\) such
			that \(\widetilde{\alpha}_2> \frac{e'_1}{a_s}\,d'_1d'_2\). Then the weights at \(v_2\) are pairwise coprime and, replacing \(\alpha_2\) by \(\widetilde{\alpha}_2\), the edge determinant is positive. Moreover, the semigroup condition at \(v_2\) in the direction of \(w\) is trivially satisfied, since the corresponding semigroup is \(\mathbb{N}\). Since \(\widetilde{\alpha}_2\) lies on the path from \(w\) to every vertex, and on the path from any node to the leaves lying away from \(w\), neither the linking numbers of \(w\), and hence \(S(\Lambda_w)\), nor the semigroup condition at the edges pointing away from \(w\) are affected. The edges of the nodes below \(v_2\) pointing towards \(w\) will be considered in the next steps.
			\item Assume \(|L_1|>1\) and \(|L'_1|>1\), and denote by \(N_1\) the number of nodes of \(\Lambda_{L_1}\). We label the nodes of \(\Lambda_{L_1}\) as \(v_2,\dots,v_{N_1+1}\) and the nodes of \(\Lambda_{L'_1}\) as \(v_{N_1+2},\dots,v_g\). We only need to elaborate on the nodes of \(\Lambda_{L_1}\), as the situation is completely symmetric for \(\Lambda_{L'_1}\) and one does not affect the other. Then, as we did before, we consider the node \(v_2\), which has weight \(\alpha_2=1\) on the edge pointing towards \(w\), which we are going to re-scale, and weights \(d'_1,d'_2\) on the edges going down (see Figure \ref{fig:rescale_base2}). At this step, we choose \(\widetilde{\alpha}_2\equiv 1 \pmod{d'_1d'_2}\) such that
			\[
			\widetilde{\alpha}_2>\max\Bigl\{c\bigl(S((\Lambda_{L'_1})_{w_{L'_1}})\bigr),\,
			\frac{e_1}{e'_1}\,d'_1d'_2\Bigr\}.
			\]
			This condition ensures that the edge determinant \(\widetilde{\alpha}_2\,\frac{e'_1}{r_1}-\frac{e_1}{r_1}\,d'_1d'_2\) of \([v_1,v_2]\) is positive, and that \(\widetilde{\alpha}_2\in S\bigl((\Lambda_{L'_1})_{w_{L'_1}}\bigr)\), thus satisfying the semigroup condition. As in the previous case, the weights at \(v_2\) are pairwise coprime, and
			neither \(S(\Lambda_w)\) nor the semigroup condition at the edges pointing
			away from \(w\) are affected.
		\end{itemize}
		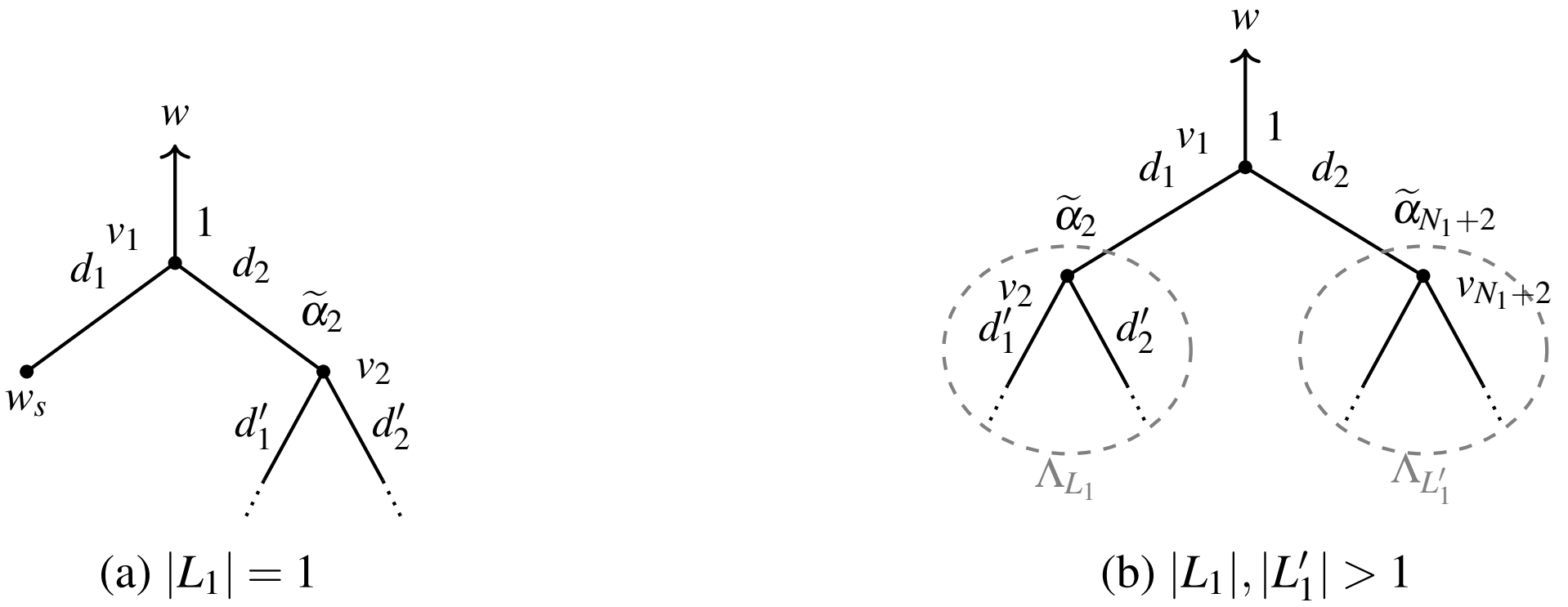
\begin{figure}[H]
			\centering
			\begin{subfigure}[b]{0.42\linewidth}
				\centering
				\begin{tikzpicture}[scale=0.8, thick, every node/.style={font=\footnotesize}]
					\coordinate (v1) at (0,0); \coordinate (ws) at (-1.5,-1.1);
					\coordinate (v2) at (1.5,-1.1);
					\draw[->] (v1) -- (0,1.2) node[pos=0.35,right] {$1$};
					\draw (v1) -- (ws) node[pos=0.35,above left=-2pt] {$d_1$};
					\draw (v1) -- (v2) node[pos=0.3,above right=-2pt] {$d_2$}
					node[pos=0.75,above right=-2pt] {$\widetilde{\alpha}_2$};
					\draw (v2) -- (0.9,-2.2) node[pos=0.5,left] {$d'_1$};
					\draw (v2) -- (2.1,-2.2) node[pos=0.5,right] {$d'_2$};
					\draw[dotted] (0.9,-2.2) -- (0.7,-2.6); \draw[dotted] (2.1,-2.2) -- (2.3,-2.6);
					\foreach \P in {v1,ws,v2} \fill (\P) circle (2pt);
					\node[above] at (0,1.2) {$w$};
					\node[left=3pt,yshift=6pt] at (v1) {$v_1$};
					\node[below] at (ws) {$w_s$};
					\node[right=3pt] at (v2) {$v_2$};
				\end{tikzpicture}
				\caption{$|L_1|=1$}
				\label{fig:rescale_base1}
			\end{subfigure}
			\hfill
			\begin{subfigure}[b]{0.52\linewidth}
				\centering
				\begin{tikzpicture}[scale=0.8, thick, every node/.style={font=\footnotesize}]
					\coordinate (v1) at (0,0); \coordinate (u) at (-1.8,-1.1);
					\coordinate (up) at (1.8,-1.1);
					\draw[->] (v1) -- (0,1.2) node[pos=0.35,right] {$1$};
					\draw (v1) -- (u) node[pos=0.3,above left=-2pt] {$d_1$}
					node[pos=0.75,above left=-2pt] {$\widetilde{\alpha}_2$};
					\draw (v1) -- (up) node[pos=0.3,above right=-2pt] {$d_2$}
					node[pos=0.75,above right=-2pt] {$\widetilde{\alpha}_{N_1+2}$};
					\draw (u) -- (-2.4,-2.2) node[pos=0.5,left] {$d'_1$};
					\draw (u) -- (-1.2,-2.2) node[pos=0.5,right] {$d'_2$};
					\draw (up) -- (1.2,-2.2); \draw (up) -- (2.4,-2.2);
					\draw[dotted] (-2.4,-2.2) -- (-2.6,-2.6); \draw[dotted] (-1.2,-2.2) -- (-1.0,-2.6);
					\draw[dotted] (1.2,-2.2) -- (1.0,-2.6); \draw[dotted] (2.4,-2.2) -- (2.6,-2.6);
					\draw[dashed,gray] (-1.8,-1.85) ellipse (1.25 and 1.05);
					\draw[dashed,gray] (1.8,-1.85) ellipse (1.25 and 1.05);
					\node[gray] at (-1.8,-3.15) {$\Lambda_{L_1}$};
					\node[gray] at (1.8,-3.15) {$\Lambda_{L'_1}$};
					\foreach \P in {v1,u,up} \fill (\P) circle (2pt);
					\node[above] at (0,1.2) {$w$};
					\node[left=3pt,yshift=6pt] at (v1) {$v_1$};
					\node[left=3pt,yshift=-4pt] at (u) {$v_2$};
					\node[right=3pt,yshift=-4pt] at (up) {$v_{N_1+2}$};
				\end{tikzpicture}
				\caption{$|L_1|,|L'_1|>1$}
				\label{fig:rescale_base2}
			\end{subfigure}
			\caption{The node \(v_1\) and the weights re-scaled in the base case.}
			\label{fig:rescale_base}
		\end{figure}
		Now, we proceed recursively. For a block \(L\) of the partition sequence with
		\(|L|>1\), we denote by \(u_L\) the node of \(\Lambda_L\) adjacent to its root,
		and by \(\overline{S}_L\) the semigroup generated by the reduced linking
		numbers of \(u_L\) with the leaves \(w_l\), \(a_l\in G\setminus L\), and with
		the arrow \(w\); that is, with the leaves seen from \(u_L\) in the direction of
		\(w\). Recall also that, by Theorem~\ref{tma:converse_rooted},
		\(S\bigl((\Lambda_L)_{w_L}\bigr)=S_{L/\gcd(L)}\). 
		
		Every node different from \(v_1\) is the node \(u_{L_i}\) or \(u_{L'_i}\) of one of the two blocks merged at some step \(i\), and the node adjacent to it in the direction of \(w\) is the node \(u_{L_i\cup L'_i}\) created at that step. Hence, interchanging the roles of \(L_i\) and \(L'_i\) if necessary, the nodes can be labelled in such a way that \(v_{j+1}=u_{L_i}\) for some \(i\), and the node adjacent to it in the direction of \(w\) is \(v_k=u_{L_i\cup L'_i}\) with \(k\leq j\).
		
		\begin{figure}[H]
			\centering
			\begin{tikzpicture}[scale=0.85, thick, every node/.style={font=\footnotesize}]
				\coordinate (vk) at (0,0); \coordinate (vj) at (1.9,-1.3);
				\draw[->,dashed] (vk) -- (0,1.4) node[pos=0.3,right] {$\widetilde{\alpha}_k$};
				\node[above] at (0,1.4) {towards $w$};
				\node[circle,draw,minimum size=0.9cm,inner sep=1pt] (Lp) at (-1.9,-1.5) {$\Lambda_{L'_i}$};
				\draw (vk) -- (Lp) node[pos=0.35,above left=-2pt] {$d_1$};
				\draw (vk) -- (vj) node[pos=0.3,above right=-2pt] {$d_0$}
				node[pos=0.78,above right=-2pt] {$\widetilde{\alpha}_{j+1}$};
				\draw (vj) -- (1.2,-2.5) node[pos=0.5,left] {$d'_1$};
				\draw (vj) -- (2.6,-2.5) node[pos=0.5,right] {$d'_2$};
				\draw[dotted] (1.2,-2.5) -- (1.0,-2.9); \draw[dotted] (2.6,-2.5) -- (2.8,-2.9);
				\fill (vk) circle (2pt); \fill (vj) circle (2pt);
				\node[left=3pt,yshift=6pt] at (vk) {$v_k$};
				\node[right=3pt] at (vj) {$v_{j+1}$};
			\end{tikzpicture}
			\caption{The inductive step: \(v_k=u_{L_i\cup L'_i}\) and
				\(v_{j+1}=u_{L_i}\), with \(d_0=e'_i/r_i\) and \(d_1=e_i/r_i\). The leaves
				seen from \(v_{j+1}\) in the direction of \(w\) are those of
				\(\Lambda_{L'_i}\) and those seen from \(v_k\) in the direction of \(w\).}
			\label{fig:rescale_step}
		\end{figure}
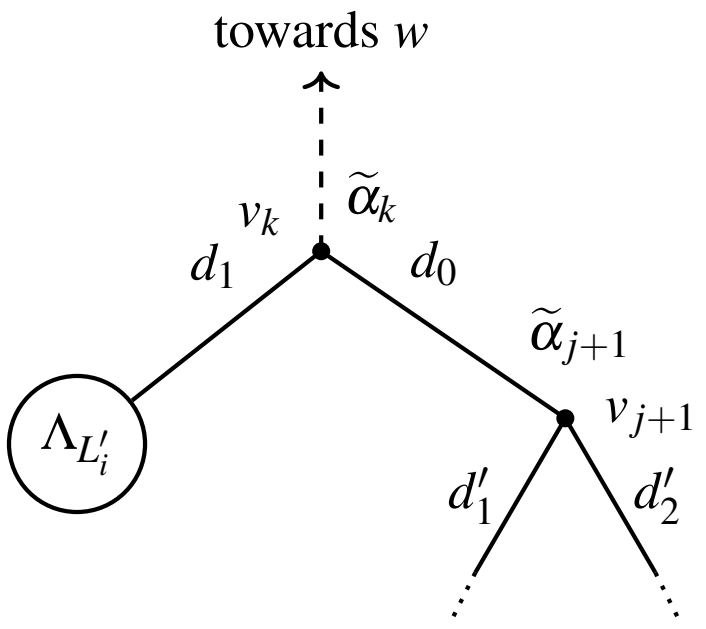

		Assume that the weights pointing towards \(w\) at the
		nodes \(v_2,\dots,v_j\) have already been re-scaled to
		\(\widetilde{\alpha}_2,\dots,\widetilde{\alpha}_j\), in such a way that they satisfy the required conditions. Let us consider the node \(v_{j+1}=u_{L_i}\) and the node
		\(v_k=u_{L_i\cup L'_i}\), \(k\leq j\). The node \(v_k\) has weight
		\(\widetilde{\alpha}_k\) on the edge pointing towards \(w\), weight
		\(d_0=e'_i/r_i\) on the edge going down to \(v_{j+1}\), and weight
		\(d_1=e_i/r_i\) on the edge going down to \(\Lambda_{L'_i}\). The node \(v_{j+1}\) has weight
		\(\alpha_{j+1}=1\) on the edge pointing towards \(w\), which we are going to
		re-scale, and weights \(d'_1,d'_2\) on the edges going down. Notice that
		\(d_0,d_1,d'_1,d'_2\) point away from \(w\), hence they are those of
		Theorem~\ref{tma:converse_rooted} and are not modified. If we replace
		\(\alpha_{j+1}\) by \(\widetilde{\alpha}_{j+1}\), the edge determinant of
		\([v_k,v_{j+1}]\) becomes
		\begin{equation}\label{eqn:aux2escala}
			\widetilde{\alpha}_{j+1}\,d_0-d_1\,\widetilde{\alpha}_k\,d'_1d'_2 .
		\end{equation}
		On the other hand, the leaves seen from \(v_{j+1}\) in the direction of \(w\)
		are the leaves of \(\Lambda_{L'_i}\) and the leaves seen from \(v_k\) in the
		direction of \(w\). The only weight at \(v_k\) not lying on the path from
		\(v_{j+1}\) to the former is \(\widetilde{\alpha}_k\), and to the latter is
		\(d_1\). Hence
		\[
		\overline{S}_{L_i}=\widetilde{\alpha}_k\,S_{L'_i/e'_i}+\frac{e_i}{r_i}\,
		\overline{S}_{L_i\cup L'_i}.
		\]
		This is a gluing of numerical semigroups: indeed,
		\(\widetilde{\alpha}_k\in\overline{S}_{L_i\cup L'_i}\) by the induction
		hypothesis, \(e_i/r_i\in S_{L'_i/e'_i}\) by the semigroup condition at \(v_k\)
		in the direction of \(\Lambda_{L'_i}\), and
		\(\gcd(\widetilde{\alpha}_k,e_i/r_i)=1\). Therefore \(\overline{S}_{L_i}\) is a
		numerical semigroup. At this step, we choose
		\(\widetilde{\alpha}_{j+1}\equiv 1 \pmod{d'_1d'_2}\) such that
		\[
		\widetilde{\alpha}_{j+1}>\max\Bigl\{c(\overline{S}_{L_i}),\,
		\frac{d_1\,\widetilde{\alpha}_k}{d_0}\,d'_1d'_2\Bigr\}.
		\]
		This condition ensures that the weights at \(v_{j+1}\) are pairwise coprime,
		that the edge determinant \eqref{eqn:aux2escala} is positive, and that
		\(\widetilde{\alpha}_{j+1}\in \overline{S}_{L_i}\), thus satisfying the
		semigroup condition at \(v_{j+1}\) in the direction of \(w\). As in the base
		case, \(\widetilde{\alpha}_{j+1}\) lies on the path from \(w\) to every vertex
		and on the path from any node to the leaves lying away from \(w\), so neither
		the linking numbers of \(w\), and hence \(S(\Lambda_w)\), nor the semigroup
		condition at the edges pointing away from \(w\) are affected. Finally,
		\(\widetilde{\alpha}_{j+1}\) only appears in the edge determinants of the edges
		at \(v_{j+1}\) and in the semigroups \(\overline{S}_L\) of the blocks
		\(L\subset L_i\), \(L\neq L_i\), so the edge determinants and the semigroups considered in
		the previous steps are not affected. Hence the induction hypothesis holds for
		\(j+1\).
		
		After finishing the re-scaling of all the nodes, the resulting rooted splice diagram has positive weights,
		pairwise coprime at each node, all its edge determinants are positive, it satisfies the semigroup condition, and its semigroup at \(w\) is \(S(\Lambda_w)=S\).
	\end{proof}
	
	\begin{rem}\label{rem: re-scalearrow}
		On the proof of Theorem \ref{thm:alteration}, the edge weight of the edge on $v_1$ towards the arrow can be any positive integer $\widetilde{\alpha}_1$ satisfying \(\widetilde{\alpha}_1\equiv 1 \pmod{d_1d_2}\). Then, if $\widetilde{\alpha}_1>1$, $\widetilde{\alpha}_2$ is chosen as follows:
		\begin{itemize}
			\item If $|L_1|=1$, we have \(\overline{S}_{L'_1}=\widetilde{\alpha}_1\mathbb{N}+d_1\mathbb{N}.\) We choose $\widetilde{\alpha}_2\equiv 1 \pmod{d'_1d'_2}$ such that
			\[
			\widetilde{\alpha}_2>\max\Bigl\{(\widetilde{\alpha}_1-1)(d_1-1),\,
			\frac{e'_1\widetilde{\alpha}_1}{a_s}\,d'_1d'_2\Bigr\}.
			\]
			
			\item If \(|L_1|>1\) and \(|L'_1|>1\), we have \(\overline{S}_{L_1}=\widetilde{\alpha}_1\,S((\Lambda_{L'_1})_{w_{L'_1}})+d_2\,\mathbb{N}.\)
			
			We choose \(\widetilde{\alpha}_2\equiv 1 \pmod{d'_1d'_2}\) such that 
			\[
			\widetilde{\alpha}_2>\max\Bigl\{c\bigl(\,\overline{S}_{L_1}\bigr),\,
			\frac{e_1\widetilde{\alpha}_1}{e'_1}\,d'_1d'_2\Bigr\}.
			\]
		\end{itemize}
	\end{rem}
	As a consequence, we obtain the following inmediate corollaries.
	\begin{corollary}
		The knot $K$ defined by the splice diagram $\Lambda$ of Theorem \ref{thm:alteration} satisfies 
		\begin{center}
			$\Delta_K(t)=(1-t)\mathcal{G}_S(t)$ 
		\end{center}
		where $\Delta_K(t)$ is the Alexander polynomial of $K$ and $\mathcal{G}_S(t)$ the generating series of $S$.
	\end{corollary}
	\begin{proof}
		The altered weights belong to the geodesics, therefore do not appear in the calculation of the linking number with the only arrow in the splice diagram.
	\end{proof}

\begin{figure}[H]
		\centering
		\begin{subfigure}[b]{0.45\linewidth}
			\centering
			\begin{tikzpicture}[scale=0.8, thick, every node/.style={font=\footnotesize}]
				\coordinate (v1) at (0,0); \coordinate (v2) at (-1.8,-1.1);
				\coordinate (v3) at (1.8,-1.1);
				\draw[->] (v1) -- (0,1.2) node[pos=0.35,right] {$1$};
				\draw (v1) -- (v2) node[pos=0.3,above left=-2pt] {$7$}
				node[pos=0.78,above left=-2pt] {$1$};
				\draw (v1) -- (v3) node[pos=0.3,above right=-2pt] {$5$}
				node[pos=0.78,above right=-2pt] {$1$};
				\draw (v2) -- (-2.5,-2.3) node[pos=0.5,left] {$3$};
				\draw (v2) -- (-1.1,-2.3) node[pos=0.5,right] {$2$};
				\draw (v3) -- (1.1,-2.3) node[pos=0.5,left] {$3$};
				\draw (v3) -- (2.5,-2.3) node[pos=0.5,right] {$2$};
				\foreach \P in {v1,v2,v3} \fill (\P) circle (2pt);
				\foreach \P/\n in {(-2.5,-2.3)/w_0,(-1.1,-2.3)/w_2,(1.1,-2.3)/w_1,(2.5,-2.3)/w_3}
				{\fill \P circle (2pt); \node[below] at \P {$\n$};}
				\node[above] at (0,1.2) {$w$};
				\node[left=3pt,yshift=6pt] at (v1) {$v_1$};
				\node[left=3pt] at (v2) {$v_2$}; \node[right=3pt] at (v3) {$v_3$};
			\end{tikzpicture}
			\caption{Before re-scaling}
		\end{subfigure}
		\hfill
		\begin{subfigure}[b]{0.45\linewidth}
			\centering
			\begin{tikzpicture}[scale=0.8, thick, every node/.style={font=\footnotesize}]
				\coordinate (v1) at (0,0); \coordinate (v2) at (-1.8,-1.1);
				\coordinate (v3) at (1.8,-1.1);
				\draw[->] (v1) -- (0,1.2) node[pos=0.35,right] {$1$};
				\draw (v1) -- (v2) node[pos=0.3,above left=-2pt] {$7$}
				node[pos=0.78,above left=-2pt] {$31$};
				\draw (v1) -- (v3) node[pos=0.3,above right=-2pt] {$5$}
				node[pos=0.78,above right=-2pt] {$43$};
				\draw (v2) -- (-2.5,-2.3) node[pos=0.5,left] {$3$};
				\draw (v2) -- (-1.1,-2.3) node[pos=0.5,right] {$2$};
				\draw (v3) -- (1.1,-2.3) node[pos=0.5,left] {$3$};
				\draw (v3) -- (2.5,-2.3) node[pos=0.5,right] {$2$};
				\foreach \P in {v1,v2,v3} \fill (\P) circle (2pt);
				\foreach \P/\n in {(-2.5,-2.3)/w_0,(-1.1,-2.3)/w_2,(1.1,-2.3)/w_1,(2.5,-2.3)/w_3}
				{\fill \P circle (2pt); \node[below] at \P {$\n$};}
				\node[above] at (0,1.2) {$w$};
				\node[left=3pt,yshift=6pt] at (v1) {$v_1$};
				\node[left=3pt] at (v2) {$v_2$}; \node[right=3pt] at (v3) {$v_3$};
			\end{tikzpicture}
			\caption{After re-scaling}
		\end{subfigure}
		\caption{Re-scaling of the rooted splice diagram of
			Example~\ref{ejemplo:negativedeterminant}.}
		\label{fig:rescaling_example}
	\end{figure}
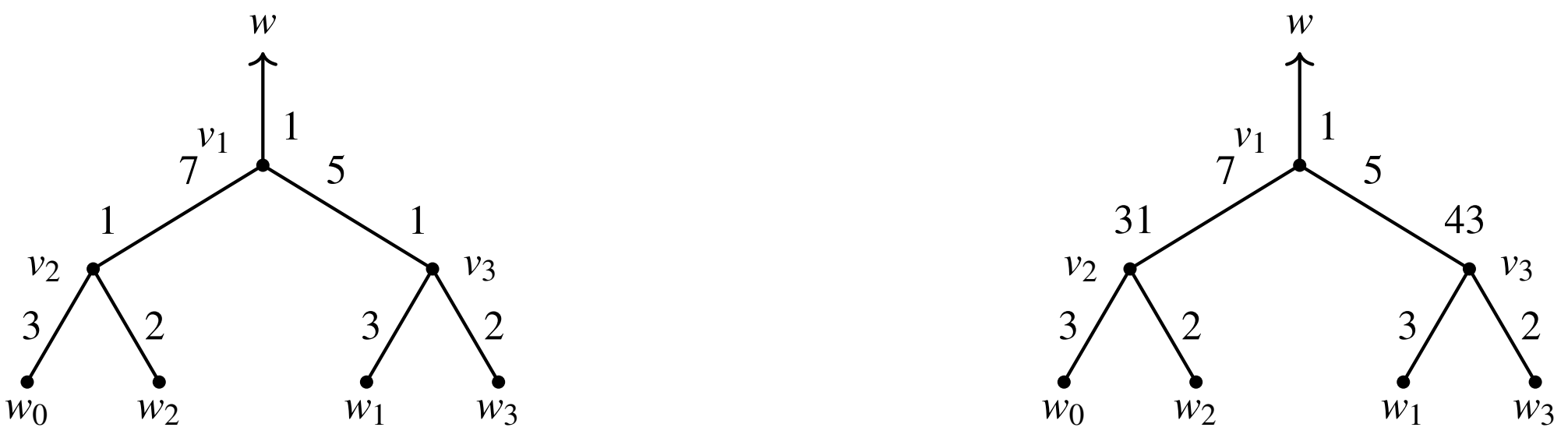

	\begin{ex}\label{ex:rescaling}
		Let us re-scale the rooted splice diagram of
		Example~\ref{ejemplo:negativedeterminant}. With the notation of the proof of
		Theorem~\ref{tma:converse_rooted}, we have \(L_1=\{10,15\}\) and
		\(L'_1=\{14,21\}\), so \(e_1=5\), \(e'_1=7\) and \(r_1=1\). The node \(v_1\)
		has weight \(1\) towards the arrow \(w\), weight \(d_1=e'_1/r_1=7\) towards
		\(\Lambda_{L_1}\) and weight \(d_2=e_1/r_1=5\) towards \(\Lambda_{L'_1}\). The
		node \(v_2\) of \(\Lambda_{L_1}\) (created at the step from \(P_3\) to \(P_2\),
		with \(r_3=5\)) has leaves \(w_0,w_2\) with weights \(15/5=3\) and \(10/5=2\),
		and the node \(v_3\) of \(\Lambda_{L'_1}\) (created at the step from \(P_2\) to
		\(P_1\), with \(r_2=7\)) has leaves \(w_1,w_3\) with weights \(21/7=3\) and
		\(14/7=2\); see Figure~\ref{fig:rescaling_example}~(a). Hence
		\(d'_1d'_2=6\) at both \(v_2\) and \(v_3\), and the two edge determinants are
		\[
		1\cdot 7-6\cdot 5=-23,\qquad 1\cdot 5-6\cdot 7=-37 .
		\]
		
		We are in the case \(|L_1|,|L'_1|>1\),  we take
		\[
		\widetilde{\alpha}_2=\max\{1,5\}\cdot 6\cdot 1+1=31,\qquad
		\widetilde{\alpha}_3=\max\{1,7\}\cdot 6\cdot 1+1=43 .
		\]
		The resulting diagram is shown in Figure~\ref{fig:rescaling_example}~(b).
		Its weights are pairwise coprime at each node, and its edge determinants are
		\[
		31\cdot 7-6\cdot 5=187>0,\qquad 43\cdot 5-6\cdot 7=173>0 .
		\]
		The leaves seen from \(v_2\) in the direction of \(w\) are \(w_1\), \(w_3\)
		and \(w\), with reduced linking numbers \(2\), \(3\) and \(5\), so the
		corresponding semigroup is \(\langle 2,3\rangle\ni 31\); similarly, the one
		seen from \(v_3\) is \(\langle 2,3,7\rangle=\langle 2,3\rangle\ni 43\). Hence
		the semigroup condition holds at every edge. Since only weights pointing
		towards \(w\) have been modified, the linking numbers of \(w\) with the leaves
		are still \(10,14,15,21\), and \(S(\Lambda_w)=\langle 10,14,15,21\rangle\).
		Finally, removing the leaf \(w\), whose weight at \(v_1\) is \(1\), we obtain
		the minimal splice diagram of the ambient integral homology sphere: two nodes
		with weights \((2,3,31)\) and \((2,3,43)\) joined by an edge with edge
		determinant \(31\cdot 43-6\cdot 6=1297>0\), that is, the splice of
		\(\Sigma(2,3,31)\) and \(\Sigma(2,3,43)\) along the fibres of degrees \(31\)
		and \(43\).
	\end{ex}
	
	To finish, we are ready to show the announced realization of complete intersection numerical semigroups as value semigroups of end curves. Let $(X,0)$ be a splice quotient singularity with $\mathbb{Z}HS^3$ link with \(\Lambda\) its splice diagram. Then
	\((X,0)\) satisfies the end curve condition \cite{neumann2005complete,neumann2010end}, i.e for each end vertex $w$ there is a reduced end curve function $f_w\in\mathcal{O}_{X,0}$,
	whose zero locus $C_w=\{f_w=0\}$ is an end curve for $w$; see also
	\cite[Section~8.5]{nemethibook}. Let $\nu_w$ be
	the discrete valuation of $\mathcal{O}_{C_w,0}=\mathcal{O}_{X,0}/(f_w)$ induced by the normalization \(\overline{C_w}\rightarrow C_w.\)
	The value semigroup of $C_w$ is
	$S(C_w):=\nu_w\bigl(\mathcal{O}_{C_w,0}\setminus\{0\}\bigr)$.
	
	\begin{theorem}\label{teo:realizationCI}
		Let \(S\) be a complete intersection numerical semigroup, then there exists a splice quotient singularity \((X,0)\) with \(\mathbb{Z}HS\) link and an end curve \((C,0)\subset (X,0)\) such that its value semigroup is \(S.\)
	\end{theorem}
	\begin{proof}
		Let $G$ be the minimal set of generators of $S$ and let
		$(\widetilde{\Lambda},w)$ be the rooted splice diagram given by
		Theorem~\ref{thm:alteration} and Remark \ref{rem: re-scalearrow}, with weight $\widetilde{\alpha}_1\geq 2$ at
		$v_1$ on the edge towards $w$. Since $G$ is minimal, every leaf weight of
		$\widetilde{\Lambda}$ is at least $2$; together with the positivity of the
		edge determinants this shows that $\widetilde{\Lambda}$ is a minimal splice
		diagram, and $w$ is one of its leaves. By \cite[Theorem~9.4]{EN} and
		\cite[Theorem~2.1]{neumann2005complex}, any system of splice diagram
		equations for $\widetilde{\Lambda}$ defines a singularity $(X,0)$ of splice
		type with $\mathbb{Z}HS^3$ link, and $C=\{z_w=0\}$ is an end curve for $w$.
		
		The local ring $\mathcal{O}_{C,0}=\mathcal{O}_{X,0}/(z_w)$ is defined by the
		splice diagram equations with $z_w=0$. With respect to the weights
		$\nu_w(z_e)=\ell_{ew}$, $e\neq w$, their leading forms define the monomial
		curve of $S(\widetilde{\Lambda}_w)$, and the remaining terms have strictly
		higher weight \cite[Section~3]{neumann2005complex}. Hence $C$ is a positive
		weight deformation of this monomial curve and
		$S(C)=S(\widetilde{\Lambda}_w)=S$.
	\end{proof}
	
	
	\section{Poincaré series of curves on splice quotient singularities}\label{section:four}
	
	In Section \ref{section:three}, we have shown that any complete intersection numerical semigroup is the semigroup of an end curve. Moreover, its generating series coincides with the Alexander polynomial of the knot associated to the end curve. Now, we will shift to the analytic setting. We will show that this identity can be carried out to an identity between the Alexander polynomial and the Poincaré series of the curve. More generally, we will show that for a reduced curve whose defining equation is a product of end curve functions, its Poincaré series coincides with the Alexander polynomial of the corresponding link. Furthermore, we provide an explicit description of the semigroup of integer multiplicities of functions along a divisor of the resolution.
	\medskip
	
	Let us recall some of the concepts from Sections \ref{ssec:Splice_diagrams} and  \ref{ssec:splice-quotient}. Let $(X,0)$ be a splice type surface singularity with $\mathbb{Z}HS^3$ link $\Sigma.$ Let $\pi:(\Tilde{X},0)\rightarrow (X,0)$ be a good resolution admitting end curve functions (one for each leaf of the dual graph). Denote $E:=\pi^{-1}(0)$ the exceptional curve, $\{E_v\}$ its irreducible components, $\Gamma$ the dual graph and $A(\Gamma)$ the matrix associated to the intersection form of the plumbing. Set $L:=H_2(\Tilde{X},\mathbb{Z})$, it is freely generated over $\mathbb{Z}$ by the classes $\{E_v\}$ and is endowed with the intersection form defined by $A(\Gamma)$. Moreover, $A(\Gamma)$ defines an embedding of abelian groups  $L\hookrightarrow L':=H^2(\Tilde{X},\mathbb{Z})$ with $L'/L\cong H_1(\Sigma,\mathbb{Z})$ and the intersection form $(\cdot,\cdot)$ on $L$ may be extended to $L'.$ Since the link of $(X,0)$ is an integral homology sphere then, $H_1(\Sigma,\mathbb{Z})$ is trivial and $L\cong L'$. Hence, the negative duals $\{E_v^*\}$ , that is $(E_i^*,E_j)=-\delta_{ij}$ (the Kronecker delta function), constitute a free basis for $L$ and $A(\Gamma)$ is the matrix of the change of basis $\{E_v\}$ to $\{-E_v^*\}.$
	
	Let $\Lambda$ and $\Lambda^{max}_\Gamma$ be the splice diagram and maximal splice diagram (obtained from $\Gamma$) for $\Sigma$. Recall that $\Lambda^{max}_\Gamma$ has the same set of vertices as $\Gamma$. Therefore, we denote $\mathcal{V}$ to the set of vertices of $\Gamma$ and of $\Lambda^{max}_\Gamma$; if necessary we will specify the object we are working with. Then, 
	\begin{equation}\label{eq:dual_base_expresion}
		E_j^*=\sum_{i\in\mathcal{V}}\hspace{1mm}\ell_{ij}E_i \text{ for each } j\in\mathcal{V}.
	\end{equation}
	
	We label the vertices of $\Gamma$ and suppose $\mathcal{V}=\{1,\dots,s\}$. Let $\mathbb{Z}[[\mathbf{t}^{\pm 1}]]:=\mathbb{Z}[[t_1^{\pm1},\dots, t_s^{\pm1}]]$ be the ring of formal power series in variables $\{t_1,\dots,t_s\}.$ One may associate a subring  of $\mathbb{Z}[[\mathbf{t}^{\pm 1}]]$ to the resolution $\pi$ as follows: Each $l\in L$ admits a unique expression of the form $l=\sum_j \hspace{1mm}l_j E_j$. Hence, the element $\mathbf{t}^l:=t_1^{l_1}\cdots t_s^{l_s}$ is well defined and we may consider the subring $\mathbb{Z}[[\mathbf{t}^L]]$ generated by the monomials $\mathbf{t}^l,\hspace{1.5mm} l\in L.$ With this notation, we introduce the topological function:
	\begin{equation}\label{eq:zeta_griega_function}
		\zeta(\mathbf{t})=\prod_{j\in\mathcal{V}}(1-\mathbf{t}^{E_j^*})^{\delta_j -2},
	\end{equation}
	where $\delta_j$ is the valency of the vertex $j$ in $\Gamma.$ Occasionally, we would prefer to focus in a particular subset of vertices. For instance, fix a vertex \(v\in \mathcal{V}\) and its corresponding exceptional divisor \(E_v,\) define the projection $pr_v:L\rightarrow L$ by $pr_v(\sum_j\hspace{1mm} l_jE_j)=l_vE_v$ and set the monomial $t^{pr_v(l)}=t^{l_v}=\mathbf{t}^{l}|_{t_j=1 \text{ for all } j\neq v}$ 
	
	\begin{defin}\label{def:zeta_vertex}
		The topological function $\zeta_v(t)$ associated to a vertex $v\in\mathcal{V}$ is the function
		\begin{equation*}
			\zeta_v(t):=\zeta(\mathbf{t})|_{t_j=1 \text{ for all } j\neq i}.
		\end{equation*}
	\end{defin}
	Moreover, Eq. \ref{eq:dual_base_expresion} implies that the function $\zeta_v(t)$ can be computed in terms of linking numbers on the maximal splice diagram,
	\begin{equation}\label{eq:zeta_vertex_linking_numbers}
		\zeta_v(t)=\prod_{j\in\mathcal{V}}(1-t^{\ell_{vj}})^{\delta_j -2}.
	\end{equation}
	
	The function $\zeta_v(t)$ has a topological interpretation in terms of Alexander polynomials as well. Recall that any vertex $v\in\mathcal{V}$ defines a knot $K_v$ embedded in the singularity link $\Sigma$; see \cite[Appendix, part 2]{neumann1981calculus} and \cite[Chapter V]{EN}; to point out this knot we add to the plumbing graph $\Gamma$ one unweighted vertex $v'$, which will be drawn as an arrow, and one edge from $v$ to $v'$. The maximal splice diagram $\Lambda^{max}_\Gamma$ is modified analogously as $\Gamma$ and the edge from $v$ to $v'$ has weight $1$ at $v$. We denote $\Gamma_v$ and $\Lambda^{max}_{\Gamma_v}$, the plumbing graph and maximal splice diagram with an arrow added to the vertex $v$. Then, $\Lambda^{max}_{\Gamma_v}$ defines a splice diagram with a unique arrow $\Tilde{\Lambda}$ which can be obtained from $\Lambda$ as follows. (We will distinguish the following cases depending on the valence of the vertex $v$ in $\Lambda^{max}_\Gamma$). 
	
	\begin{itemize}
		\item Suppose $\delta_v=1$. If $\Gamma$ has no vertex of valence greater than two then, $\Lambda$ was either the empty splice diagram or consists of two vertices joined with an edge. Then, $\Tilde{\Lambda}$ consists of one vertex $w$, an arrow $v'$ and one edge with edge weight $1$ joining $w$ and $v'$. Otherwise, let $w$ be the nearest vertex to $v$ of valence greater than two, let $e$ be the edge adjacent to $w$ pointing in the direction of $v$ and let $d_1$ be the edge weight of $e$ on $w$. Since  $\delta_w\geq 2$, $w$ corresponds to a node in $\Lambda$ and one of its edges joins $w$ with a leaf $v'$ and has edge weight $d_1$. Then, $\Tilde{\Lambda}$ is obtained from $\Lambda$ replacing the vertex $v'$ with an arrow (also denoted by $v'$). 
		\begin{figure}[H]
			\centering
			\begin{subfigure}{0.28\textwidth}
				\centering
				\begin{tikzpicture}[scale=1, thick]
					\coordinate (V3) at (7.5, 0);
					
					\coordinate (TL2_solid)  at (6.4, 0.75);
					\coordinate (TL2_dotted) at (6.0, 1.05);
					
					\coordinate (BL2_solid)  at (6.4, -0.75);
					\coordinate (BL2_dotted) at (6.0, -1.05);
					
					\coordinate (R2_end)    at (8.5, 0);
					
					\draw (V3) -- (TL2_solid) node[pos=0.4, above =1pt] {$d_2$};
					\draw[dotted] (TL2_solid) -- (TL2_dotted);
					
					\draw (V3) -- (BL2_solid) node[pos=0.4, below =1pt] {$d_n$};
					\draw[dotted] (BL2_solid) -- (BL2_dotted);
					
					\fill (6.2,  0.3) circle (1.2pt);
					\fill (6.1,  0.0) circle (1.2pt);
					\fill (6.2, -0.3) circle (1.2pt);
					
					\draw (V3) -- (R2_end) node[midway, above=1pt] {$d_1$};
					
					\fill (V3) circle (2pt) node[below=1pt] {$w$};
					\fill (R2_end) circle (2pt);
				\end{tikzpicture}
				\caption{ Piece of the splice diagram $\Lambda$.}
				\label{fig:lambda_1_1}
			\end{subfigure}
			\hfill
			\begin{subfigure}{0.28\textwidth}
				\centering
				\begin{tikzpicture}[scale=1, thick]
					\coordinate (V1) at (0, 0);
					
					\coordinate (TL1_solid)  at (-1.1, 0.75);
					\coordinate (TL1_dotted) at (-1.5, 1.05);
					
					\coordinate (BL1_solid)  at (-1.1, -0.75);
					\coordinate (BL1_dotted) at (-1.5, -1.05);
					
					\coordinate (R1_dot_L)  at (0.9, 0);
					\coordinate (R1_dot_R)  at (1.5, 0);
					\coordinate (V2)        at (2.2, 0);
					\coordinate (R1_end)    at (3, 0);
					
					\draw (V1) -- (TL1_solid) node[pos=0.4, above =1pt] {$d_2$};
					\draw[dotted] (TL1_solid) -- (TL1_dotted);
					
					\draw (V1) -- (BL1_solid) node[pos=0.4, below =1pt] {$d_n$};
					\draw[dotted] (BL1_solid) -- (BL1_dotted);
					
					\fill (-1.4,  0.3) circle (1.2pt);
					\fill (-1.5,  0.0) circle (1.2pt);
					\fill (-1.4, -0.3) circle (1.2pt);
					
					\draw (V1) -- (R1_dot_L) node[midway, above=1pt] {$d_1$};
					\draw[dotted] (R1_dot_L) -- (R1_dot_R);
					\draw (R1_dot_R) -- (V2);
					\draw[->] (V2) -- (R1_end);
					
					\fill (V1) circle (2pt) node[below=1pt] {$w$};
					\fill (V2) circle (2pt) node[above=1pt] {$v$};
				\end{tikzpicture}
				\caption{ Piece of the maximal splice diagram $\Lambda^{\text{max}}_{\Gamma_v}$.}
				\label{fig:lambda_2_1}
			\end{subfigure}
			\hfill
			\begin{subfigure}{0.28\textwidth}
				\centering
				\begin{tikzpicture}[scale=1, thick]
					\coordinate (V3) at (7.5, 0);
					
					\coordinate (TL2_solid)  at (6.4, 0.75);
					\coordinate (TL2_dotted) at (6.0, 1.05);
					
					\coordinate (BL2_solid)  at (6.4, -0.75);
					\coordinate (BL2_dotted) at (6.0, -1.05);
					
					\coordinate (R2_end)    at (8.5, 0);
					
					\draw (V3) -- (TL2_solid) node[pos=0.4, above =1pt] {$d_2$};
					\draw[dotted] (TL2_solid) -- (TL2_dotted);
					
					\draw (V3) -- (BL2_solid) node[pos=0.4, below =1pt] {$d_n$};
					\draw[dotted] (BL2_solid) -- (BL2_dotted);
					
					\fill (6.2,  0.3) circle (1.2pt);
					\fill (6.1,  0.0) circle (1.2pt);
					\fill (6.2, -0.3) circle (1.2pt);
					
					\draw[->] (V3) -- (R2_end) node[midway, above=1pt] {$d_1$};
					
					\fill (V3) circle (2pt) node[below=1pt] {$w$};
					\fill (R2_end) node[right=1pt] {$v'$};
				\end{tikzpicture}
				\caption{ Piece of the splice diagram $\Tilde{\Lambda}$ .}
				\label{fig:lambda_3_1}
			\end{subfigure}
			
			\caption{Construction of the splice diagram $\Tilde{\Lambda}$ when $\delta_v=1$}
			\label{fig:arrow_leaf_maximal_splice_1}
		\end{figure}
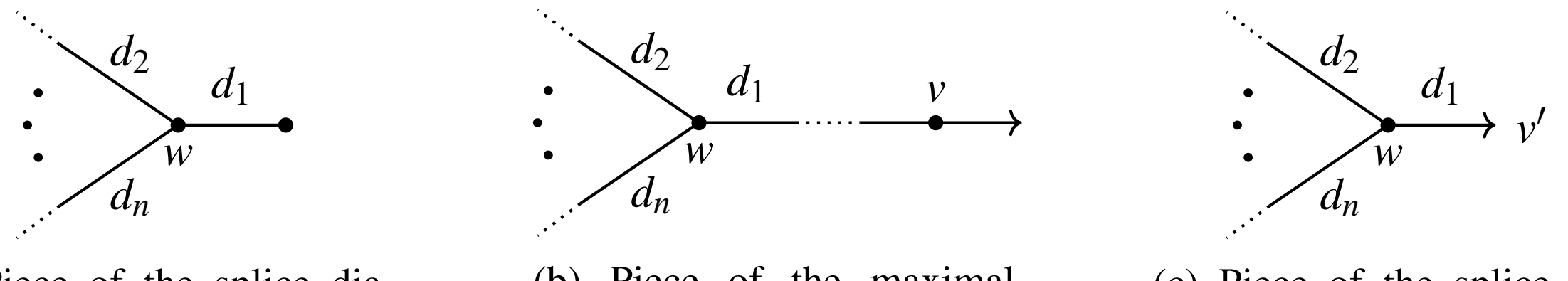
		
		\item Suppose $\delta_v\geq 3$. In this case, $v$ corresponds to a node $w$ in $\Lambda$ (so $\Lambda$ has at least one node). One forms $\Tilde{\Lambda}$ by adding in $\Lambda$ an arrow vertex $v'$ and an edge $e$  connecting $w$ to $v'$; the edge weight of $e$ at $w$ is $1$. \\
		
		\begin{figure}[H]
			\centering
			\begin{subfigure}{0.28\textwidth}
				\centering
				\begin{tikzpicture}[scale=1, thick]
					\coordinate (V3) at (7.5, 0);
					
					\coordinate (TL2_solid)  at (6.4, 0.75);
					\coordinate (TL2_dotted) at (6.0, 1.05);
					
					\coordinate (BL2_solid)  at (6.4, -0.75);
					\coordinate (BL2_dotted) at (6.0, -1.05);
					
					\coordinate (R2_end)    at (8.5, 0);
					
					\draw (V3) -- (TL2_solid) node[pos=0.4, above =1pt] {$d_1$};
					\draw[dotted] (TL2_solid) -- (TL2_dotted);
					
					\draw (V3) -- (BL2_solid) node[pos=0.4, below =1pt] {$d_n$};
					\draw[dotted] (BL2_solid) -- (BL2_dotted);
					
					\fill (6.2,  0.3) circle (1.2pt);
					\fill (6.1,  0.0) circle (1.2pt);
					\fill (6.2, -0.3) circle (1.2pt);

					\fill (V3) circle (2pt) node[below=1pt] {$w$};
				\end{tikzpicture}
				\caption{ Piece of the splice diagram $\Lambda$.}
				\label{fig:lambda_1_2}
			\end{subfigure}
			\hfill
			\begin{subfigure}{0.28\textwidth}
				\centering
				\begin{tikzpicture}[scale=1, thick]
					\coordinate (V1) at (0, 0);
					
					\coordinate (TL1_solid)  at (-1.1, 0.75);
					\coordinate (TL1_dotted) at (-1.5, 1.05);
					
					\coordinate (BL1_solid)  at (-1.1, -0.75);
					\coordinate (BL1_dotted) at (-1.5, -1.05);
					
					\coordinate (R1_end)    at (1.0, 0);
					
					\draw (V1) -- (TL1_solid) node[pos=0.3, above =1pt] {$d_1$};
					\draw[dotted] (TL1_solid) -- (TL1_dotted);
					
					\draw (V1) -- (BL1_solid) node[pos=0.3, below =1pt] {$d_n$};
					\draw[dotted] (BL1_solid) -- (BL1_dotted);
					
					\fill (-1.4,  0.3) circle (1.2pt);
					\fill (-1.5,  0.0) circle (1.2pt);
					\fill (-1.4, -0.3) circle (1.2pt);
					
					\draw[->] (V1) -- (R1_end);
					
					\fill (V1) circle (2pt);
					
					\fill (V1) circle (2pt) node[below=1pt] {$v$};
					
				\end{tikzpicture}
				\caption{ Piece of the maximal splice diagram $\Lambda^{\text{max}}_{\Gamma_v}$.}
				\label{fig:lambda_2_3}
			\end{subfigure}
			\hfill
			\begin{subfigure}{0.28\textwidth}
				\centering
				\begin{tikzpicture}[scale=1, thick]
					\coordinate (V3) at (7.5, 0);
					
					\coordinate (TL2_solid)  at (6.4, 0.75);
					\coordinate (TL2_dotted) at (6.0, 1.05);
					
					\coordinate (BL2_solid)  at (6.4, -0.75);
					\coordinate (BL2_dotted) at (6.0, -1.05);
					
					\coordinate (R2_end)    at (8.5, 0);
					
					\draw (V3) -- (TL2_solid) node[pos=0.4, above =1pt] {$d_1$};
					\draw[dotted] (TL2_solid) -- (TL2_dotted);
					
					\draw (V3) -- (BL2_solid) node[pos=0.4, below =1pt] {$d_n$};
					\draw[dotted] (BL2_solid) -- (BL2_dotted);
					
					\fill (6.2,  0.3) circle (1.2pt);
					\fill (6.1,  0.0) circle (1.2pt);
					\fill (6.2, -0.3) circle (1.2pt);
					
					\draw[->] (V3) -- (R2_end) node[midway, above=1pt] {$d_1$};
					
					\fill (V3) circle (2pt) node[below=1pt] {$w$};
					\fill (R2_end) node[right=1pt] {$v'$};
				\end{tikzpicture}
				\caption{ Piece of the splice diagram $\Tilde{\Lambda}$ .}
				\label{fig:lambda_3_3}
			\end{subfigure}
			
			\caption{Construction of the splice diagram $\Tilde{\Lambda}$ when $\delta_v\geq 3$}
			\label{fig:arrow_leaf_maximal_splice_3}
		\end{figure}
		
		\item Suppose $\delta_v=2$. Then, $v$ has two adjacent edges in $\Lambda^{max}_\Gamma$ with (coprime) edge weights, say $a,b$ respectively. If $\Gamma$ has no vertices of valence greater than two then, $\Lambda$ is the empty splice diagram or consists of two vertices of valence one joined with an edge. In this case, $\Tilde{\Lambda}$ consists of one node of valence three, two leaves (each of them joined with the node by an edge) and one arrow appended to the node. The edge weights are $a,b$ and $1$, and the $1$ corresponds to the edge pointing towards the arrow $v'$. Otherwise, let $w_1$ (resp. $w_2$) be the nearest vertex to $v$ with valence not two in $\Lambda^{max}_\Gamma$ in the direction of the edge with edge weight $a$ (resp. $b$). Since $\delta_{w_1},\delta_{w_2}\neq 2$, $w_1$ and $w_2$ correspond to vertices in $\Lambda$ and there is an edge $e$ connecting them. Hence, $\Tilde{\Lambda}$ is obtained from $\Lambda$ by halving $e$ with a vertex $w$ of valence three. The edge joining $w$ with $w_1$ has edge weight $a$ at the end adjacent to $w$ and has edge weight the edge weight of $e$ at $w_1$ on the other end. The edge joining $w$ with $w_2$ has edge weight $b$ at the end adjacent to $w$ and has edge weight the edge weight of $e$ at $w_2$ on the other end. The remaining edge joins $w$ to an arrow $v'$ and has edge weight $1$.\\
		
		\begin{figure}[H]
			\centering
			\begin{subfigure}{0.28\textwidth}
				\centering
				\begin{tikzpicture}[scale=1, thick]
					\coordinate (V3) at (-3.5, 0);
					
					\coordinate (TL2_solid)  at (-4.6, 0.75);
					\coordinate (TL2_dotted) at (-5, 1.05);
					
					\coordinate (BL2_solid)  at (-4.6, -0.75);
					\coordinate (BL2_dotted) at (-5.0, -1.05);
					
					\coordinate (R2)    at (-1.5, 0);
					
					\draw (V3) -- (TL2_solid) node[pos=0.4, above =1pt] {$d_2$};
					\draw[dotted] (TL2_solid) -- (TL2_dotted);
					
					\draw (V3) -- (BL2_solid) node[pos=0.4, below =1pt] {$d_n$};
					\draw[dotted] (BL2_solid) -- (BL2_dotted);
					
					\fill (-4.8,  0.3) circle (1.2pt);
					\fill (-4.9,  0.0) circle (1.2pt);
					\fill (-4.8, -0.3) circle (1.2pt);
					
					\draw (V3) -- (R2) node[pos=0.2, above=0pt] {$d_1$} node[pos=0.8, above=0pt] {$d_1'$};
					\coordinate (TR2_solid)  at (-0.4, 0.75);
					\coordinate (TR2_dotted) at (0, 1.05);
					
					\coordinate (BR2_solid)  at (-0.4, -0.75);
					\coordinate (BR2_dotted) at (0, -1.05);
					\draw (R2) -- (TR2_solid) node[pos=0.4, above =1pt] {$d_2'$};
					\draw[dotted] (TR2_solid) -- (TR2_dotted);
					
					\draw (R2) -- (BR2_solid) node[pos=0.4, below =1pt] {$d_s'$};
					\draw[dotted] (BR2_solid) -- (BR2_dotted);
					\fill (-0.2,  0.3) circle (1.2pt);
					\fill (-0.1,  0.0) circle (1.2pt);
					\fill (-0.2, -0.3) circle (1.2pt);
					\fill (V3) circle (2pt) node[below=1.1pt] {\small $w_1$};
					\fill (R2) circle (2pt) node[below=1.1pt] {\small $w_2$};
				\end{tikzpicture}
				\caption{ Piece of the splice diagram $\Lambda$.}
				\label{fig:lambda_1_3}
			\end{subfigure}
			\hfill
			\begin{subfigure}{0.28\textwidth}
				\centering
				\begin{tikzpicture}[scale=1, thick]

					\coordinate (V1_top) at (0, 1.0);
					\coordinate (V1_bot) at (0, -1.0);
					\coordinate (R1)     at (2.0, 0);
					\coordinate (R1_end) at (2.8, 0);
					
					\coordinate (TL1_top_solid)  at (-0.8, 1.5);
					\coordinate (TL1_top_dotted) at (-1.2, 1.75);
					\coordinate (TL1_bot_solid)  at (-0.8, 0.5);
					\coordinate (TL1_bot_dotted) at (-1.2, 0.25);
					
					\draw (V1_top) -- (TL1_top_solid) node[pos=0.4, above=1pt] {$d_2$};
					\draw[dotted] (TL1_top_solid) -- (TL1_top_dotted);
					
					\draw (V1_top) -- (TL1_bot_solid) node[pos=0.6, below =1pt] {$d_n$};
					\draw[dotted] (TL1_bot_solid) -- (TL1_bot_dotted);
					
					\fill (-1.0, 1.2) circle (1.2pt);
					\fill (-1.1, 1.0) circle (1.2pt);
					\fill (-1.0, 0.8) circle (1.2pt);
					
					\coordinate (BL1_top_solid)  at (-0.8, -0.5);
					\coordinate (BL1_top_dotted) at (-1.2, -0.25);
					\coordinate (BL1_bot_solid)  at (-0.8, -1.5);
					\coordinate (BL1_bot_dotted) at (-1.2, -1.75);
					
					\draw (V1_bot) -- (BL1_top_solid) node[pos=0.6, above =1pt] {$d'_2$};
					\draw[dotted] (BL1_top_solid) -- (BL1_top_dotted);
					
					\draw (V1_bot) -- (BL1_bot_solid) node[pos=0.4, below =1pt] {$d'_s$};
					\draw[dotted] (BL1_bot_solid) -- (BL1_bot_dotted);
					
					\fill (-1.0, -0.8) circle (1.2pt);
					\fill (-1.1, -1.0) circle (1.2pt);
					\fill (-1.0, -1.2) circle (1.2pt);
					
					\draw (V1_top) -- (0.6, 0.7) node[pos=1.2, above =1pt] {$d_1$};
					\draw[dotted] (0.6, 0.7) -- (1.4, 0.3);
					\draw (1.4, 0.3) -- (R1) node[pos=0.4, above =1.5pt] {$a$};
					
					\draw (V1_bot) -- (0.6, -0.7) node[pos=1.2, below =1pt] {$d'_1$};
					\draw[dotted] (0.6, -0.7) -- (1.4, -0.3);
					\draw (1.4, -0.3) -- (R1) node[pos=0.4, below =1.5pt] {$b$};
					
					\draw[->] (R1) -- (R1_end) node[pos=0.1, above=1pt] {$v$};
					
					\fill (V1_top) circle (2pt) node[below=0.5pt] {\small $w_1$};
					\fill (V1_bot) circle (2pt) node[above=0pt] {\small $w_2$};
					\fill (R1)     circle (2pt);

				\end{tikzpicture}
				\caption{ Piece of the maximal splice diagram $\Lambda^{\text{max}}_{\Gamma_v}$.}
				\label{fig:lambda_2_2}
			\end{subfigure}
			\hfill
			\begin{subfigure}{0.28\textwidth}
				\centering
				\begin{tikzpicture}[scale=1, thick]
					\coordinate (V2_top) at (7.0, 1.0);
					\coordinate (V2_bot) at (7.0, -1.0);
					\coordinate (R2)     at (9.0, 0);
					\coordinate (R2_end) at (10, 0);
					
					\coordinate (TL2_top_solid)  at (6.2, 1.5);
					\coordinate (TL2_top_dotted) at (5.8, 1.75);
					\coordinate (TL2_bot_solid)  at (6.2, 0.5);
					\coordinate (TL2_bot_dotted) at (5.8, 0.25);
					
					\draw (V2_top) -- (TL2_top_solid) node[pos=0.4, above =1pt] {$d_2$};
					\draw[dotted] (TL2_top_solid) -- (TL2_top_dotted);
					
					\draw (V2_top) -- (TL2_bot_solid) node[pos=0.6, below =1pt] {$d_n$};
					\draw[dotted] (TL2_bot_solid) -- (TL2_bot_dotted);
					
					\fill (6.0, 1.2) circle (1.2pt);
					\fill (5.9, 1.0) circle (1.2pt);
					\fill (6.0, 0.8) circle (1.2pt);
					
					\coordinate (BL2_top_solid)  at (6.2, -0.5);
					\coordinate (BL2_top_dotted) at (5.8, -0.25);
					\coordinate (BL2_bot_solid)  at (6.2, -1.5);
					\coordinate (BL2_bot_dotted) at (5.8, -1.75);
					
					\draw (V2_bot) -- (BL2_top_solid) node[pos=0.6, above =1pt] {$d'_2$};
					\draw[dotted] (BL2_top_solid) -- (BL2_top_dotted);
					
					\draw (V2_bot) -- (BL2_bot_solid) node[pos=0.4, below =1pt] {$d'_s$};
					\draw[dotted] (BL2_bot_solid) -- (BL2_bot_dotted);
					
					\fill (6.0, -0.8) circle (1.2pt);
					\fill (5.9, -1.0) circle (1.2pt);
					\fill (6.0, -1.2) circle (1.2pt);
					
					\draw (V2_top) -- (R2) node[pos=0.2, above =1pt] {$d_1$} node[pos=0.8, above =1pt] {$a$};
					\draw (V2_bot) -- (R2) node[pos=0.3, below =1pt] {$d'_1$} node[pos=0.8, below=1pt] {$b$};
					
					\draw[->] (R2) -- (R2_end) node[pos=0.3, above=1pt] {$1$};
					
					\fill (V2_top) circle (2pt) node[below=1pt] {\small $w_1$};
					\fill (V2_bot) circle (2pt) node[above=1pt] {\small $w_2$};
					\fill (R2)     circle (2pt) node[below=1pt] {\small $w$};
					\fill (R2_end) node[right=1pt] {$v'$};
					
				\end{tikzpicture}
				\caption{ Piece of the splice diagram $\Tilde{\Lambda}$ .}
				\label{fig:lambda_3_2}
			\end{subfigure}
			
			\caption{Construction of the splice diagram $\Tilde{\Lambda}$ when $\delta_v=2$}
			\label{fig:arrow_leaf_maximal_splice_2}
		\end{figure}
	\end{itemize}
	
	Observe that $\Tilde{\Lambda}$ is also a splice diagram for $\Sigma$ and $\Tilde{\Lambda}$ has a unique arrow which stands for the knot $K_v$. Hence, \cite{EN} shows that the Alexander polynomial of $K$ can be computed from $\Tilde{\Lambda}$ as in Theorem \ref{teo:alexpol_splice}.
	
	\begin{theorem}\label{teo:zeta_alex_poly}
		Let $(X,0)$ be a splice type surface singularity with $\mathbb{Z}HS^3$ link $\Sigma.$ Let $\pi:(\Tilde{X},0)\rightarrow (X,0)$ be a good resolution admitting end curve functions. Let $\Gamma$ be the dual graph of the resolution $\pi$, let $v$ be a vertex of $\Gamma$, let $\zeta_v(t)$ be the topological function from Definition \ref{def:zeta_vertex}, let $K_v$ be the knot embedded in $\Sigma$ defined by $v$ and let $\Delta_{K_v}(t)$ be its Alexander polynomial. Then,
		\begin{equation*}
			(1-t)(1-t^{\ell_{vv}})\zeta_v(t)=\Delta_{K_v}(t)
		\end{equation*}
	\end{theorem}
	\begin{proof}
		Let $\Lambda^{max}_\Gamma$ be the maximal splice diagram obtained from $\Gamma.$ Modify $\Lambda^{max}_\Gamma$ to indicate the knot $K_v$, that is, add a vertex $v'$ (as an arrow) and an edge from $v$ to $v'$ with edge weight $1$ on the end adjacent to $v.$ Use the modification of $\Lambda^{max}_\Gamma$ to construct the splice diagram $\Tilde{\Lambda}$. Let $\mathcal{E}$ (resp. $\mathcal{N}$) be the number of leaves (resp. nodes) of $\Tilde{\Lambda}$; we count the arrow as a leaf as well. Then, Theorem \ref{teo:alexpol_splice} and Equation \ref{eq:zeta_vertex_linking_numbers} yield
		\begin{equation*}
			(t-1)(1-t^{\ell_{vv}})\zeta_v(t)(-1)^s=\Delta_{K_v}(t),
		\end{equation*}
		where $s=\sum_{j\in \mathcal{N}} (\delta_j-2)-(|\mathcal{E}|-1)$. A simple computation shows that $\sum_{j\in \mathcal{N}} (\delta_j-2)=|\mathcal{E}|-2,$ which concludes the proof.
	\end{proof}
	
	\subsection{Poincaré series and the semigroup associated to a divisor}\label{subsection:41}
	
	Let $(X,0)$ be a splice quotient-singularity with $\mathbb{Z}HS^3$ link and let $\pi: (\Tilde{X},0)\rightarrow (X,0)$ be the good resolution associated with the splice-quotient analytic structure. Then, one considers the filtration 
	\begin{equation*}
		\mathcal{F}(l):=\{f\in \mathcal{O}_{X,0}\hspace{1.5mm} | \hspace{1.5mm} div(f\circ\pi)\geq l\},\hspace{1.5mm} l\in L,
	\end{equation*}
	
	where $L$ is the set of integral divisors and $\mathcal{O}_{X,0}$  is the local ring of holomorphic functions. Denote $E:=\sum_{j\in\mathcal{V}} E_j \in L$. The quotient $ \mathcal{F}(l)/ \mathcal{F}(l+E)$ is a finite dimensional $\mathbb{C}-$vector space. Therefore, the multi-index filtration  $\mathcal{F}(l)$ defines a Poicaré series (see \cite{CDG04,ZPNemethi}) as follows. Set the series \begin{equation*}
		L(\mathbf{t}):=\sum_{l\in L}\dim_{\mathbb{C}}(\mathcal{F}(l)/\mathcal{F}(l+E))\mathbf{t}^l\in \mathbb{Z}[[\mathbf{t}^L]]
	\end{equation*}  
	
	Moreover, $ \mathcal{F}(l)= \mathcal{F}(l+E_j)$ if $(l,E_j^*)>0$ (i.e if $l_j<0$ in $l=\sum l_jE_j$) thus, a straightforward computation shows that $L(\mathbf{t})\prod_j (t_j-1)$ belongs to $\mathbb{Z}[[\mathbf{t}^{L_{\geq 0}}]]$. The Poincaré series associated with the filtration $\mathcal{F}(l)$ is
	\begin{equation*}
		P(\mathbf{t}):=-\frac{L(\mathbf{t})\prod_j (t_j-1)}{(1-\mathbf{t}^E)}\in \mathbb{Z}[[\mathbf{t}^{L_{\geq 0}}]].
	\end{equation*}
	
	As with  $\zeta(\mathbf{t}),$ instead of working with the full set $\{E_j\}_{j\in\mathcal{V}}$ of divisors, we sometimes rather work with a concrete exceptional divisor $E_v$. Hence, the filtration needs to be adapted to this context. Fix $v$ a vertex and let $pr_v$ be the projection. For any $l_vE_v\in pr_v(L)$, we define the filtration
	\begin{equation}\label{eq:filtration_vertex}
		\mathcal{F}_v(l_vE_v):=\{f\in \mathcal{O}_{X,0}\hspace{1.5mm} | \hspace{1.5mm} div(f\circ\pi)\geq l_vE_v\}.
	\end{equation}
	
	Generally, we will write $\mathcal{F}_v(l_v)$ instead of $\mathcal{F}_v(l_vE_v).$ In this framework, $L_v(t)$ and $P_v(t)$ are series in one variable and both series coincide, that is,
	\begin{equation*}
		P_v(t)=L_v(t)=\sum_{l\in \mathbb{Z}_{\geq0}}\dim_{\mathbb{C}}(\mathcal{F}_v(l)/\mathcal{F}_v(l+1))t^l.
	\end{equation*}
	
	From the definition can be seen that $P(\mathbf{t})$ is an analytic invariant, nevertheless is common in normal surface singularities to analize if the Poincaré series are just analytic or in contrary, they only depend on the topological aspects of the surface singularity. Although the latter is generally false, for splice-quotient singularities the Poincaré series is topological. 
	\begin{theorem}\cite[Theorem 5.1.1]{nemethi2012cohomology}\label{teo:Z=P}
		Let $(X,0)$ be a splice-quotient singularity with $\mathbb{Z}HS^3-$link, $\pi:(\Tilde{X},0)\rightarrow (X,0)$ the good resolution that agrees with the splice-quotient analytic structure and $P(\mathbf{t})$ the Poincaré series associated with this resolution. Then,
		\begin{equation*}
			P(\mathbf{t})=\zeta(\mathbf{t}),
		\end{equation*}
		where $\zeta(\mathbf{t})$ is defined as in Equation \ref{eq:zeta_griega_function}.
	\end{theorem}
	Theorem \ref{teo:Z=P} is stated for splice-quotients with $\mathbb{Z}HS^3-$link, nonetheless the result is true for splice quotient singularities in general. Moreover, Theorem \ref{teo:Z=P} sheds light on the relation between $P(\mathbf{t}) \text{ and } \zeta(\mathbf{t})$ when reducing variables. Fix $v$ a vertex in $\Gamma$. Let $\zeta_v(t):=\zeta(\mathbf{t})|_{t_j=1 \text{ for all } j\neq v}$ be the topological series associated to $v$ and $P_v(t)$ the Poincaré series defined with the filtration $\mathcal{F}_v$ (Equation \ref{eq:filtration_vertex}). Therefore,
	\begin{theorem}\cite[Theorem 6.1.7]{ZPNemethi}\label{teo:Z=P_vertex}
		For $(X,0)$ a splice-quotient with $\mathbb{Z}HS^3-$link and $E_v$ any exceptional divisor of the resolution
		\begin{equation*}
			P_v(t)=P(\mathbf{t})|_{t_j=1 \text{ for all } j\neq v}.
		\end{equation*}
		Hence, $P_v(t)=\zeta_v(t)$ and the Poincaré series associated to any exceptional divisor is also topogical.
	\end{theorem}
	
	Write $P_v(t)=\sum_{i\geq 0} c_it^i$. The elements of the set
	\begin{equation*}
		S_v=\{\; i\;|\; c_i\neq0\;\}
	\end{equation*}
	
	are the integer multiplicities of functions along the exceptional divisor $E_v$; see \cite{nemethi2012cohomology}. In particular, this set has a semigroup structure and, although Poincaré series are in general analytical, Theorem \ref{teo:Z=P_vertex} shows that the semigroup $S_v$ is topological. Furthermore, as Némethi states in \cite{nemethi2012cohomology}, ``The semigroups $S_v$ are extremely hard to determine or to characterize its properties". We devote the rest of the section to prove the following result.
	
	\begin{theorem}\label{thm:semigroup_divisors}
		Let $(X,0)$ be a splice type singularity with $\mathbb{Z}HS^3-$link, $\pi:(\Tilde{X},0)\rightarrow (X,0)$ a good resolution that compatible with the splice type analytic structure, let $\Gamma$ be the dual graph of $\pi$ and $\Lambda$ the splice diagram constructed from $\Gamma$. Fix a vertex $v$ of $\Gamma$, let $K_v$ be the knot defined by $v$ in the link of the singularity and let $\Tilde{\Lambda}$ be the modified splice diagram with a unique arrow $v'$. Recall that $v'$ represents the knot $K_v$ in the splice diagram (we also consider $v'$ as a vertex of $\Tilde{\Lambda}$). 
		\begin{itemize}
			\item If $\delta_v=1$ or $\delta_v\geq3$. Let $(\Tilde{\Lambda},v')$ be the rooted splice diagram whose root is the vertex corresponding to the arrow . Then, 
			\begin{align*}
				S_v=&S(\Tilde{\Lambda}_{v'}) \text{ if }\delta_v\geq 3,\\
				S_v=&(\lceil d_w/d_1^2 \rceil+s)\mathbb{N}+S(\Tilde{\Lambda}_{v'}) \text{ if }\delta_v=1, 
			\end{align*}
			where:
			\begin{itemize}
				\item $d_w$ is the product of all the edge weights adjacent to $w$ (the nearest node to $v'$ in $\Tilde{\Lambda}$).\\
				
				\item $d_1$ is the edge weight pointing towards $v'$ at $w$.\\
				
				\item $s$ is: 
				\begin{enumerate}
					\item If $d_1=1$: the number of blow ups in the resolution between the nearest node to $v$ and $v$; $v$ included\\
					
					\item If $d_1>1$: the number of blow ups, starting from the minimal good resolution, performed in the direction of $v$ to obtain $v$.  \\
				\end{enumerate} 
			\end{itemize}

			\item If $\delta_v=2$. Let $a,b$ be the (coprime) edge weights of the adjacent edges to $v$ in $\Lambda^{max}_\Gamma$. Let $w$ be the nearest node to $v'$ in $\Tilde{\Lambda}$, $e_1$ the edge adjacent to $w$ with edge weight $a$ at $w$ and $e_2$ the edge adjacent to $w$ with edge weight $b$ at $w$. Consider the rooted splice diagrams $(\Tilde{\Lambda}_{we_1},w)$ and $(\Tilde{\Lambda}_{we_1},w)$. Then,
			\begin{equation*}
				S_v=aS((\Tilde{\Lambda}_{we_2})_w)+bS((\Tilde{\Lambda}_{we_1})_w)
			\end{equation*}
			
		\end{itemize}
	\end{theorem}
	\begin{proof}
		The following equation will be useful throughout the proof, 
		\begin{equation}\label{eq:geometric_sum}
			\frac{1}{1-t^n}=\sum_{i\in\mathbb{Z}_{\geq 0}}t^{i\cdot n}.
		\end{equation}
		
		Recall that $(X,0)$ is a splice type singularity with integral homology sphere and $\pi$ is a good resolution compatible with the splice type analytic structure. Therefore, $\Lambda$ satisfies the semigroup condition.\\
		
		First, we address the case $\delta_v=1$ or $\delta_v\geq3$. Note that $\Tilde{\Lambda}$ is $\Lambda$ but with a boundary vertex replaced with an arrowhead (if $\delta_v=1$) or with a new edge with edge weight $1$ joining a node with $v'$ (if $\delta_v\geq 3$). Hence, $\Tilde{\Lambda}$ and the rooted splice diagram $(\Tilde{\Lambda},v')$ satisfies the semigroup condition as well. Then, Corollary \ref{cor:teosemigroup_rooted} and Theorem \ref{teo:alexpoly_generating} yield
		\begin{equation}\label{eq:alex_poly_semigroup_1}
			\Delta_{K_v}(t)=(1-t)\left(\sum_{s\in S(\Tilde{\Lambda}_{v'})} t^s \right).
		\end{equation}
		
		Therefore, Theorem \ref{teo:zeta_alex_poly}, Theorem \ref{teo:Z=P_vertex} and Equation \ref{eq:alex_poly_semigroup_1} imply
		\begin{equation*}
			\frac{\sum_{s\in S(\Tilde{\Lambda}_{v'})} t^s }{1-t^{\ell_{vv}}}=P_v(t)
		\end{equation*}
		
		and Equation \ref{eq:geometric_sum} shows for $\delta_v=1 \text{ and } \delta_v\geq 3$
		\[S_v=\ell_{vv}\mathbb{N}+S(\Tilde{\Lambda}_{v'}).\]
		
		If $\delta_v\geq 3$, the linking number $\ell_{vv}$ is the product of all the edge weights adjacent to $v$. Thus, since the edge joining $v$ with $v'$ has edge weight $1$, $\ell_{vv}=\ell_{v'v}$ and is a syzygy of a distinguished sequence for $S(\Tilde{\Lambda}_{v'})$; see Section \ref{section:three}. Therefore, $\ell_{vv}\in S(\Tilde{\Lambda}_{v'}).$ If $\delta_v=1$, Lemma \ref{lemma:linking_ww} provides the result.\\
		
		Proceed with the case $\delta_v=2.$ Delete $w$ and the adjacent edge from $\Tilde{\Lambda}_{we_1}$ and from $\Tilde{\Lambda}_{we_2}$, the resultant splice diagrams are subdiagrams of $\Lambda$. Thus, the rooted splice diagrams $(\Tilde{\Lambda}_{we_1},w)$ and $(\Tilde{\Lambda}_{we_1},w)$ satisfy the semigroup condition. Consider $K_1$(resp. $K_2$) the knot defined by $w$ in (the integral homology sphere defined by) the splice diagram $\Tilde{\Lambda}_{we_1}$ (resp. $\Tilde{\Lambda}_{we_2}$). Corollary \ref{cor:teosemigroup_rooted} and Theorem \ref{teo:alexpoly_generating} yield
		
		\begin{equation}\label{eq:alex_poly_semigroup_2}
			\Delta_{K_i}(t)=(1-t)\left(\sum_{s\in S((\Tilde{\Lambda}_{we_i})_w)} t^s \right),\text{ for }i=1,2.
		\end{equation}
		
		Let $\mathcal{V}_1$ (resp. $\mathcal{V}_2$) be the set of vertices other than $w$ in $\Tilde{\Lambda}_{we_1}$ (resp. $\Tilde{\Lambda}_{we_2}$), the linking numbers in $\Tilde{\Lambda}$ of $v'$ are
		$$\begin{array}{cc}
			\ell_{v'x}=b\ell'_{wx} &\text{ if }x\in\mathcal{V}_1 \\
			\ell_{v'x}=a\ell'_{wx} &\text{ if }x\in\mathcal{V}_2\\
			\ell_{v'w}=ab,&
		\end{array}$$
		and is straightforward to show that
		\begin{equation*}
			\Delta_{K_v}=(1-t)(1-t^{ab})\frac{\Delta_{K_1}(t^b)}{(1-t^b)}\frac{\Delta_{K_2}(t^a)}{(1-t^a)}.
		\end{equation*}
		
		To conclude, Equation \ref{eq:alex_poly_semigroup_2}, Theorem \ref{teo:zeta_alex_poly} and Theorem \ref{teo:Z=P_vertex} reveal that
		\begin{equation*}
			P_v(t)=\left(\sum_{s\in S((\Tilde{\Lambda}_{we_1})_w)} t^{b\cdot s}\right)\left(\sum_{s\in S((\Tilde{\Lambda}_{we_2})_w)} t^{a\cdot s} \right),
		\end{equation*}
		which implies  
		\[S_v=aS((\Tilde{\Lambda}_{we_2})_w)+bS((\Tilde{\Lambda}_{we_1})_w).\]
	\end{proof}
	
	\begin{rem}
		Let $v$ be a vertex of $\Gamma$ with $\delta_v=1 \text{ or } \delta_v=2$. The expressions 
		\begin{align*}
			(\lceil d_w/d^2_1 \rceil+s)\mathbb{N}+S(\Tilde{\Lambda}_{v'}),\\
			aS((\Tilde{\Lambda}_{we_2})_w)+bS((\Tilde{\Lambda}_{we_1})_w)
		\end{align*}
		provide a description of the semigroup $S_v$. Note that $\gcd(\lceil d_w/d^2_1 \rceil+s,1)=\gcd(a,b)=1$. However, $\lceil d_w/d^2_1 \rceil+s$ is not necessary an element of $S(\Tilde{\Lambda}_{v'})$,  anagously  $a$ is not necessary an element of $S((\Tilde{\Lambda}_{we_1})_w)$ and $b$ is not necessary an element of $S((\Tilde{\Lambda}_{we_2})_w)$. Therefore, in general, the aforementioned expressions do not represent a gluing of semigroups.\\
		
		For instance, consider the splice type singularity $(X,0)\subset (\mathbb{C}^3,0)$ defined by the equation $\{x^2+y^3+z^7=0\}$. Let $\pi$ be the minimal good resolution of $(X,0)$ and let $\Gamma$ be its dual graph (Figure \ref{fig:not_gluing_2}) and let $\Lambda$ be its splice diagram (Figure \ref{fig:not_gluing_1}). Let $v$ be the leaf in $\Gamma$ in the direction of $v'$. Then, 
		$$\begin{array}{c}
			\ell_{vv}=\lceil 42/49 \rceil=\lceil 6/7 \rceil=1\\
			\\
			S(\Tilde{\Lambda}_{v'})=2\mathbb{N}+3\mathbb{N}
		\end{array}$$
		Hence, $1$ does not belong to $S(\Tilde{\Lambda}_{v'})$ and the expression $\mathbb{N}+S(\Tilde{\Lambda}_{v'})$ does not represent a gluing of semigroups. Note that $S_v=\mathbb{N}.$
		\begin{figure}[H]
			\centering
			\begin{subfigure}{0.4\linewidth}
				\centering
				\begin{tikzpicture}[scale=1, thick]
					
					\coordinate (C)  at (0, 0);
					
					\coordinate (TL) at (-1, 1);
					\coordinate (BL) at (-1, -1);
					
					\coordinate (R)  at (1.2, 0);
					
					\draw (TL) -- (C) node[pos=0.6, above =1pt] {$3$};
					
					\draw (BL) -- (C) node[pos=0.4, above =1pt] {$2$};
					
					\draw (C) -- (R) node[pos=0.3, below=1pt] {$7$};
					
					\fill (C)  circle (2pt) node[above =1pt]{$w$};
					\fill (TL) circle (2pt);
					\fill (BL) circle (2pt);
					\fill (R) circle (2pt) node[right=1pt]{$v'$};
					
				\end{tikzpicture}
				\caption{Splice diagram $\Lambda$ for the link of $(X,0)$}
				\label{fig:not_gluing_1}
			\end{subfigure}
			\hfill
			\begin{subfigure}{0.58\linewidth}
				\centering
				\begin{tikzpicture}[scale=0.9, thick]
					
					\coordinate (C)  at (0, 0);
					
					\coordinate (TL) at (-1, 1);
					\coordinate (BL) at (-1, -1);
					
					\coordinate (R)  at (1.2, 0);
					
					\draw (TL) -- (C) ;
					
					\draw (BL) -- (C) ;
					
					\draw (C) -- (R) ;
					
					\fill (C)  circle (2pt) node[above =1pt]{$-1$};
					\fill (TL) circle (2pt) node[above=1pt]{$-3$};
					\fill (BL) circle (2pt) node[above left=0.5pt]{$-2$};
					\fill (R) circle (2pt) node[right=1pt]{$v$} node[above=1pt]{$-7$};
					
				\end{tikzpicture}
				
				\caption{Resolution graph $\Gamma$ of the minimal good resolution.}
				\label{fig:not_gluing_2}
			\end{subfigure}
			
			\caption{\centering Splice type singularity with $\mathbb{Z}HS^3$-link defined by the equation $\{x^2+y^3+z^7=0\}$.}
		\end{figure}
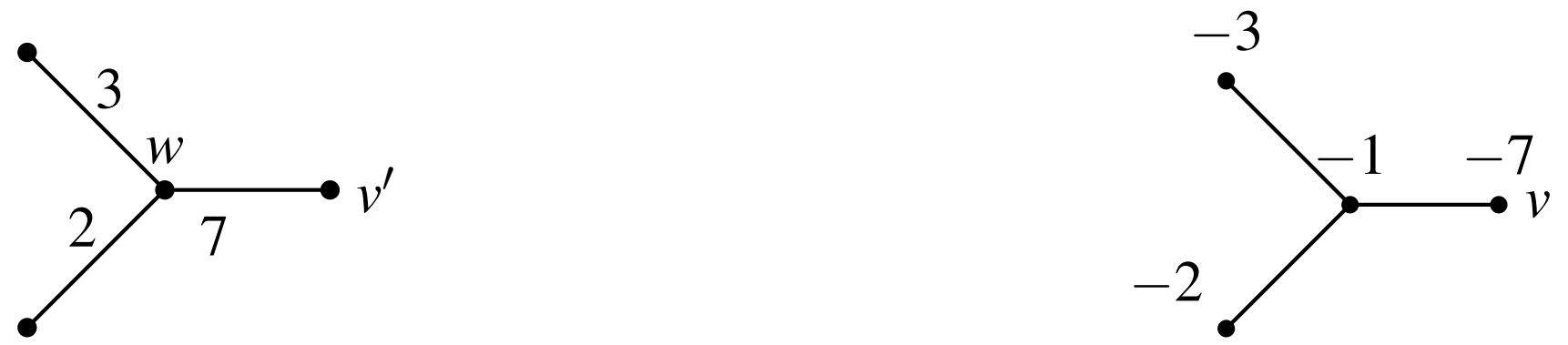
	\end{rem}

	\subsection{The Poincaré series of an end curve function}\label{subsection:42}
	
	Up to this point our attention has been focused on the filtration defined by the set of divisors from the resolution $\pi: (\Tilde{X},0)\rightarrow (X,0).$ Nevertheless, let $(C,0)\subset (X,0)$ be an embedded curve such that $\pi$ is an embedded resolution of $(C,0)$ and no exceptional divisor intersects two different components of the strict transform. Under the previous conditions, the curve $(C,0)$ defines a filtration over the local ring $\mathcal{O}_X.$ Recall that the splice quotient structure of $(X,0)$ comes equipped with a set of functions that play an essential role on the singularity, the end curve functions. Hence, is worth studying the filtration defined by the zero locus of an end curve function and if the filtration can be related to the combinatoric aspects of the singularity. Throughout the rest of the section, $(X,0)$ is a splice type singularity with integral homology sphere, unless otherwise stated.\\
	
	Let $\Gamma$ be the dual graph of the resolution and let $\Lambda$ be the splice diagram constructed from $\Gamma$. Let $z_1,\dots,z_{s}\in\mathcal{O}_X$ be the end curve functions, one for each leaf of $\Gamma$, and denote $v_i$ the leaf in $\Lambda$ associated with $z_i$. Let $C_i:\{z_i=0\}$ be the end curve defined by $z_i$ in $X$. The strict transform $\bar{C}_i$ of $C_i$ in $\Tilde{X}$ is smooth an intersects transversely in one point the exceptional divisor corresponding to a leaf of the resolution graph. The filtration defined by $C_i$ is
	\begin{equation*}
		\mathcal{F}_{C_i}(l):=\{\hspace{1mm} f\in\mathcal{O}_X \hspace{1mm}|\hspace{1mm} (\hspace{1mm} div(f\circ\pi)\hspace{1mm} ,\hspace{1mm} \bar{C}_i\hspace{1mm} )\geq l\hspace{1mm}\},\hspace{1.5mm} l\in\mathbb{Z},
	\end{equation*}
	where $(,)$ denotes intersection multiplicity.\\
	
	Moreover, $\mathcal{F}_{C_i}(l)$ can be interpreted in terms of the curve. Denote $I_{C_i}:=(z_i)$ the ideal defined by $C_i$ in $\mathcal{O}_X$, any function $f\in I_{C_i}$ satisfies $(\hspace{1mm} div(f\circ\pi)\hspace{1mm} ,\hspace{1mm} \bar{C}_i\hspace{1mm} )= \infty$. Therefore, the functions in $I_{C_i}$ belong to $\mathcal{F}_{C_i}(l)$ for all $l\in\mathbb{Z}$ and the ideal of the curve (i.e functions containing $z_i$ as a factor) is irrelevant in the computation of the dimensions
	\begin{equation*}
		\dim_{\mathbb{C}}\left( \mathcal{F}_{C_i}(l) /\mathcal{F}_{C_i}(l+1) \right) \text{ for all } l\in\mathbb{Z}.
	\end{equation*}
	Hence, the filtration depends only on the local ring $\mathcal{O}_{C_i}:=\mathcal{O}_X/I_{C_i}$ of the curve. Since $\pi$ resolves $C_i$, $\pi: \bar{C}_i\rightarrow C_i$ is the normalization of $C_i$ and defines an inclusion of rings
	\[\pi^*: \mathcal{O}_{C _i}\hookrightarrow \bar{\mathcal{O}}_{C_i}.\]
	
	Recall that $\bar{\mathcal{O}}_{C_i}$ is a discrete valuation ring. Thus, there exists an element $t\in \bar{\mathcal{O}}_{C_i}$ such that any element $\tau \in\bar{\mathcal{O}}_{C_i}\setminus\{0\}$ can be expressed as 
	\[\tau=ut^n,\]
	where $u$ is a unit in $\bar{\mathcal{O}}_{C_i}$ and $n\in \mathbb{Z}_{\geq 0}.$ Consequently, $\pi^*$ defines a valuation 
	$\nu_i:\mathcal{O}_{C_i} \rightarrow \mathbb{Z}_{\geq 0} \cup \{\infty\}$
	that sets $\nu_i(f)=\ord_{t}(\pi^*(f))$, if $f$ is non-zero, and $\nu_i(0)=\infty$. In particular, for $f\in \mathcal{O}_{C_i}$ the value $\nu_i(f)$ represents the intersection multiplicity between the curves $D:\{f=0\}$ and $C_i$. Hence, the filtration $\mathcal{F}_{C_i}$ can be expresssed in terms of the normalization, i.e,
	\begin{equation*}
		\mathcal{F}_{C_i}(l):=\{\hspace{1mm} f\in\mathcal{O}_{C_i} \hspace{1mm}|\hspace{1mm} \nu_i(f)\geq l\hspace{1mm}\}
	\end{equation*}
	
	Furthermore, the unique maximal ideal $\mathfrak{m}$ of the local ring $\mathcal{O}_X$ is generated by the set $\{z_1,\dots,z_{s}\}$ of end curve functions; see \cite[Theorem 4.1]{neumann2005complex}. Hence, any $f\in\mathcal{O}_X$ can be writen as a convergent power series in $\{z_1,\dots,z_{s}\}$. Now that the ingredients have been presented, we proceed with the computation of the Poincaré series defined by the filtration $\mathcal{F}_{C_i}(l)$.\\
	
	Fix an end curve function, for simplicity we choose $z_s$. A first approach to figure which pieces $\mathcal{F}_{C_s}(l) /\mathcal{F}_{C_s}(l+1)$ have positive dimension is computing the intersection number between $C_s$ and the curve $C_i:\{z_i=0\}$, $i=1,\dots,s-1.$ Since $z_i$ is an end curve function, 
	$$div(z_i\circ \pi)=E_i^*+\bar{C}_i,$$
	where $E_i$ is the exceptional divisor intersected by the strict transform $\bar{C}_i$. Then,
	\begin{equation*}
		(\hspace{1mm} div(z_i\circ\pi)\hspace{1mm} ,\hspace{1mm} \bar{C}_s\hspace{1mm} )= (\hspace{1mm}E_i^*+\bar{C}_i\hspace{1mm},\hspace{1mm} \bar{C}_s\hspace{1mm})= (\hspace{1mm}E_i^*\hspace{1mm},\hspace{1mm} \bar{C}_s\hspace{1mm})= \ell_{v_iv_s}.
	\end{equation*}
	
	Consecuently, $\nu_s(z_i)=\ell_{v_iv_s}$. The monomials $\prod_{i\neq s}z_i^{\alpha_i}$ ($\alpha_i\in \mathbb{Z}_{\geq0} \text{ for all $i$ })$ are elements of $\mathcal{O}_{X}$, therefore,
	\begin{equation*}
		\dim_{\mathbb{C}}\left(\mathcal{F}_{C_s}(l) /\mathcal{F}_{C_s}(l+1)\right)\geq 1, \text{ for all } l\in S(\Lambda_{v_s}).
	\end{equation*}
	Furthermore, \cite[Theorem 4.1]{neumann2005complex} ensures that the semigroup of values of the curve $C_s$ is indeed $S(\Lambda_{v_s})$, ergo the values of $\nu_s$ belong to $S(\Lambda_{v_s})$, i.e, 
	\begin{equation*}
		\dim_{\mathbb{C}}\left(\mathcal{F}_{C_s}(l) /\mathcal{F}_{C_s}(l+1)\right)\geq 1 \text{ if and only if } l\in S(\Lambda_{v_s}).
	\end{equation*}
	
	To conclude, we show that the dimension is exactly one. Let $f,g\in \mathcal{F}_{C_s}(l) /\mathcal{F}_{C_s}(l+1)$ for some $l\in S(\Lambda_{v_s})$,
	\begin{equation*}
		\pi^*(f)= \alpha\hspace{0.2mm} t^l+\sum_{i>l}\alpha_i\hspace{0.2mm} t^i, \hspace{1cm} \pi^*(g)= \beta\hspace{0.2mm} t^l+\sum_{i>l}\beta_i\hspace{0.2mm} t^i,\hspace{2.5mm}\alpha,\beta\neq 0.
	\end{equation*}
	
	Consider $h:=f-\frac{\alpha}{\beta}g\in\mathcal{O}_{C_i}.$ Then,
	\begin{equation*}
		\pi^*(h)=\pi^*(f)-\frac{\alpha}{\beta} \pi^*(g)=\sum_{i>l}\gamma_i \hspace{0.2mm}t^i. 
	\end{equation*}
	
	Therefore, $f\text{ and }g$ are in the same $\mathbb{C}$-vector subspace of $\mathcal{F}_{C_s}(l) /\mathcal{F}_{C_s}(l+1)$ and that the dimension has to be $1.$\\
	
	The argument for the remaining end curve functions can be carried out \textit{mutatis mutandis}. Hence,
	\begin{equation}\label{eq:P=S_irreducible}
		P_{C_i}(t)=\sum_{l\geq 0}  \dim_{\mathbb{C}}\left(\mathcal{F}_{C_i}(l) /\mathcal{F}_{C_i}(l+1)\right)t^l=\sum_{l\in S(\Lambda_{v_i})} t^l= \mathcal{G}_{S(\Lambda_{v_i})}(t),
	\end{equation}
	where $\mathcal{G}_{S(\Lambda_{v_i})}(t)$ denotes the generating series of $S(\Lambda_{v_i}).$ Moreover,
	
	\begin{theorem}\label{teo:P=Alexander=G_irreducible}
		Let $(X,0)$ be a splice-quotient singularity with integral homology sphere link. Let $\pi: (\Tilde{X},0)\rightarrow (X,0)$ be a good resolution compatible with the splice quotient analytic structure, let $\Gamma$ be the dual graph and let $\Lambda$ be the splice diagram obtained from $\Gamma$. Let $\{z_1,\dots,z_s\}$ be the set of end curve functions (one for each leaf of $\Gamma$), let $v_i$ be the vertex in $\Lambda$ associated with $z_i$, let $C_i:\{z_i=0\}$ be the end curve defined by $z_i$ in $X$ and let $K_i$ be the knot cut out by $C_i$ in $\Sigma.$ Then,
		\begin{equation*}
			(1-t)P_{C_i}(t)=(1-t)\mathcal{G}_{S(\Lambda_{v_i})}(t)=\Delta_{K_i}(t),
		\end{equation*}
		
	\end{theorem}
	\begin{proof}
		Equation \ref{eq:P=S_irreducible} shows $P_{C_i}(t)=\mathcal{G}_{S(\Lambda_{v_i})}(t).$ Note that $\Lambda$ satisfies the semigroup condition, therefore, the rooted splice diagram $(\Lambda,v_i)$ satisfies the condition as well. The relation 
		$$(1-t)\mathcal{G}_{S(\Lambda_{v_i})}(t)=\Delta_{K_i}(t)$$ 
		is a direct consequence of Corollary \ref{cor:teosemigroup_rooted} and Theorem \ref{teo:alexpoly_generating}.
	\end{proof}
	
	\subsection{The non-irreducible case}
	We follow the notation from Section \ref{subsection:42}. Consider a subset of the end curve functions, say $\{ z_{1},\dots, z_k\}$ with $1<k\leq s$). The function $z_{1}\cdots z_{k}$ defines a reduced curve germ $(C,0)\subset (X,0)$ embedded in the surface, with $k$ irreducible components $C_1, \ldots , C_k$, where $C_i$ is the end curve defined by $z_i$. In particular, $C$ has $k$ irreducible components $\{C_1,\dots,C_k\}$ and each $C_{i}$ is the end curve defined by the end curve function $z_i$. The strict transform $\bar{C}$ of $C$ has irreducible components $\bar{C}_1,\dots,\bar{C}_k$, and each of them intersects transversely a unique exceptional divisor, which is a leaf of $\Gamma$.
	
	Since different end curve functions correspond to different leaves, the resolution $\pi: (\Tilde{X},0)\rightarrow (X,0)$ is an embedded good resolution of $(C,0)$ in which no exceptional divisor is intersected by two different components of $\bar{C}$. Hence, the curve defines a filtration on the local ring of the surface:
	\begin{equation*}
		\mathcal{F}_{C}(l):=\left\{f\in \mathcal{O}_{X,0}\hspace{1.5mm} | \hspace{1.5mm} \left( div(f\circ\pi), \bar{C}_i\right)\geq l_i,\hspace{2.5mm} i=1,\dots,k\; \right\} ,\hspace{1.5mm} l:=(l_1,\dots,l_{k})\in \mathbb{Z}^{k}.
	\end{equation*}
	
	As in the irreducible case,  $I_C:=(z_{1}\cdots z_{k})$, the ideal of the curve in $\mathcal{O}_{X,0}$, is contained in $\mathcal{F}_C(l)$, for every $l\in \mathbb{Z}^{k}$. Therefore, $I_C$ does not affect the dimensions of the quotients
	\begin{equation*}
		\mathcal{F}_{C}(l) /\mathcal{F}_{C}(l+(1,\dots,1))  ,\; l\in \mathbb{Z}^{k}
	\end{equation*}
	and the filtration can be considered on the local ring $\mathcal{O}_C:=\mathcal{O}_{X}/I_C$ of the curve. Moreover, the filtration can be recovered from the curve alone. The restriction of $\pi$ to $\bar{C}_i$ is actually the normalization of the curve $C_i$; fix a local coordinate $t_i\in\bar{\mathcal{O}}_{C_i}$ and let $\nu_i:\mathcal{O}_{C_i} \rightarrow \mathbb{Z}_{\geq 0} \cup \{\infty\} $ be the valuation it defines.

	Consequently, the filtration can be defined as
	\begin{equation*}
		\mathcal{F}_C(l):=\left\{f\in \mathcal{O}_{C}\hspace{1.5mm} | \hspace{1.5mm} \nu_i(f)\geq l_i,\hspace{2.5mm} i=1,\dots,k\;\right\} ,\hspace{1.5mm} l\in \mathbb{Z}^{k}.
	\end{equation*}
	
	Let $t_i$ be the variable associated to $C_i$ in $\mathbb{Z}[t_1^{\pm 1},\dots, t_k^{\pm 1}]$. Denote $\underline{\mathbf{1}}:=(1,\dots,1)\in\mathbb{Z}^{k}$  and set $\underline{\mathbf{t}}^l=\prod_{i=1}^k t_i^{l_i}$. Define the series
	\begin{equation}\label{eq:L_series_non-irreducible_curve}
		L_C(\mathbf{t})=\sum_{l\in \mathbb{Z}^{k}}\dim_{\mathbb{C}}(\mathcal{F}_C(l)/\mathcal{F}_C(l+\underline{\mathbf{1}}))\underline{\mathbf{t}}^l.
	\end{equation}
	Then, the Poincaré series of the filtration $\mathcal{F}_C$ is:
	\begin{equation}\label{eq:Poincare_series_non-irreducible_curve}
		P_C(\mathbf{t})=\frac{L_C(\mathbf{t})\prod_{i=1}^k(t_i-1)}{(\prod_{i=1}^k t_i-1)}.
	\end{equation}

	We will provide a generalization of Theorem \ref{teo:P=Alexander=G_irreducible} for the (reduced) curve $C:\{z_1\cdots z_k=0\}$, a curve such that each branch is an end curve. To achive the result, the plan is to describe the curve valuation as a limit of divisorial valuations following the procedure exposed in \cite{Nemethi-Veys_filtrations}. The procedure goes as follows:\\
	
	Let $\mathcal{V}_C:=\{1,\dots,k\}$ be the set of vertices in $\Gamma$ whose corresponding exceptional divisor is intersected by the strict transform of the curve. Each component $\bar{C}_i$ of $\bar{C}$ intersects a unique exceptional divisor $E_i$ in a point $p_i , \hspace{2.5mm} i\in \mathcal{V}_C$. Moreover, different components intersect different exceptional divisors. Hence, blowing up $\Tilde{X}$ with centers $\{p_i\}$ produces a good resolution $\pi_1:(\Tilde{X}_1,0)\rightarrow (X,0)$, with new exceptional divisors $\{E_{i,1}\}_{i\in\mathcal{V}_C}$ and strict transforms $\bar{C}$ and $\bar{C}_i$; we will maintain the same notation for the strict transform and its irreducible components. The strict transform $\bar{C}_i\subset \Tilde{X}_1$  intersects transversely the exceptional divisor $E_{i, 1}$ in a point. Thus, blowing up $\Tilde{X}_1$ with centers the intersection points $\{E_{i, 1}\cap \bar{C}_i\}$ defines a good resolution $\pi_2:(\Tilde{X}_2,0)\rightarrow (X,0)$, with new exceptional divisors $\{E_{i, 2}\}_{i\in\mathcal{V}_C}$ and strict transforms $\bar{C}$ and $\bar{C}_i$. After $m$ blow ups, $m\in\mathbb{Z}_{>0}$, one gets a manifold $\Tilde{X}_m$ with new exceptional divisors $\{E_{i, m}\}_{i\in\mathcal{V}_C}$ and a good resolution $\pi_m:(\Tilde{X},0)\rightarrow (X,0)$. Since $\pi$ was an embedded resolution for $(C,0)$ and a good resolution for $(X,0)$, each resolution $\pi_i$ satisfies both properties as well.\\
	
	For each interger $m\in\mathbb{Z}_{>0}$, one can define a filtration on $\mathcal{O}_{X,0}$ associated to the set $\{E_{i,m}\}_{i\in\mathcal{V}_C}$:
	
	\begin{equation*}
		\mathcal{F}_{\mathcal{V}_C, m}(l):=\{f\in \mathcal{O}_{X,0}\hspace{1.5mm} | \hspace{1.5mm} div(f\circ\pi_m)\geq \sum_{i\in\mathcal{V}_C} l_iE_{i,m}\},\hspace{1.5mm} l\in \mathbb{Z}^{k}.
	\end{equation*}
	
	Take a variable $t_i$ for each $E_{i,m},\; i\in\mathcal{V}_C$. The series $L_{\mathcal{V}_C, m}(\mathbf{t})$ and $P_{\mathcal{V}_C, m}(\mathbf{t})$ are defined as in Eq. \ref{eq:L_series_non-irreducible_curve} and \ref{eq:Poincare_series_non-irreducible_curve} respectively, but with the filtration $\mathcal{F}_{\mathcal{V}_C, m}$ instead. On this setting, 
	Némethi and Veys in \cite{Nemethi-Veys_filtrations} proof
	\begin{theorem}\cite[Theorem 2.3.1 and Corollary 3.1.2]{Nemethi-Veys_filtrations}\label{teo:NV-filtrations}
		
		\begin{enumerate}
			
			\item For each $l\in\mathbb{Z}^{k}$, the limit filtration $\lim_m \mathcal{F}_{\mathcal{V}_C, m}$ given by
			\begin{equation*}
				(\lim_m \mathcal{F}_{\mathcal{V}_C, m})(l):=\lim_{m\rightarrow\infty}(\mathcal{F}_{\mathcal{V}_C, m}(l))
			\end{equation*}
			is a well defined filtration on $\mathcal{O}_{X,0}$.\\
			
			\item We have the equality of filtrations 
			\begin{equation*}
				\lim_m \mathcal{F}_{\mathcal{V}_C, m}=\mathcal{F}_C.
			\end{equation*}
			
			\item We have that 
			\begin{equation*}
				\lim_m P_{\mathcal{V}_C, m}(\mathbf{t})=P_C(\mathbf{t}).
			\end{equation*}
		\end{enumerate}
		
	\end{theorem}
	
	In the sequel we write $\pi_0\text{ for the resolution }\pi,$ and write $\{E_{i,0}\}_{i\in\mathcal{V}_C}$ for the set exceptional divisors of the resolution $\pi$ intersected by the strict transform of $C.$ For each $m\in\mathbb{Z}_{\geq 0}$, let $\pi_m$ be the resolution obtained after $m$ steps of the aforementioned procedure and let $\Gamma_m$ be its dual graph. Let $L_m$ be the lattice of integral cycles associated with $\pi_m$, generated by the classes of the exceptional divisors $\{E_v\}_{v\text{ vertex in } \Gamma_m}$, and let $L_m^*$ be its dual, generated by the negative duals $\{E_v^*\}_{v\text{ vertex in }\Gamma_m}.$ Define for $\pi_m$ the topological function $\zeta_m(\mathbf{t})$ as in Eq. \ref{eq:zeta_griega_function}.
	\begin{defin}
		For each resolution resolution $\pi_m$, let $\{E_{i,m}\}_{i\in\mathcal{V}_C}$ be the set of set of exceptional divisors intersecting the strict transform of the curve and let $(\mathcal{V}_C,m)$ be the corresponding vertices in $\Gamma_m$. Then, define the zeta function restricted to the exceptional divisors $\{E_{i,m}\}_{i\in\mathcal{V}_C}$ as:
		\[\zeta_{{\mathcal{V}_C},m}(\mathbf{t}):=\zeta_m(\mathbf{t})|_{t_j=1 \text{ for all } j\notin (\mathcal{V}_C,m)}.
		\]
		
	\end{defin}
	Let $\Lambda^{\text{max}}_{\Gamma_m}$ (resp. $\Lambda_m$) be the maximal splice diagram (resp. splice diagram) obtained from $\Gamma_m$ and let $\Sigma$ be the integral homology sphere associated with $\Lambda$. On the step $m$ of the aforementioned procedure we blow up points in (the smooth part of) the exceptional divisors $\{E_{i,m}\}_{i\in \mathcal{V}_C}$. Consequently, the blow ups do not add nodes to the dual graph and the splice diagram $\Lambda_m$ is not altered, i.e, $\Lambda_m=\Lambda_n$ for all $m,n\in\mathbb{Z}_{\geq 0}.$ Let $\mathcal{E}$ (resp. $\mathcal{N}$) be the set of leaves (resp. nodes) of $\Lambda.$ Divide $\mathcal{E}$ into two subsets: the subset of leaves whose corresponding knot is cut out in $\Sigma$ by a component of $C$ (denoted $\mathcal{V}_C$) and  the remaining leaves (denoted $\overline{\mathcal{E}}$). Recall that vertices in $\Gamma_m$ of valence greater than two and nodes in $\Lambda$ are in one to one correspondence. Therefore, we use $\mathcal{N}$ to denote the vertices of valence greater than two in $\Gamma_m$ and we write $\overline{\mathcal{E}}$ for the leaves in $\Gamma_m$ whose divisors do not intersect $\bar{C}$. Then, in $pr_{\mathcal{V}_C,m}(L_m)$ we have:
	
	\begin{alignat}{2}\label{eq:dual_multivariable}
		E_v^*&= \sum_{i\in\mathcal{V}_C} \ell_{vi} E_{i,m}, & \qquad &v\in\mathcal{N} \; \text{ or } \; v\in \overline{\mathcal{E}},\\
		E_{i,m}^*&=\sum_{j\in\mathcal{V}_C, j\neq i} \ell_{vj} E_{j,m} + \ell_{vv,m}E_{i,m}, & \qquad & v\in\mathcal{V}_C \nonumber 
	\end{alignat}
	where $\ell_{vi}$ is the linking number in $\Lambda$ and $\ell_{vv,m}$ is the linking number in $\Lambda_{\Gamma_m}^{\text{ max }}$ of the leaf $v$ (the leaf corresponding to $E_{v,m}$). Furthermore, as computed in Section \ref{ssec:splice-quotient}, 
	\begin{equation*}
		\ell_{vv,m}=\ell_{vv}+m,
	\end{equation*}
	where $\ell_{vv}$ is the linking number in $\Lambda_{\Gamma_0}^{\text{ max }}$ of the leaf $v$ that corresponds to the exceptional divisor $E_{i,0}$.\\ 
	
	Recall that $(X,0)$ is a splice quotient singularity, $\pi:(\Tilde{X},0)\rightarrow (X,0)$ is a good resolution compatible with the splice quotient structure and $C$ is a curve whose each component is an end curve. Therefore, the set $\{E_i\}_{i\in\mathcal{V}_C}$ of exceptional divisors intersected by the strict transform $\bar{C}$ consists of leafs of the dual graph $\Gamma$. Let $\pi_1: (\Tilde{X},0)\rightarrow (X,0)$ be the good resolution obtained by blowing up $\Tilde{X}$ along the intersection points $\{E_i\cap \bar{C}_i\}_{i\in \mathcal{V}_C}$ and let $\Gamma_1$ be its dual graph. The set of exceptional divisors $\{E_{i,1}\}_{i\in\mathcal{V}_C}$ corresponds leaves in $\Gamma_1.$ We claim that the resolution $\pi_1$ is compatible with the splice quotient structure of $(X,0)$ and $\{z_1,\dots, z_s\}$ act as end curve functions for $\pi_1$ as well.\\
	
	Fix the function $z_i$, which is an end curve function for $\pi$. If $C_i:\{z_i=0\}$ is a component of $C$, the strict transform of $C_i$ under $\pi_1$ is smooth and intersects transversely the exceptional divisor $E_{i,1}.$ If $C_i:\{z_i=0\}$ is not a component of $C$, the strict transform of $C_i$ under $\pi_1$ is smooth and intersects transversely an exceptional divisor represented by a leaf in $\Gamma_1$ (the strict transform of $C_i$ under $\pi_1$ intersects the same exceptional divisor as the strict transform of $C_i$ under $\pi$). The leaves of $\Gamma_1$ consist of the leaves corresponding to the exceptional divisors $\{E_{i,1}\}_{i\in\mathcal{V}_C}$ and the leaves which are leaves in $\Gamma$ as well. Consecuently, there exists an end curve function for each of the leaves of $\Gamma_1$ and the End Curve Theorem (see  \cite{neumann2005complex}, \cite{neumann2010end}) of Neumann and Wahl implies that the resolution $\pi_1$ is compatible with the splice-quotient structure. In particular, the splice quotient structure defined on $\pi_1$ coincides with the pull-back of the splice quotient structure defined on $\pi$.\\
	
	Therefore, for each $m\in\mathbb{Z}_{\geq 0}$, the resolution $\pi_m$ is compatible with the splice quotient structure of the singularity. Therefore, the resolution $\pi_m$  satisfies Theorem \ref{teo:Z=P}. Note that the filtration $\mathcal{F}_{\mathcal{V}_C,m}$ is the divisorial filtration defined by $\pi_m$ on $\mathcal{O}_X$ restricted to the set of exceptional divisors $\{E_{i,m}\}_{i\in\mathcal{V}_C}$. Then, \cite[Theorem 6.1.7]{ZPNemethi} implies that 
	\begin{equation}\label{eq:Z=P_filtration_curve}
		P_{\mathcal{V}_C,m}(\mathbf{t})=\zeta_{\mathcal{V}_C,m}(\mathbf{t}).
	\end{equation}
	For $i\in\mathcal{V}_C$ write the coefficients of $E_{i,m}^*$ in $pr_{\mathcal{V}^C}(L_m)$ as a vector,
	\[\underline{\nu}_{\;i,m}:=(\ell_{1i},\dots,\ell_{ii,m},\dots,\ell_{ki}),\] 
	\begin{theorem}\label{teo:P=alexpoly_curve_divisorial}
		Let $(X,0)$ be a splice-quotient singularity with integral homology sphere link $\Sigma$. Let $\pi_0: (\Tilde{X},0)\rightarrow (X,0)$ be a good resolution compatible with the splice quotient analytic structure. Let $\{z_1,\dots,z_s\}$ be the set of end curve functions from $\pi_0$, define $C:\{z_1\cdots z_k=0\}$ and let $K$ be the link cut out by $C$ in $\Sigma.$ Then, for each $m\in\mathbb{Z}_{\geq 0}$ 
		\begin{equation*}
			\prod_{i\in\mathcal{V}_C}(\mathbf{t}^{\underline{\nu}_{\;i,m}}-1)P_{{\mathcal{V}^C},m}(\mathbf{t})=\Delta_L(\mathbf{t}).
		\end{equation*}
	\end{theorem}
	\begin{proof}
		First, Equation \ref{eq:dual_multivariable} shows that
		\[\prod_{i\in\mathcal{V}_C}(\mathbf{t}^{\underline{\nu}_{\;i,m}}-1)\zeta_{{\mathcal{V}^C},m}(\mathbf{t})=\Delta_L(\mathbf{t}).\]
		Then, Equation \ref{eq:Z=P_filtration_curve} concludes the proof. 
	\end{proof}
	\medskip
	
	\begin{corollary}\label{coro:relation_poincaré_divisorial}
		For each $i\in\{1,\dots,k\}$ let $\underline{1}_i$ be the element of $\mathbb{Z}^k$ with all entries zero but a $1$ on the $i$-th coordinate. Then, under the hypothesis from Theorem \ref{teo:P=alexpoly_curve_divisorial},
		\begin{equation*}
			\prod_{i\in\mathcal{V}_C}(\mathbf{t}^{\underline{\nu}_{\;i,m}}-1)P_{{\mathcal{V}^C},m}(\mathbf{t})=\prod_{i\in\mathcal{V}_C}(\mathbf{t}^{\underline{\nu}_{\;i,m}+\underline{1}_i}-1)P_{{\mathcal{V}^C},m+1}(\mathbf{t})
		\end{equation*}
	\end{corollary}
	\medskip

	\begin{corollary}\label{coro:multivariable_P=pol.Alex}
		Let $(X,0)$ be a splice quotient singularity with $\mathbb{Z}HS^3-$link $\Sigma$ and let $\pi: (\Tilde{X},0)\rightarrow (X,0)$ be a good resolution compatible with the splice-quotient structure. Let $\{z_1,\dots, z_s\}$ be the set end curve functions (one for each leaf of the resolution graph) for $\pi$ and let $(C,0)\subset (X,0)$ be the (reduced) curve defined by the function $f:= z_{1}\cdots z_{k}$, $1<k\leq s$. Let $K$ be the link cutted out by $C$ in $\Sigma$, let $\Delta_K(t)$ be the Alexander polynomial of $K$ and let $P_{C}(t)$ be the Poincaré series defined by the filtration $\mathcal{F}_{C}$. Then,
		\begin{equation*}
			P_C(\mathbf{t})=\Delta_K(\mathbf{t})
		\end{equation*}
	\end{corollary}
	\begin{proof}
		Corollary \ref{coro:relation_poincaré_divisorial} shows that the Poincaré series $P_{{\mathcal{V}_C},m}(\mathbf{t})$ can be expressed in terms of $P_{{\mathcal{V}_C},0}(\mathbf{t})$ and the values $\{\underline{\nu}_{\;i,0}\}_{i\in\mathcal{V}_C}$.
		
		\begin{equation*}
			P_{{\mathcal{V}_C},m}(\mathbf{t})=P_{{\mathcal{V}_C},0}(\mathbf{t})\prod_{n=0}^{m-1} \left(\frac{\prod_{i\in\mathcal{V}_C}(\mathbf{t}^{\underline{\nu}_{i,0}+n\underline{1}_i}-1)}{\prod_{i\in\mathcal{V}^C}(\mathbf{t}^{\underline{\nu}_{i,0}+(n+1)\underline{1}_i}-1)}\right)
		\end{equation*}
		
		Therefore, applying Theorem \ref{teo:NV-filtrations} yields
		
		\begin{equation*}\label{eq:P_C=P_V0-multiple}
			P_C(\mathbf{t})=\lim_m P_{\mathcal{V}_C,m}(\mathbf{t})=P_{\mathcal{V}_C,0}(\mathbf{t})\prod_{n=0}^{\infty} \left(\frac{\prod_{i\in\mathcal{V}_C}(\mathbf{t}^{\underline{\nu}_{i,0}+n\underline{1}_i}-1)}{\prod_{i\in\mathcal{V}_C}(\mathbf{t}^{\underline{\nu}_{i,0}+(n+1)\underline{1}_i}-1)}\right).
		\end{equation*}
		
		Now, write the previous equation in a suitable manner,
		\begin{alignat*}{1}\label{eq:P_V0-multiple}
			P_C(t)=P_{{\mathcal{V}_C},0}(\mathbf{t})&\prod_{i\in\mathcal{V}_C}(\mathbf{t}^{\underline{\nu}_{i,0}}-1)\prod_{n=1}^\infty\left(\frac{\prod_{i\in\mathcal{V}_C}(\mathbf{t}^{\underline{\nu}_{i,0}+n\underline{1}_i}-1)}{\prod_{i\in\mathcal{V}_C}(\mathbf{t}^{\underline{\nu}_{i,0}+n\underline{1}_i}-1)}\right)=P_{{\mathcal{V}_C},0}(\mathbf{t})\prod_{i\in\mathcal{V}_C}(\mathbf{t}^{\underline{\nu}_{i,0}}-1)=\Delta_K(\mathbf{t}).
		\end{alignat*}
	\end{proof}
	\begin{rem}\label{rem:plane_curves}
		Let $(C,0)\subset(\mathbb{C}^2,0)$ be a reduced plane curve germ with branches
		$C_i=\{f_i=0\}$, $i=1,\dots,k$. The smooth germ $(\mathbb{C}^2,0)$ is a splice
		quotient singularity with $\mathbb{Z}HS^3$ link $S^3$. Take an embedded good
		resolution of $C$ so that different branches meet different
		leaves. Every leaf admits an
		end curve function, and the $f_i$ are end curve functions for the leaves met by
		the $C_i$. Hence Theorem~\ref{teo:P=Alexander=G_irreducible} (for $k=1$) and
		Corollary~\ref{coro:multivariable_P=pol.Alex} (for $k\geq 2$) apply to $C$, and
		recover the result of Campillo, Delgado and Gusein-Zade \cite{CDGduke}.
	\end{rem}
	
	
	\medskip
	\noindent\textbf{Tool Use}
	The mathematical ideas, results and proofs presented in this paper were conceived and developed by the authors. Claude and Google Gemini were used for language editing, including grammar, syntax, and create figures. The authors take full responsibility for the mathematical correctness and content of the paper.

	\bibliographystyle{plainurl}
	\bibliography{biblioknots}
\end{document}